%% file: CP_LanglandsIII.tex
\documentclass[dvips,12pt]{article}
\usepackage{amsmath,amssymb,amsfonts,amscd,amsthm,pstricks,pst-node,epsfig
,soul
}
\usepackage[latin1]{inputenc}
\input{cmds}

\def\prod{\mathop{\Pi}\limits}
\def\sum{\mathop{\Sigma}\limits}
\def\bigoplus{\mathop{\oplus}\limits}

\def\vec{\overrightarrow}

\def\Sm{\operatorname{Sm}}
\def\SP{\operatorname{SP}}
\def\CY{\operatorname{CY}}
\def\CORR{\operatorname{CORR}}
\def\Corr{\operatorname{Corr}}
\def\DM{\operatorname{DM}}
\def\AdFRepsp{\operatorname{AdFRepsp}}
\def\wt{\widetilde}
\def\SV{Suslin-Voevod\-sky}
\def\fl{f_\ell \cdot\ell}
\def\fnl{f_n \cdot\ell}
\def\nfl{n_\ell [f_\ell \cdot\ell]}
\def\2nfl{2n_\ell [2f_\ell \cdot\ell]}
\def\Nfl{n [f_\ell \cdot\ell]}

\def\Abstract{\small
\begin{center}
{\bf Abstract\vspace{-.5em}\vspace{0pt}}
\end{center}
\quotation}

\newcommand{\myendofpart}{\ifodd\value{page}\vfill\eject
\thispagestyle{empty}
\fi}

\def\Bi{\begin{itemize}}

\def\Aa{{\mathbb{A}\,}}
\def\SS{{\mathbb{S}\,}}
\def\zit{{\mathbb{Z}\,}}
\def\bZ{{\mbox{\st{Z}}\,}}

\def\cit{{\mathbb{C}\,}}
\def\CC{{\mathbb{C}\,}}
\def\rit{{\mathbb{R}\,}}

\def\HH{{\mathbb{H}\,}}

\def\labelenumi{\alph{enumi})}
\def\theenumi{\alph{enumi}}
\def\labelenumii{\arabic{enumii}.}
\def\theenumii{\arabic{enumii}}

\def\Aut{\operatorname{Aut}}
\def\Hom{\operatorname{Hom}}
\def\hol{\operatorname{hol}}
\def\cusp{\operatorname{cusp}}
\def\Tan{\operatorname{Tan}}
\def\TAN{\operatorname{TAN}}
\def\Lie{\operatorname{Lie}}
\def\ppr{\operatorname{pr}}

\def\Ind{\operatorname{Ind}}

\def\Gal{\operatorname{Gal}}

\def\Rep{\operatorname{Rep}}
\def\IrrRep{\operatorname{IrrRep}}
\def\RedRep{\operatorname{RedRep}}
\def\Redcusp{\operatorname{Redcusp}}
\def\Irrcusp{\operatorname{Irrcusp}}
\def\RedELLIP{\operatorname{RedELLIP}}
\def\Repsp{\operatorname{Repsp}}

\def\resp#1{(resp. #1)}
\def\rresp#1{\qquad \mbox{(resp.} \quad #1\ )}

\def\ELLIP{\operatorname{ELLIP}}

\def\RL{_{R\times L}}

\def\FRepsp{\operatorname{FRepsp}}
\def\FREPSP{\operatorname{FREPSP}}
\def\GL{\operatorname{GL}}

\def\ELLIP{\operatorname{ELLIP}}
\def\Irr{\operatorname{Irr}}

\def\TAN{\operatorname{TAN}}
\def\EIS{\operatorname{EIS}}

\def\lr{left (resp. right) }
\def\rl{right (resp. left) }

\def\bbf{\boldmath\bf}

\def\To{\begin{CD} @>>>\end{CD}}

\begin{document}

\renewcommand{\qedsymbol}{{\small{$\blacksquare$}}}
\pagestyle{myheadings}
\parindent=0pt
\def\thepage{\arabic{page}}
\markright{}

\include{LanglandsIII_00}

\tableofcontents
\eject
\setcounter{page}{1}
\include{LanglandsIII_0}
\include{LanglandsIII_1}
\include{LanglandsIII_2}
\include{LanglandsIII_3}
\include{LanglandsIII_4}
\include{LanglandsIII_5}

\include{LanglandsIII_Bib}
\myendofpart

\markright{Contents}
\def\thepage{\roman{page}}
\setcounter{page}{1}


\end{document}

%% file: cmds.tex
\def\simarrow{\mathrel{\raise -0.5mm\hbox{$\sim$}}\hspace{-1.8mm}{\rightarrow} } 

\def\bsimarrow{\leftarrow\hspace{-0.7mm}\mathrel{\raise -0.5mm\hbox{$\backsim$}} }

\def\bt{\begin{tabular}}
\def\te{\end{tabular}}

\def\Abstract{\small
\begin{center}
{\bf Abstract\vspace{-.5em}\vspace{0pt}}
\end{center}
\quotation
}

\def\lettrine#1#2#3{\noindent\hangindent#1\hangafter-#2
\hskip-#1\smash{\hbox to #1{#3\hfill}}\ignorespaces}

\newcommand{\To}[1]{\mathop{\to}\limits_{#1}}

\def\BM{\begin{pmatrix}}
\def\EM{\end{pmatrix}}

\def\pt{\mathrel{\scriptstyle \bullet}}

\def\d=f{\buildrel\hbox{\scriptsize d\'{e}f}\over \Longleftrightarrow}

\def\cit{\text{\it I\hskip -6ptC\/}}

\def\square{\hfill\hbox{\vrule height .9ex width .8ex depth -.1ex}}
\def\rit{\text{\it I\hskip -2pt  R}}

\def\zit{\text{\it Z\hskip -4pt  Z}}

\def\nit{\text{\it I\hskip -2pt  N}}

\def\rl {\rit^{\hskip 1pt\ell}}

\def\Bd{{\text B}}

\def\Ds{{\cal D}}
\def\Ed{{\text E}}

\def\Fs{{\cal F}}

\def\Hs{{\cal H}}

\def\Ls{{\cal L}}

\def\Os{{\cal O}}

\def\be{\begin{equation}}
\def\ee{\end{equation}}
\def\beqn{\begin{eqnarray}}
\def\eeqn{\end{eqnarray}}
\def\nobeqn{\begin{eqnarray*}}
\def\noeeqn{\end{eqnarray*}}
\def\ba{\left(\begin{array}}
\def\ea{\end{array} \right) }

\def\bpr{\paragraph{Proof.}}

\def\epr{\square\vskip 6pt}

\def\u{\underline}
\def\o{\overline}
\def\and{\; \mbox{and} \;}

\newcommand{\half}{\frac{1}{2}}

\def\hfl#1#2{\smash{\mathop{\hbox to 12mm{\rightarrowfill}}
\limits^{\scriptstyle #1}_{\scriptstyle #2}}}

\def\mod{\mathop{\rm mod}\nolimits}

\def\Be{\begin{enumerate}}
\def\Ee{\end{enumerate}}

\def\Bena{\begin{enumerate}
\def\labelenumi{\theenumi)}
\def\theenumi{\arabic{enumi}}
\def\labelenumii{\theenumii)}
\def\theenumii{\alph{enumii}}}

\def\Bean{\begin{enumerate}
\def\labelenumii{\theenumii)}
\def\theenumii{\arabic{enumii}}
\def\labelenumi{\theenumi)}
\def\theenumi{\alph{enumi}}}

\def\Bero{\begin{enumerate}
\def\labelenumii{\theenumii)}
\def\theenumii{\arabic{enumii}}
\def\labelenumi{(\theenumi)}
\def\theenumi{\roman{enumi}}}

\def\BeRo{\begin{enumerate}
\def\labelenumii{\theenumii)}
\def\theenumii{\arabic{enumii}}
\def\labelenumi{(\theenumi)}
\def\theenumi{\Roman{enumi}}}

\def\Bi{\vskip 11pt\begin{itemize}\itemsep=18pt}
\def\Ei{\end{itemize}\vskip 11pt}

\def\Bd{\begin{description}}
\def\Ed{\end{description}}

\def\R{\right}
\def\L{\left}

%% file: LanglandsIII_00.tex
{{\pagestyle{empty}

\null
 \vfill
 
\begin{center} 
{\LARGE \bbf $n$-dimensional geometric-shifted global bilinear correspondences of Langlands on mixed motives ---  III\par}%
   \vskip 3em
   {\large C. Pierre
}

\end{center}

\vfill


\vfill

{\Abstract{\noindent 
This third paper, devoted to global correspondences of Langlands, bears more particularly on geometric-shifted bilinear correspondences on mixed (bi)motives generated under the action of the products, right by left,  of differential elliptic operators.

The mathematical frame, underlying these correspondences, deals with the categories of the Suslin-Voevodsky mixed (bi)motives and of the Chow mixed (bi)motives which are both in one-to-one correspondence with the functional representation spaces of the shifted algebraic bilinear semigroups.

A bilinear holomorphic and supercuspidal spectral representation of an elliptic bioperator is then developed.

}} \vfill\eject

}}


%% file: LanglandsIII_0.tex
\section*{Introduction}

\addcontentsline{toc}{section}{Introduction}

\addtocontents{toc}{\protect\thispagestyle{empty}}

This paper constitutes the third part of the $n$-dimensional global correspondences of Langlands \cite{Pie3}, \cite{Pie4} and is particularly devoted to the study of the Langlands correspondences on mixed (bi)motives.

When the $n$-dimensional global correspondences of Langlands \cite{Cara} bear on pure bimotives, the related geometric-shifted correspondences deal with mixed bimotives which are assumed to be generated under the action of the products, right by left, of differential (elliptic) operators.

This leads us, more particularly, to:
\Bean
\item work in the frame of the {\bbf categories of the Suslin-Voevodsky (mixed) (bi)mo\-tives and of the Chow (mixed) (bi)motives\/} which are both in one-to-one correspondence with the functional representation spaces of the (shifted) algebraic bilinear semigroups, as introduced in \cite{Pie3}.

\item envisage a {\bbf bilinear version for the index theorem\/}.

\item develop a {\bbf bilinear holomorphic and (super)cuspidal spectral representation of an elliptic bioperator\/}.
\Ee
\vskip 11pt

A first step then consists in introducing {\bbf triangulated categories of mixed (bi)motives\/} which are built on corresponding pure (bi)motives being in one-to-one correspondence with the {\bbf functional representation spaces of algebraic bilinear semigroups\/} recalled in chapter 1 and hereafter.

\Bi
\item Let $F_\omega $ \resp{$F_{\o\omega }$} denote the set of \lr {\bbf pseudo-ramified complex completions\/} $F_{\omega _{j,m_j}}$
\resp{$F_{\o\omega _{j,m_j}}$} 
corresponding to transcendental extensions
restricted to the upper \resp{lower} half space and 
being in one-to-one correspondence with the corresponding complex splitting subsemifields
$\widetilde F_{\omega _{j,m_j}}$
\resp{$\widetilde F_{\o\omega _{j,m_j}}$}
characterized by Galois extension degrees given by integers modulo $N$.

Similarly, let 
 $F^+_v $ \resp{$F^+_{\o v }$} be the set of \lr {\bbf real pseudo-ramified  completions\/} 
$F^+_{v _{j_\delta ,m_{j_\delta }}}$
\resp{$F^+_{\o v _{j_\delta ,m_{j_\delta }}}$}
in one-to-one correspondence with
 the corresponding  real splitting subsemifields $\wt F^+_{v_{j_\delta ,m_{j_\delta }}}$
\resp{$\wt F^+_{\o v_{j_\delta ,m_{j_\delta }}}$}
of the (finite) extension semifield 
$\wt F^+_L$ \resp{$\wt F^+_R$}
of a number field $k$ of characteristic zero.

The set $F^+_v$ \resp{$F^+_{\o v}$} of \lr real pseudoramified completions covers the corresponding set
$F_\omega $ \resp{$F_{\o\omega }$} of complex completions.

\item {\bbf The bilinear algebraic semigroup\/} of matrices $\GL_n(\wt F_{\o\omega} \times \wt F_\omega )\equiv T^t_n(\wt F_{\o\omega })\times T_n(\wt F_\omega )$ has its representation space in the   bilinear affine semigroup $G^{(2n)}(\wt F_{\o\omega} \times \wt F_\omega )$ given by the (bilinear) tensor product 
$\wt M^{2n}_R\otimes \wt M^{2n}_L$ of a right $T^t_n(\wt F_{\o\omega })$-semimodule $\wt M^{2n}_R$ by its left equivalent $\wt M^{2n}_L$.

The bilinear algebraic semigroup $G^{(2n)}(\wt F_{\o\omega }\times\wt F_\omega )$, with entries in the product
$\wt F_{\o\omega }\times\wt F_\omega $ of right pseudoramified complex extensions 
$\wt F_{\o\omega }$ by the left equivalent set $\wt F_\omega $, gives rise by an isomorphism of compactification to the complete algebraic bilinear semigroup $G^{(2n)}( F_{\o\omega }\times F_\omega )$ (which is an abstract bisemivariety) over the product
$F_{\o\omega }\times F_\omega $ of sets of completions.

The linear algebraic semigroup 
$G^{(2n)}(\wt F_\omega )\equiv \wt M^{2n}_L$
\resp{$G^{(2n)}(\wt F_{\o\omega} )\equiv \wt M^{2n}_R$} decomposes into the set 
$\{\wt g^{(2n)}_L(j,m_j)\}_{j,m_j}$
\resp{$\{\wt g^{(2n)}_R(j,m_j)\}_{j,m_j}$} of $r$ packets, $1\le j\le r$, of complex equivalent conjugacy class representatives 
$\wt g^{(2n)}_L(j,m_j)$
\resp{$\wt g^{(2n)}_R(j,m_j)$}.
\vskip 11pt

\item The {\bbf functional representation space 
$\FRepsp (\GL_n(F_{\o\omega} \times F_\omega ))$ of the complete bilinear semigroup
$\GL_n(F_{\o\omega} \times F_\omega )$
\/} is the bisemisheaf $(\widehat M^{2n}_R\otimes \widehat M^{2n}_L)$ of $C^\infty$-differentiable bifunctions on
$( M^{2n}_R\otimes  M^{2n}_L)$, i.e. the (tensor) product of the semisheaf
$\widehat M^{2n}_L$ of $C^\infty$-differentiable functions on $ M^{2n}_L$ by the semisheaf $\widehat M^{2n}_R$ of $C^\infty$-differentiable cofunctions on $ M^{2n}_R$.
\vskip 11pt

\item Let ${\rm CY}^{2n_\ell}(Y_L)
\subset \bZ^{2n_\ell}(Y_L)
\subset {\rm CH}^{2n_\ell}(Y_L)$
\resp{${\rm CY}^{2n_\ell}(Y_R)
\subset \bZ^{2n_\ell}(Y_R)
\subset \linebreak
{\rm CH}^{2n_\ell}(Y_R)
$}
be a \lr algebraic semicycle of dimension $2n_\ell$ on the \lr algebraic semigroup $Y_L\equiv G^{(2n)}(\wt F_\omega )$ \resp{$Y_R\equiv G^{(2n)}(\wt F_{\o\omega })$} of complex dimension $n$, $n_\ell< n$,
where $\bZ^{2n_\ell}(Y_L)$
\resp{$\bZ^{2n_\ell}(Y_R)$} is the semigroup of algebraic semicycles of codimension $2n_\ell$ and
${\rm CH}^{2n_\ell}(Y_L)$
\resp{${\rm CH}^{2n_\ell}(Y_R)$} is the Chow semigroup of algebraic semicycles of codimension $2n_\ell$ on $Y_L$ \resp{$Y_R$} \cite{Jan1}.
\vskip 11pt

\item In this context, it is recalled that {\bbf a Suslin-Voevodsky \lr\linebreak presheaf
$M(X^{\rm sv}_L)$
\resp{$M(X^{\rm sv}_R)$}} on the smooth semischeme
$X^{\rm sv}_L$
\resp{$X^{\rm sv}_R$} of complex dimension $\ell$ on the category 
$\Sm_L(k )$
\resp{$\Sm_R(k)$} of smooth semischemes over 
$k $ 
 is a functor from
$X^{\rm sv}_L$
\resp{$X^{\rm sv}_R$} to the chain complex associated with the abelian semigroup 
$\bigsqcup\limits_{i_{\ell}}\Hom_{\Sm_L(k)}(\dot\Sigma_L,\SP^{i_\ell}(X^{\rm sv}_L))$
\resp{\linebreak$\bigsqcup\limits_{i_{\ell}}\Hom_{\Sm_R(k)}(\dot\Sigma_R,\SP^{i_\ell}(X^{\rm sv}_R))$}
\cite{Mor}
where:
\Bi
\item $\dot\Sigma_L$
\resp{$\dot\Sigma_R$} is a cosimplicial object from the collection of the \lr complex topological $2n_\ell$-simplices
$\Sigma^{2n_\ell}_L$
\resp{$\Sigma^{2n_\ell}_R$}.

\item $\SP^{i_\ell}(X^{\rm sv}_L)$ denotes the $i_\ell$-th symmetric product of 
$X^{\rm sv}_L$.
\Ei
\Ei
\vskip 11pt

A {\bbf Suslin-Voevodsky submotive\/} of dimension $2n_\ell=i_\ell\times 2\ell$ is noted 
$Z_L(2n_\ell)$
\resp{$Z_R(2n_\ell)$} and corresponds to a \lr semicycle 
${\rm CY}^{2n_\ell}(Y_L)$
\resp{${\rm CY}^{2n_\ell}(Y_R)$} in 
$\bZ^{2n_\ell}(Y_L)$
\resp{$\bZ^{2n_\ell}(Y_R)$}.

In order to define Suslin-Voevodsky mixed motives, {\bbf shifted correspondences\/} must be introduced by the homomorphism:
\begin{align*}
\CORR^S_L: \quad \Corr(\SP^{f_\ell}(X^{\rm sv}_L),X^{2n_\ell-2f_\ell\cdot\ell}_L)&\To\Corr^S(\Delta ^{2f_\ell\cdot\ell}_L,X^{2n_\ell-2f_\ell\cdot\ell}_L)\\
\rresp{\CORR^S_R: \quad \Corr(\SP^{f_\ell}(X^{\rm sv}_R),X^{2n_\ell-2f_\ell\cdot\ell}_R)&\To\Corr^S(\Delta ^{2f_\ell\cdot\ell}_R,X^{2n_\ell-2f_\ell\cdot\ell}_R)}\end{align*}
\Bi
\item from correspondences
$\Corr(\SP^{f_\ell}(X^{\rm sv}_L),X^{2n_\ell-2f_\ell\cdot\ell}_L)$ sending the $i_\ell$-th submotive\linebreak
$\SP^{i_\ell}(X^{\rm sv}_L)$ of dimension $2n_\ell=i_\ell\times 2\ell$ to the product
$X^{2n_\ell}_L=\SP^{f_\ell}(X^{\rm sv}_L)\times X^{2n_\ell-2f_\ell\cdot\ell}_L$ where $X^{2n_\ell-2f_\ell\cdot\ell}_L$ is a smooth scheme of complex dimension
$(n_\ell-f_\ell\cdot\ell)$, $f_\ell\cdot\ell\le n_\ell\le n$, $i_\ell\cdot\ell=n_\ell=(f_\ell\cdot\ell)+(n_\ell-f_\ell\cdot\ell)$,

\item to shifted correspondences
$\Corr^S(\Delta ^{2f_\ell\cdot\ell}_L,X^{2n_\ell-2f_\ell\cdot\ell}_L)$ where the smooth semischemes
$\SP^{f_\ell}(X^{\rm sv}_L)$ has been sent to the corresponding smooth semischeme
\[ \Delta^{2f_\ell\cdot\ell}_L=\SP^{f_\ell}(X^{\rm sv}_L)\times \AdFRepsp(T_{f_\ell\cdot\ell}(\CC))\]
where $\Delta^{2f_\ell\cdot\ell}_L$ is the total space of the tangent bundle $\TAN(\SP^{f_\ell}(X^{\rm sv}_L))$ with base space $\SP^{f_\ell}(X^{\rm sv}_L)$ and fibre given by the adjoint functional representation space of the group $T_{f_\ell\cdot\ell}(\CC)$ of triangular matrices.
\Ei
\vskip 11pt

A {\bbf Suslin-Voevodsky \lr mixed semimotive\/} 
$M_{\DM_L(k )}(X^{\rm sv}_L)$
\resp{$M_{\DM_R(k )}(X^{\rm sv}_R)$} can be defined as the functor
\begin{align*}
M_{\DM_L(k )}(X^{\rm sv}_L)&=\bigsqcup_{i_\ell,f_\ell}\Hom_{\Sm_L( k ) }(\SP^{i_\ell}(X^{\rm sv}_L),X^{2n_\ell}_L[2f_\ell\cdot\ell])\\
\rresp{M_{\DM_R(k )}(X^{\rm sv}_R)&=\bigsqcup_{i_\ell,f_\ell}\Hom_{\Sm_R( k) }(\SP^{i_\ell}(X^{\rm sv}_R),X^{2n_\ell}_R[2f_\ell\cdot\ell])}
\end{align*}
where:
\Bi
\item $X^{2n_\ell}_L[2f_\ell\cdot\ell]=\Delta ^{2f_\ell\cdot\ell}_L\times X^{2n_\ell-2f_\ell\cdot\ell}_L$ is the smooth semischeme of dimension $2n_\ell$ shifted in $2f_\ell\cdot\ell$ dimensions by means of the shifted correspondences sending $\SP^{f_\ell}(X^{\rm sv}_L)$ into $\Delta^{2f_\ell\cdot\ell}_L$.

\item $\DM_L(k )$ is the triangulated category of Suslin-Voevodsky left mixed semimotives.
\Ei
\vskip 11pt

The following {\bf proposition} can then be stated ({\bf propositions 2.11 and 2.12}):

Under the action of the adjoint functional representation space
\[
\AdFRepsp(\GL_{f_\ell\cdot\ell}(\cit\times \cit))
=\AdFRepsp(T_{f_\ell\cdot\ell}(\cit))
\times\AdFRepsp(T^t_{f_\ell\cdot\ell}(\cit))\;,
\]
the bilinear cohomology of the \SV\ pure bimotive
$M(X^{\rm sv}\RL)
=M(X^{\rm sv}_R)
\otimes M(X^{\rm sv}_L)$
is transformed into the {\bf the bilinear cohomology of the corresponding \SV\ mixed bimotive\/}
$M_{{\rm DM}\RL}(X^{\rm sv}\RL)
=M_{{\rm DM}_R}(X^{\rm sv}_R)
\otimes M_{{\rm DM}_L}(X^{\rm sv}_L)$ by the isomorphism:
\begin{multline*}
\HH D_{2f_\ell\cdot\ell}:\quad
H^{2n_\ell  }(M(X^{\rm sv}\RL),Z\RL(2n_\ell ))\\
\To H^{2n_\ell -{2f_\ell\cdot\ell}}
(M_{{\rm DM} \RL}(X^{\rm sv}\RL),Z\RL(2n_\ell [{2f_\ell\cdot\ell}]))
\end{multline*}
where:
\Bi
\item $Z\RL(2n_\ell  )=Z_R(2n_\ell )\times Z_L(2n_\ell )$ is the product of \SV\ pure semisubmotives of dimension $2n_\ell $;

\item $Z\RL(2n_\ell [2\fl])
=Z_R(2n_\ell [2\fl])\times Z_L(2n_\ell [2\fl])
\equiv X^{2n_\ell }_R[2\fl])\times 
 X^{2n_\ell }_L [2\fl]$ is the product of \SV\ mixed subsemimotives of dimension $2n_\ell $ shifted in $2\fl$ dimensions;
 
 \item  $H^{2n_\ell -2\fl }(\centerdot,\centerdot)$ is the mixed bilinear cohomology defined in proposition 2.17;
 \Ei
 in such a way that:
 \Bean
 \item $H^{2n_\ell -{2f_\ell\cdot\ell}}
(M_{{\rm DM}\RL}(X^{\rm sv}\RL),Z\RL(2n_\ell [{2f_\ell\cdot\ell}]))
=H_{2\fl}(\Delta ^{2\fl}_R\times \Delta ^{2\fl}_L, \Fs^{2\fl}\RL({\rm TAN}))
\linebreak\times H^{2n_\ell  }(M(X^{\rm sv}\RL),Z\RL(2n_\ell  ))$
where $\Fs^{2\fl}\RL({\rm TAN})$ is the bilinear fibre of the tangent bibundle   
${\rm TAN}({\rm SP}^{f_\ell }(X^{\rm sv}\RL))
={\rm TAN}({\rm SP}^{f_\ell }(X^{\rm sv}_R))
\times{\rm TAN}({\rm SP}^{f_\ell }(X^{\rm sv}_L))$;

\item $H^{2n_\ell -2\fl }(M_{\rm DM\RL}(X^{\rm sv}\RL),Z\RL(2n_\ell [2\fl]))
\simeq\FRepsp(\GL_{n_\ell [\fl]}
((F_{\o\omega }\otimes\cit)
\times(F_{\omega }\otimes\cit)))$
 where
 $\FRepsp(\GL_{n_\ell [\fl]}
((F_{\o\omega }\otimes\cit)
\times(F_{\omega }\otimes\cit)))=
\AdFRepsp(\GL_{\fl}(\cit\otimes\cit))\linebreak
\times\FRepsp(\GL_{n_\ell}(F_{\o\omega }\times F_\omega ))$ is the functional representation space of the bilinear complete semigroup
$\GL_{n_\ell [\fl]}((F_{\o\omega }\otimes \cit)
\times (F_{\omega }\otimes \cit))$ shifted in $(2\fl)$ complex dimensions.
\Ee
\vskip 11pt

To be more explicit, let $D^{2\fl}_R\otimes D^{2\fl}_L$ be the product of a (right) differential (elliptic) operator $D^{2\fl}_R$ acting on $2\fl$ variables by its left equivalent.  {\bbf This bioperator is defined by its action}:
\[ D_R^{2\fl}\otimes D_L^{2\fl}: \qquad
Z_R(2n_\ell )\times Z_L(2n_\ell )
\To Z_R(2n_\ell [2\fl])\times Z_L(2n_\ell [2\fl] )
\]
transforming the \SV\ pure subbisemimotive of dimension $(2n_\ell )$ into the corresponding \SV\ mixed subbisemimotive of dimension $(2n_\ell)$ shifted in $(2\fl)$ dimensions.

Indeed, it is seen in chapter 2 that $H_{2\fl}(\Delta ^{2\fl}_R\times \Delta ^{2\fl}_L,\Fs^{2\fl}\RL({\rm TAN}))\simeq \linebreak \AdFRepsp (\GL_{\fl}(\cit\times\cit))$ is the bilinear homology with coefficients in the bilinear fibre $\Fs^{2\fl}\RL({\rm TAN})$ of the tangent bibundle ${\rm TAN}[{\rm SP}^{f_\ell }(X^{\rm sv}_R)\times {\rm SP}^{f_\ell }(X^{\rm sv}_L)]$.
\vskip 11pt

In connection with the work of G. Kasparov \cite{Kas}, we shall introduce in chapter 3 {\bbf a $K_*K^*$ functor} on the categories of elliptic bioperators and products, right by left, of \SV\ pure motives allowing to set up a {\bbf bilinear version of the index theorem}.

\Bean
\item If $H^*(M(X^{\rm sv}\RL))=
\bigoplus\limits^{n_s}_{n_\ell =n_1}
H^{2n_\ell  }(M(X^{\rm sv}\RL),Z\RL(2n_\ell  ))$ denotes the total bilinear cohomology of the pure bimotive $M(X^{\rm sv}\RL)$ and if $K^*(X^{\rm sv}\RL)$, introduced as the product, right by left, of abelian semigroups generated by the complex vector bundles over $X^{\rm sv}\RL=
X^{\rm sv}_R\times
X^{\rm sv}_L$, is the $K$-cohomology associated with the pure bimotive $M(X^{\rm sv}\RL)$, {\bbf the total Chern character in the bilinear $K$-cohomology \cite{W-R} of the pure bimotive $M(X^{\rm sv}\RL)$} is given by the homomorphism:
\[ {\rm ch}^*(M(X^{\rm sv}\RL )):\qquad
K^*(X^{\rm sv}\RL )\To H^*(M(X^{\rm sv}\RL ))\;.\]

\item Similarly, if
$H_*(\Delta ^*_R\times\Delta ^*_L,\Fs^*\RL({\rm TAN}))
=\bigoplus\limits_{\fl} H_{2\fl}(\Delta ^{2\fl}_R\times \Delta ^{2\fl}_L,\Fs^{2\fl}\RL({\rm TAN}))
\simeq \bigoplus\limits_{\fl}\AdFRepsp (\GL_{\fl}(\cit\times \cit))$ is the total bilinear homology with coefficients in the set of bilinear fibres $\Fs^{2\fl}\RL({\rm TAN})$ and if $K_*({\rm SP}^{\rm FL}(X^{\rm sv}\RL	))$ is the bilinear $K$-homology, introduced as the product, right by left, of abelian semigroups generated by the set of tangent bibundles ${\rm TAN}({\rm SP}^{f_\ell }(X^{\rm sv}\RL))$, {\bbf the Chern character in the bilinear $K$-homology, associated with the pure bimotive $M(X^{\rm sv}\RL)$}, is given by the homomorphism:
\[ {\rm ch}_*(M(X^{\rm sv}\RL)): \qquad
K_*({\rm SP}^{\rm FL} ( X^{\rm sv}\RL	))
\To H_*(\Delta ^*_R\times\Delta ^*_L,\Fs^*\RL(\TAN) );\]

\item {\bbf The total Chern character
${\rm ch}^*(M_{{\rm DM}\RL}(X^{\rm sv}\RL))$ of the \SV\ mixed bisemimotive
$M_{{\rm DM}\RL}(X^{\rm sv}\RL)$} in the mixed bilinear $K$-homology-$K$-\linebreak cohomology is given by the homomorphism:
\begin{multline*} {\rm ch}^*(M_{{\rm DM}\RL}(X^{\rm sv}\RL)):
\quad K_*({\rm SP}^{\rm FL} ( X^{\rm sv}\RL	))\times
K^*( X^{\rm sv}\RL	)\qquad\\
\qquad\To H_*(\Delta ^*_R\times\Delta ^*_L,\Fs^*\RL(\TAN) )\times H^*
(M(X^{\rm sv}\RL));
\end{multline*}
in such a way that
\[{\rm ch}_*(M(X^{\rm sv}\RL)\times
{\rm ch}^*(M(X^{\rm sv}\RL)
\To 
{\rm ch}^*(M_{{\rm DM}\RL}(X^{\rm sv}\RL))\]
corresponds to a bilinear version of the index theorem.
\Ee
\vskip 11pt

Chapter 4 deals with the {\bbf holomorphic and toroidal spectral representations of an elliptic bioperator} associated with the functional representation space\linebreak
$\FRepsp (\GL_{n_\ell [\fl]}
(F_{\o\omega }\otimes\cit)\times
(F_{\omega }\otimes\cit))$ of the complete  bilinear semigroup
$\GL_{n_\ell [\fl]}
(F_{\o\omega }\otimes\cit)\times
(F_{\omega }\otimes\cit)$ shifted in
$(2\fl)$ dimensions.

Taking into account that:
\Bena
\item the functional representation space
\cite{Del1}, \cite{Vog} $
\FRepsp (\GL_{n_\ell }
(F_{\o\omega }\otimes F_{\omega }))$ of the complete bilinear semigroup
$\GL_{n_\ell }
(F_{\o\omega }\times
F_{\omega })$ is the bisemisheaf
$(\widehat M^{2n_\ell }_R\otimes \widehat M^{2n_\ell }_L)$ of differentiable bifunctions over
$\GL_{n_\ell }
(F_{\o\omega }\times
(F_{\omega })$,

\item there exists a toroidal isomorphism of compactification:
\[ \gamma \RL: \qquad \widehat M^{2n_\ell }_R\otimes
\widehat M^{2n_\ell }_L
\To \widehat M^{2n_\ell }_{T_R}\otimes
\widehat M^{2n_\ell }_{T_L}\]
sending $(\widehat M^{2n_\ell }_R\otimes
\widehat M^{2n_\ell }_L)$ into its toroidal equivalent
$(\widehat M^{2n_\ell }_{T_R}\otimes
\widehat M^{2n_\ell }_{T_L})
=\linebreak \FRepsp(\GL_{n_\ell }(F^T_{\o\omega }\times F^T_\omega ))
$ where $F^T_\omega $ and $F^T_{\o\omega }$ are sets of toroidal completions,

\item there exists a correspondence:
\[ \HH_{\bigoplus}: \qquad
\widehat M^{2n_\ell }_{T_R}\otimes
\widehat M^{2n_\ell }_{T_L} \To
\widehat M^{2n_\ell }_{T_{R_{\bigoplus}}}\otimes
\widehat M^{2n_\ell }_{T_{L_{\bigoplus}}}  \]
in such a way that
\[ \widehat M^{2n_\ell }_{T_{R_{\bigoplus}}}\otimes
\widehat M^{2n_\ell }_{T_{L_{\bigoplus}}}
=\bigoplus_j\bigoplus_{m_j}
(\widehat M^{2n_\ell }_{T_{\o\omega _{j,m_j}}}\otimes
\widehat M^{2n_\ell }_{T_{\omega _{j,m_j}}})\]
decomposes into the sum of bisections
$(\widehat M^{2n_\ell }_{T_{\o\omega _{j,m_j}}}\otimes
\widehat M^{2n_\ell }_{T_{\omega _{j,m_j}}})$
of $(\widehat M^{2n_\ell }_{T_R}\otimes
\widehat M^{2n_\ell }_{T_L})$
according to the conjugacy class representatives of
$(  \widehat M^{2n_\ell }_{T_R}\otimes
 \widehat M^{2n_\ell }_{T_L})$,
\Ee
{\bbf the elliptic bioperator
$(D^{2\fl}_R\otimes D^{2\fl}_L)$ maps
$(\widehat M^{2n_\ell }_{T_{R_{\bigoplus}}}\otimes
\widehat M^{2n_\ell }_{T_{L_{\bigoplus}}})$ into its shifted equivalent according to}:
\[D^{2\fl}_R\otimes D^{2\fl}_L: \qquad
\widehat M^{2n_\ell }_{T_{R_{\bigoplus}}}\otimes
\widehat M^{2n_\ell }_{T_{L_{\bigoplus}}}
\To
\widehat M^{2n_\ell }_{T_{R_{\bigoplus}}}[2\fl]\otimes
\widehat M^{2n_\ell }_{T_{L_{\bigoplus}}}[2\fl]\]
where:
\[(\widehat M^{2n_\ell }_{T_{R}}[2\fl]\otimes
\widehat M^{2n_\ell }_{T_{L}}[2\fl])
= \FRepsp(\GL_{n_\ell [\fl]}
((F^T_{\o\omega }\otimes \cit)\times
((F^T_{\omega }\otimes \cit))\]
is the perverse bisemisheaf of differentiable bifunctions over
$\GL_{n_\ell [\fl]}
((F^T_{\o\omega }\otimes \cit)\times
((F^T_{\omega }\otimes \cit))$ shifted in $(2\fl)$ dimensions: it is thus a $(\Ds_R\otimes \Ds_L)$-bisemimodule in such a way that $\Ds_R$ \resp{$\Ds_L$} is a \rl sheaf of differentiable operators of finite order with holomorphic coefficients \cite{M-T}.
\vskip 11pt

\Bi
\item Referring to \cite{Pie3}, we see that each bifunction 
$(\widehat M^{2n_\ell }_{T_{\o\omega _{j,m_j}}}\otimes
\widehat M^{2n_\ell }_{T_{\o\omega _{j,m_j}}})\in
(\widehat M^{2n_\ell }_{T_{R_{\bigoplus}}}\otimes
\widehat M^{2n_\ell }_{T_{L_{\bigoplus}}})$ is the product, right by left , of $n_\ell $-dimensional complex semitori:
\begin{align*}
\widehat M^{2n_\ell }_{T_{\o\omega _{j,m_j}}}
&=
T^{2n_\ell }_R(j,m_j)=\lambda ^{\half}(2n_\ell ,j,m_j)\ e^{-2\pi ijz_{n_\ell} }\\
{\rm and} \quad 
\widehat M^{2n_\ell }_{T_{\omega _{j,m_j}}}
&=
T^{2n_\ell }_L(j,m_j)=\lambda ^{\half}(2n_\ell ,j,m_j)\ e^{2\pi ijz_{n_\ell} }\;, \quad z_{n_\ell }\in\cit^{n_\ell }\;,\end{align*}
in such a way that
$\widehat M^{2n_\ell }_{T_{L_{\bigoplus}}}$
\resp{$\widehat M^{2n_\ell }_{T_{R_{\bigoplus}}}$}
is the (truncated) {\bbf Fourier development of a normalized $2n_\ell $-dimensional left \resp{right} cusp form of weight 2} restricted to the upper \resp{lower} half space:
\begin{align*}
\widehat M^{2n_\ell }_{T_{L_{\bigoplus}}}
\equiv \EIS_L(2n_\ell ,j,m_j)
&=\bigoplus^r_{j=1}
\bigoplus_{m_j}\lambda ^{\half}(2n_\ell ,j,m_j)\ 
e^{2\pi ijz_{n_\ell }}\\
\rresp{
\widehat M^{2n_\ell }_{T_{R_{\bigoplus}}}
\equiv \EIS_R(2n_\ell ,j,m_j)
&=\bigoplus^r_{j=1}
\bigoplus_{m_j}\lambda ^{\half}(2n_\ell ,j,m_j)\ 
e^{-2\pi ijz_{n_\ell }}}.
\end{align*}

\item Similarly, $\widehat M^{2n_\ell }_{T_{L_{\bigoplus}}} [2\fl]$ \resp{$\widehat M^{2n_\ell }_{T_{R_{\bigoplus}}} [2\fl]$} decomposes into sums of $2n_\ell $-dimensional semitori shifted in $2\fl$ dimensions in such a way that:
\begin{align*}
\widehat M^{2n_\ell }_{T_{L_{\bigoplus}}}[2\fl]
&\equiv \EIS_L(2n_\ell [2\fl],j,m_j)\\
&\qquad =\bigoplus^r_{j=1}
\bigoplus_{m_j}
E_{2\fl}(2n_\ell ,j,m_j)\centerdot
\lambda ^{\half}(2n_\ell ,j,m_j)\ 
e^{2\pi ijz_{n_\ell }}\\
\rresp{
\widehat M^{2n_\ell }_{T_{R_{\bigoplus}}}[2\fl]
&\equiv \EIS_R(2n_\ell [2\fl],j,m_j)\\
&\qquad =\bigoplus^r_{j=1}
\bigoplus_{m_j}
E_{2\fl}(2n_\ell ,j,m_j)\centerdot
\lambda ^{\half}(2n_\ell ,j,m_j)\ 
e^{-2\pi ijz_{n_\ell }}}
\end{align*}
be the (truncated) {\bbf Fourier development of a normalized \lr $2n_\ell $-dimensional mixed cusp form shifted in $2\fl$ dimensions, where $E_{2\fl}(2n_\ell ,j,m_j)$ are shifts of generalized global Hecke characters\linebreak $\lambda ^{\half}(2n_\ell ,j,m_j)$.}
\vskip 11pt

\item This allows to set up the {\bf bieigenvalue equation}:
\begin{multline*}
(D_R^{2\fl}\otimes D_L^{2\fl} )
(\EIS_R (2n_\ell ,j^{\rm up}=j,m_j))\otimes
(\EIS_L (2n_\ell ,j^{\rm up}=j,m_j))\\
=E^2_{2\fl} (2n_\ell ,j,m_j)
(\EIS_R (2n_\ell ,j^{\rm up}=j,m_j))\otimes
(\EIS_L (2n_\ell ,j^{\rm up}=j,m_j))
\end{multline*}
{\bbf of which spectral representation is given by the set of $r$-bituples}:
\begin{multline*}
 \bigl\{ 
(\EIS_R (2n_\ell ,j^{\rm up}=1,m_1))\otimes
(\EIS_L (2n_\ell ,j^{\rm up}=1,m_1)),\cdots,\bigr.\\
(\EIS_R (2n_\ell ,j^{\rm up}=j,m_j))\otimes
(\EIS_L (2n_\ell ,j^{\rm up}=j,m_j)),\cdots,\\
\bigl.(\EIS_R (2n_\ell ,j^{\rm up}=r,m_r))\otimes
(\EIS_L (2n_\ell ,j^{\rm up}=r,m_r))
\bigr\}
\end{multline*}
{\bbf where
$\EIS_L (2n_\ell ,j^{\rm up}=j,m_j))$ \resp{$\EIS_R (2n_\ell ,j^{\rm up}=j,m_j))$} is the truncated Fourier development at the  $j$ classes of the $2n_\ell $-dimensional cusp form.}
\vskip 11pt

\item It then appears that 
$\EIS_R (2n_\ell [2\fl],j,m_j)\otimes
\EIS_L (2n_\ell [2\fl],j,m_j)$ constitutes a {\bbf supercuspidal representation of the shifted algebraic complete semigroup\linebreak
$\GL_{n_\ell [\fl]}(
(F_{\o\omega _{\oplus}}\otimes \cit)\times
(F_{\omega _{\oplus}}\otimes \cit))$}.
\Ei
\vskip 11pt

{\bbf The origin of the (bilinear) spectral theory then results from geometric-shifted global (bilinear) correspondences of Langlands as it will be seen hereafter.}
\vskip 11pt

This leads us to develop {\bbf in 
chapter 5 geometric-shifted global bilinear correspondences of Langlands}.

\Bi
\item If 
$(W^{ab}_{F^{S_{\cit}}_{\o\omega }}\times
W^{ab}_{F^{S_{\cit}}_{\omega }})$ is the product, right by left, of the shifted global Weil group 
$W^{ab}_{F^{S_{\cit}}_{\o\omega }}$ and $
W^{ab}_{F^{S_{\cit}}_{\omega }}$ introduced in chapter 1, there exists an irreducible representation
$\IrrRep^{(2n_\ell [2\fl])}_{W_{F\RL}}
(W^{ab}_{F^{S_{\cit}}_{\o\omega }}\times
W^{ab}_{F^{S_{\cit}}_{\omega }})$
of $(W^{ab}_{F^{S_{\cit}}_{\o\omega }}\times
W^{ab}_{F^{S_{\cit}}_{\omega }})$ given by the representation space
\[
\Repsp (\GL_{n_\ell [\fl]}
( ( F_{\o\omega _{\bigoplus}}\otimes\cit )\times
( F_{\omega _{\bigoplus}}\otimes\cit ) ))
\begin{aligned}[t]
&\equiv G^{(2n_\ell [2\fl])}
( ( F_{\o\omega _{\bigoplus}}\otimes\cit )\times
( F_{\omega _{\bigoplus}}\otimes\cit ) )\\
&\equiv M^{2n_\ell }_{R_{\bigoplus}}[2\fl]\otimes
M^{2n_\ell }_{L_{\bigoplus}}[2\fl]
\end{aligned}\]
of the shifted bilinear complete semigroup
$\GL_{n_\ell [\fl]}
( ( F_{\o\omega _{\bigoplus}}\otimes\cit )\times
( F_{\omega _{\bigoplus}}\otimes\cit ) )$.

So, {\bbf on the shifted irreducible bilinear complete semigroup 
$G^{(2n_\ell [2\fl])}
( ( F_{\o\omega _{\bigoplus}}\linebreak \otimes\cit )\times
( F_{\omega _{\bigoplus}}\otimes\cit ) )$, the geometric-shifted global bilinear correspondence of Langlands is}:

{ \begin{small}
\[\begin{psmatrix}[colsep=.5cm,rowsep=1cm]
\IrrRep^{(2n_\ell [2\fl])}_{W_{F\RL}}(W^{ab}_{F^{S_{\cit}}_{\o\omega }}\times
W^{ab}_{F^{S_{\cit}}_{\omega }}) 
& \raisebox{3mm}{$\sim$}
&
\Irrcusp(\GL_{n_\ell [\fl]}(({F^T_{\o\omega }}\otimes\cit)\times
({F^T_{\omega }}\otimes\cit)))
 \\
G^{(2n_\ell [2\fl])}((F_{\o\omega _{\bigoplus}}\otimes\cit)\times
(F_{\omega _{\bigoplus}}\otimes\cit))
 & &
\EIS\RL(2n_\ell [2\fl],j,m_j)
\\
 && G^{(2n_\ell [2\fl])}
( ( {F^T_{\o\omega }}\otimes\cit )\times
( {F^T_{\omega }}  \otimes\cit ))
\ncline[arrows=->,nodesep=5pt]{1,1}{1,3}
\ncline[doubleline=true,nodesep=5pt]{1,1}{2,1}
\ncline[doubleline=true,nodesep=5pt]{1,3}{2,3}
\ncline[arrows=->,nodesep=5pt]{2,1}{3,3}^{\rotatebox{-30}{$\sim$}}
\ncline[arrows=->,nodesep=5pt]{3,3}{2,3}<{\rotatebox{90}{$\sim$}}
\end{psmatrix}\]
\end{small}}
where $\Irrcusp(\GL_{n_\ell [\fl]}(({F^T_{\o\omega }}\otimes\cit)\times
({F^T_{\omega }}\otimes\cit)))$ is the shifted irreducible supercuspidal representation of
$\GL_{n_\ell [\fl]}(({F^T_{\o\omega }}\otimes\cit)\times
({F^T_{\omega }}\otimes\cit))$ over the product of toroidal shifted completions.
\vskip 11pt

\item Similarly, {\bbf on the reducible shifted $2n$-dimensional bilinear complete algebraic semigroup} 
$G^{2n [2\fnl]}
( ( F_{\o\omega _{\bigoplus}}\otimes\cit )\times
( F_{\omega _{\bigoplus}}\otimes\cit ) )
= \mathop{\boxplus}\limits_{n_\ell =n_1}^{n_s}G^{2n_\ell [2\fl]}
( ( F_{\o\omega _{\bigoplus}}\otimes\cit )\times
( F_{\omega _{\bigoplus}}\otimes\cit ) )$, there exists the geometric-shifted global bilinear reducible correspondence of Langlands:
{ \begin{small}
\[\begin{psmatrix}[colsep=.5cm,rowsep=1cm]
\RedRep^{(2n [2\fnl])}_{W_{F\RL}}(W^{ab}_{F^{S_{\cit}}_{\o\omega }}\times
W^{ab}_{F^{S_{\cit}}_{\omega }}) 
& \raisebox{3mm}{$\sim$}
&
\Redcusp(\GL_{n [\fnl]}(({F^T_{\o\omega }}\otimes\cit)\times
({F^T_{\omega }}\otimes\cit)))
 \\
G^{(2n  [2\fnl])}((F_{\o\omega _{\bigoplus}}\otimes\cit)\times
(F_{\omega _{\bigoplus}}\otimes\cit))
 & &
\EIS\RL(2n_\ell [2\fl],j,m_j)
\\
 && G^{(2n [2\fnl])}
( ( {F^T_{\o\omega }}\otimes\cit )\times
( {F^T_{\omega }} \otimes\cit) )
\ncline[arrows=->,nodesep=5pt]{1,1}{1,3}
\ncline[doubleline=true,nodesep=5pt]{1,1}{2,1}
\ncline[doubleline=true,nodesep=5pt]{1,3}{2,3}
\ncline[arrows=->,nodesep=5pt]{2,1}{3,3}^{\rotatebox{-30}{$\sim$}}
\ncline[arrows=->,nodesep=5pt]{3,3}{2,3}<{\rotatebox{90}{$\sim$}}
\end{psmatrix}\]
\end{small}}
\vskip 11pt

\item Geometric-shifted global bilinear correspondences of Langlands are also established on real shifted irreducible and reducible bilinear complete algebraic semigroups.
\vskip 11pt

\item Remark that the geometric-shifted global bilinear correspondences of Langlands considered in this paper differ from {\bf the geometric correspondences\/} initiated by V. Drinfeld and G. Laumon.

Indeed, these deal with an $\ell$-adic $n$-dimensional irreducible local system $E$ on a smooth algebraic curve $X$ over a ground field $K$ and say that it can be associated to $E$ an automorphic sheaf $S_E$ which is a perverse sheaf on the moduli stack  ${\rm Bun}_n(X)$ of vector bundles of rank $n$ on $X$ \cite{Lau}, \cite{F-G-V}, \cite{Fre}, \cite{Gai}.

\item The last version of this paper was motivated to precise the nature of a general bilinear mixed cohomology theory.
\Ei

%% file: LanglandsIII_1.tex
\section{Global class field concepts and pure motivic cohomologies}

\subsection[Pseudo-ramified and pseudo-unramified infinite places of semifields]{Pseudo-ramified and pseudo-unramified infinite places of\linebreak semifields}

\Bi
\item Let $k$ be a  number field of characteristic $0$ and let $\wt F$ denote a finite extension set of $k$ such that $\wt F$ is assumed to be a {\bbf symmetric splitting field} $\wt{F}=\wt F_R\cup \wt F_L$ composed of the right and left algebraic extension semifields $\wt F_R$ and $\wt F_L$ in one-to-one correspondence.
\vskip 11pt

\item $\wt F_L$
\resp{$\wt F_R$} is assumed to be  composed of a set of complex (resp. conjugate complex) simple roots of a polynomial ring over $k$.  If the algebraic extension field of $k$ is real, then the symmetric splitting field will be noted $\wt{F}^+=\wt F^+_R \cup \wt F^+_L$ where the \lr algebraic extension semifield 
$\wt{F}^+_L$
\resp{$\wt{F}^+_R$} is composed of the set of positive \resp{symmetric negative} simple real roots.
\vskip 11pt

\item The left and right equivalence classes of the global completions of 
$\wt{F}^{(+)}_L$ and
$\wt{F}^{(+)}_R$
(which correspond to transcendental extensions of $k$), obtained by an isomorphism of compactification of the corresponding finite extensions, are the left and right infinite real \resp{complex} places of 
${F}^{(+)}_L$ and
${F}^{(+)}_R$: they are noted
$v=\{ v_1,\dots,v_{j_\delta },\dots,v_{r_\delta }\}$ and
$\o v=\{ \o v_1,\dots,\o v_{j_\delta },\dots,\o v_{r_\delta }\}$ in the real case and
$\omega =\{ \omega _1,\dots,\omega _{j},\dots,\omega _{r}\}$ and
$\o \omega =\{ \o \omega _1,\dots,\o \omega _{j},\dots,\o \omega _{r}\}$ in the complex case.
\vskip 11pt

\item {\bbf The pseudo-unramified real places\/} are characterized algebraically by their global class residue degrees 
$f_{v_{j_\delta }}$ and
$f_{\o v_{j_\delta }}$ given by 
$f_{v_{j_\delta }}=[\wt F^{+,nr}_{v_{j_\delta }}:k]= j$ and
$f_{\o v_{j_\delta }}=[\wt F^{+,nr}_{\o v_{j_\delta }}:k]= j$, $j\in\nit$, $1\le j_\delta \le r_\delta $, where
$\wt F^{+,nr}_{v_{j_\delta }}$ and
$\wt F^{+,nr}_{\o v_{j_\delta }}$ denote  basic real pseudo-unramified extensions (splitting subsemifields) of
$k$ in one-to-one correspondence with the corresponding completions
$F_{v_{j_\delta }}^{+,nr}$ and
$ F_{\o v_{j_\delta }}^{+,nr}$ at the places $v_{j_\delta}$ and $\o v{j_\delta }$.

Similarly, {\bbf pseudo-unramified complex places} are characterized by their global class residue degrees $f_{\omega _j}$ and
$f_{\o \omega _j}$ given by:
\[ f_{\omega _j}=[\wt F^{nr}_{\omega _j}:k]= j \qquad
\text{and} \qquad
 f_{\o \omega _j}=[\wt F^{nr}_{\o \omega _j}:k]= j\]
 where  $\wt F^{nr}_{\omega _j}$ and
  $\wt F^{nr}_{\o \omega _j}$ denote complex basic pseudo-unramified extensions of $k$ in one-to-one correspondence with the corresponding completions 
 $ F^{nr}_{\omega _j}$ and
 $ F^{nr}_{\o \omega _j}$ at the places $\omega _j$ and $\o\omega _j$.
 \vskip 11pt
 
 {\bbf Infinite pseudo-ramified real places} are assumed to be also characterized by Galois extension degrees: they are in fact classes of completions of which degrees are given by integers modulo $N$, $\zit/N\ \zit$, as follows:
 \[ [\wt F^+_{v_{j_\delta }}:k] 
= *+j\ N \qquad \text{and} \qquad
[\wt F^+_{\o v_{j_\delta }}:k] 
= *+j\ N\]
 where:
 \Bi
 \item $F^+_{v_{j_\delta }}$ and
$F^+_{\o v_{j_\delta }}$ denote respectively a real basic ramified completion of $F^+_L$ and of $F^+_R$ in one-to-one correspondence with the splitting subsemifields
$\wt F^+_{v_{j_\delta }}$ and
$\wt F^+_{\o v_{j_\delta }}$;
\item $*$ denotes an integer inferior to $N$.
\Ei

And, {\bbf infinite pseudo-ramified complex places} are similarly characterized by degrees given by the integers modulo $N$, $\zit/N\ \zit$, according to
 \[ [\wt F_{\omega _{j}}:k] = (*+j\ N)m^{(j)} \qquad \text{and} \qquad
 [\wt F_{\o \omega _{j}}:k] = (*+j\ N) m^{(j)}\]
 where:
 \Bi
 \item $F_{\omega _j}$ and
  $F_{\o \omega _j}$ are respectively the basic complex pseudo-ramified completions of $  F_L$ and of $  F_R$ in one-to-one correspondence with the corresponding splitting subsemifields 
$\wt F_{\omega _j}$ and
$\wt F_{\o \omega _j}$;
\item $m^{(j)}=\sup(m_{j_\delta }+1)$ is the {\bbf multiplicity of the $j_\delta $-th real completion covering its $j$-th complex equivalent} or the number of compactified divisors of $F_{\omega _j}$ and of $F_{\o\omega _j}$.
\Ei
\vskip 11pt

\item {\bbf The origin of the integer $N$} in the real case results from the fact that the real pseudo-ramified completions 
$F^+_{v_{j_\delta }}$ and
$F^+_{\o v_{j_\delta }}$ are assumed to be generated respectively from the irreducible central subcompletions
$F^+_{v^1_{j_\delta }}$ and
$F^+_{\o v^1_{j_\delta }}$ characterized by a (Galois extension) degree $[\wt F^+_{v^1_{j_\delta }}:k]=N$ and
$[\wt F^+_{\o v^1_{j_\delta }}:k]=N$.

Similarly, the complex pseudo-ramified completions
$F_{\omega _j}$ and
$F_{\o\omega _j}$ are generated respectively from equivalent subcompletions 
$F_{\omega _j^1}$ and
$F_{\o\omega _j^1}$ having a degree or rank equal to $N\centerdot m^{(j)}$.
\vskip 11pt

\item On the other hand, as {\bf a place} is an equivalence class of completions, we have to consider a set of:
\Bi
\item real pseudo-ramified \resp{pseudo-unramified} completions
$\{F^{+(,nr)}_{v_{j_\delta ,m_{j_\delta }}}\}_{m_{j_\delta }}$ and
$\{F^{+(,nr)}_{\o v_{j_\delta ,m_{j_\delta }}}\}_{m_{j_\delta }}$, $1\le j_\delta \le r_\delta $, equivalent respectively to the corresponding basic completions
$F^{+(,nr)}_{v_{j_\delta }}$ and
$F^{+(,nr)}_{\o v_{j_\delta }}$, where $m_{j_\delta }\ge 1$ is an increasing integer such that $m^{(j_\delta )}=\sup(m_{j_\delta })$ denotes the multiplicity of 
$F^{+(,nr)}_{v_{j_\delta }}$ and of
$F^{+(,nr)}_{\o v_{j_\delta }}$;

\item complex pseudo-ramified \resp{pseudo-unramified} completions
$\{F^{(nr)}_{\omega _{j,m_j }}\}_{m_j}$ and
$\{F^{(nr)}_{\o \omega _{j,m_j }}\}_{m_j}$, $1\le j\le r$, equivalent respectively to the corresponding basic completions $F^{(nr)}_{\omega _j}$ and
$F^{(nr)}_{\o \omega _j}$ where $m_j\ge 1$ is an increasing integer such that $m_\omega ^{(j)}=\sup(m_j)$ refers to the multiplicity of 
$F^{(nr)}_{\omega _j}$ and
$F^{(nr)}_{\o\omega _j}$.
\Ei
\vskip 11pt

\item All the {\bbf real pseudo-ramified completions}
$F^{+}_{v_{j_\delta },m_{j_\delta }}$ 
\resp{$F^{+}_{\o v_{j_\delta },m_{j_\delta }}$}, $m_{j_\delta }\ge1$, in a place 
$v_{j_\delta }$
\resp{$\o v_{j_\delta }$}, are characterized by the same (Galois extension) degree $\simeq j\centerdot N$ and {\bbf and are cut into $j$ irreducible equivalent real subcompletions}
$F_{v^{j'_\delta }_{j'_\delta }}$, $1\le j'_\delta \le j_\delta $, having a degree equal to $N$.  

In the same manner, the complex pseudo-ramified completions
$F_{\omega _{j,m_j}}$
\resp{$F_{\o \omega _{j,m_j}}$}, $m_j\ge 1$, in a place $\omega _j$ \resp{$\o\omega _j$} are characterized by the same degree $\simeq j\centerdot m^{(j)}\centerdot N$ and are cut into $j$ equivalent complex subcompletions $F_{\omega ^{j' }_j}$, $1\le j'\le j$, having a degree equal to $m^{(j  )}\centerdot N$.
 \Ei
 \vskip 11pt
 
 \subsection[Definition: Infinite pseudo-ramified adele semirings and semigroups $F_{\omega _\oplus}$ and $F_{v_\oplus}$]{\bbf Definition: Infinite pseudo-ramified adele semirings and semi\-groups $F_{\omega _\oplus}$ and $F_{v_\oplus}$}
 
 \Bi
 \item Infinite pseudo-ramified adele semirings $\Aa^\infty _{F_v}$, $\Aa^\infty _{F_{\o v}}$, $\Aa^\infty _{F_\omega }$ and $\Aa^\infty _{F_{\o \omega }}$ can be introduced by considering the products of the basic completions  over primary places of respectively $F^+_L$, $F^+_R$, $F_L$ and  $F_R$ according to:
 \begin{align*}
 \Aa^\infty _{F^+_v} &= \prod_{j_{\delta _p}} F^+_{v_{j_{\delta _p}}}
 \;, &\Aa^\infty _{F^+_{\o v}} &= \prod_{j_{\delta _p}} F^+_{\o v_{j_{\delta _p}}}
 \;, && 1\le j_{\delta _p}\le r_\delta \le \infty\;,  \\[15pt]
 \Aa^\infty _{F_\omega } &= \prod_{j_{p}} F_{\omega _{j_p}}
 \;, & \Aa^\infty _{F_{\o \omega }} &= \prod_{j_{p}} F_{\o \omega _{j_p}}\;,&& 1\le j_{p}\le r \le \infty\;.
\end{align*}
 
 \item Similarly, direct sums of completions will be given by:
 \begin{align*}
 F^+_{v_\oplus} &= \bigoplus_{j_\delta } \bigoplus_{m_{j_\delta }} F^+_{v_{j_\delta ,m_{j_\delta }}}\;, 
 & F^+_{\o v_\oplus} &= \bigoplus_{j_\delta } \bigoplus_{m_{j_\delta }} F^+_{\o v_{j_\delta ,m_{j_\delta }}}\;, 
 \\[15pt]
 F_{\omega _\oplus} &= \bigoplus_{j } \bigoplus_{m_{j }} F_{\omega _{j ,m_{j }}}\;, 
 & F_{\o \omega _\oplus} &= \bigoplus_{j } \bigoplus_{m_{j }} F_{\o \omega _{j  ,m_{j  }}}\;, \end{align*}
 \Ei

\subsection{Global inertia subgroups}

\Bi
\item Let 
$\wt F_{\omega _{j,m_j}}$
\resp{$\wt F_{\o\omega _{j,m_j}}$}, $m_j>1$, denote a complex pseudo-ramified extension corresponding to the respective pseudo-ramified completion
$  F_{\omega _{j,m_j}}$
\resp{$ F_{\o\omega _{j,m_j}}$} and approximatively equivalent to $j$ basic complex pseudo-ramified extension
$\wt F_{\omega _{j,1}}$
\resp{$\wt F_{\o\omega _{j,1}}$}, $m_j=1$.
\vskip 11pt

\item Respectively, let
$\{\wt F^{nr}_{\omega _{j,m_j}}\}^{m_j}_{m_j=1}$
\resp{$\{\wt F^{nr}_{\o\omega _{j,m_j}}\}^{m_j}_{m_j=1}$} denote the set of complex pseudo-unramified extensions corresponding to the respective pseudo-unramified completions at the $j$-th pseudo-unramified complex place.
\vskip 11pt

\item Let $\Gal(\wt F^{nr}_{\omega _{j,m_j}}/k)$
\resp{$\Gal(\wt F^{nr}_{\o\omega _{j,m_j}}/k)$}
 be the Galois subgroup of the pseudo-unramified complex extension
$\wt F^{nr}_{\omega _{j,m_j}}$
\resp{$\wt F^{nr}_{\o\omega _{j,m_j}}$} of $k$ and let
$\Gal(\wt F_{\omega _{j,m_j}}/k)$
\resp{$\Gal(\wt F_{\o\omega _{j,m_j}}/k)$} be the Galois subgroup of the pseudo-ramified complex extension 
$\wt F_{\omega _{j,m_j}}$
\resp{$\wt F_{\o\omega _{j,m_j}}$} of $k$.
\vskip 11pt

\item Then, the {\bbf global inertia subgroup}
$I_{\wt F_{\omega _{j,m_j}}}$
\resp{$I_{\wt F_{\o\omega _{j,m_j}}}$} of 
$\Gal(\wt F_{\omega _{j,m_j}}/k)$
\resp{$\Gal(\wt F_{\o\omega _{j,m_j}}/k)$} will be defined by
\begin{align*}
\Gal(\wt F_{\omega _{j,m_j}}/k)&=
\Gal(\wt F^{nr}_{\omega _{j,m_j}}/k)
\times I_{\wt F_{\omega _{j,m_j}}}\\
\rresp{\Gal(\wt F_{\o\omega _{j,m_j}}/k)&=
\Gal(\wt F^{nr}_{\o\omega _{j,m_j}}/k)
\times I_{\wt F_{\o\omega _{j,m_j}}}}
\end{align*}
which leads to the exact sequence
\begin{align*}
1 &\longrightarrow I_{\wt F_{\omega _{j,m_j}}}
\longrightarrow \Gal(\wt F_{\omega _{j,m_j}}/k)
\longrightarrow \Gal(\wt F^{nr}_{\omega _{j,m_j}}/k)
\longrightarrow 1\\
\rresp{1 &\longrightarrow I_{\wt F_{\o\omega _{j,m_j}}}
\longrightarrow \Gal(\wt F_{\o\omega _{j,m_j}}/k)
\longrightarrow \Gal(\wt F^{nr}_{\o\omega _{j,m_j}}/k)
\longrightarrow 1}.
\end{align*}
The global inertia subgroup
$I_{\wt F_{\omega _{j,m_j}}}$
\resp{$I_{\wt F_{\o\omega _{j,m_j}}}$} of order $N\centerdot m^{(j  )}$ can then be considered as the subgroup of inner automorphisms of Galois while the Galois subgroup
$\Gal(\wt F _{\omega _{j,m_j}}/k)$
\resp{$\Gal(\wt F _{\o\omega _{j,m_j}}/k)$} can be viewed as a subgroup of modular automorphisms of Galois with respect to
$I_{\wt F_{\omega _{j,m_j}}}$
\resp{$I_{\wt F_{\o\omega _{j,m_j}}}$}.
\Ei

\subsection{Shifted completions}

\Bi
\item In the context of this paper, we have to introduce  the shifted completions
$F^{S_{\cit}}_{\omega _{j,m_j}}=F_{\omega _{j,m_j}}\otimes\cit$
\resp{$F^{S_{\cit}}_{\o\omega _{j,m_j}}=F_{\o\omega _{j,m_j}}\otimes\cit$} where
$F_{\omega _{j,m_j}}$
\resp{$F_{\o\omega _{j,m_j}}$} denotes the corresponding unshifted \lr complex pseudo-ramified completion.

Notice that a set of shifted completions deals equivalently with
\Bean
\item a difference ring \cite{Coh} consisting in the ring of the unshifted completions and an isomorphism of this one onto the subring of shifted completions;
\item a $\GL_1(\cit)$-fibre bundle whose basis is the set of unshifted completions and total space the set of shifted completions: this case is the one-dimensional equivalent of the one envisaged in chapter 2 (for instance, see proposition 2.10).
\Ee

So, the unshifted completions $F_{\omega _{j,m_j}}$ are in one-to-one correspondence with the pseudo-ramified extensions
$\wt F_{\omega _{j,m_j}}$ of the symmetric splitting field
$\wt F=\wt F_R\cup \wt F_L$ of the polynomial ring
$k[x]$ while the shifted 
completions $F^{S_{\cit}}_{\omega _{j,m_j}}$ are in one-to-one correspondence with
pseudo-ramified extensions
$\wt F^{S_{\cit}}_{\omega _{j,m_j}}$ of the shifted symmetric splitting field $\wt F^{S_{\cit}}=\wt F^{S_{\cit}}_R\cup \wt F^{S_{\cit}}_L$  of the difference polynomial subring
$S(k[x])$.

\vskip 11pt

\item The {\bf sum}, over $j$, of the set of {\bbf equivalent complex pseudo-ramified shifted completions}
$F^{S_{\cit}}_{\omega _{j,m_j}}$
\resp{$F^{S_{\cit}}_{\o\omega _{j,m_j}}$}, $m_j\ge1$ is given by:
\[
F^{S_{\cit}}_{\omega _{\oplus}} = \bigoplus_j \bigoplus_{m_j} F^{S_{\cit}}_{\omega _{j,m_j}}\qquad
\rresp{F^{S_{\cit}}_{\o\omega _{\oplus}} = \bigoplus_j \bigoplus_{m_j} F^{S_{\cit}}_{\o\omega _{j,m_j}}}.\]
\Ei

\subsection{Weil shifted global bilinear (semi)groups}

\Bi
\item Let $\Gal(\wt F^{S_{\cit}}_{\omega _j}/k)$
\resp{$\Gal(\wt F^{S_{\cit}}_{\o\omega _j}/k)$} denote the Galois subgroup of the shifted extension
$\wt F^{S_{\cit}}_{\omega _j}$
\resp{$\wt F^{S_{\cit}}_{\o\omega _j}$}. Similarly,
$\Gal(\wt F^{nr;S_{\cit}}_{\omega _j}/k)$
\resp{$\Gal(\wt F^{nr;S_{\cit}}_{\o\omega _j}/k)$} will denote the Galois subgroup of the shifted pseudo-unramified extension
$\wt F^{nr;S_{\cit}}_{\omega _j}$
\resp{$\wt F^{nr;S_{\cit}}_{\o\omega _j}$}.
\vskip 11pt

\item If 
$I_{\wt F^{S_{\cit}}_{\omega _j}}$
\resp{$I_{\wt F^{S_{\cit}}_{\o\omega _j}}$} is the shifted global inertia subgroup of
$\Gal(\wt F^{S_{\cit}}_{\omega _j}/k)$
\resp{\linebreak $\Gal(\wt F^{S_{\cit}}_{\o\omega _j}/k)$}, then we have that:
\begin{align*}
\Gal(\wt F^{S_{\cit}}_{\omega _j}/k)&=
\Gal(\wt F^{nr;S_{\cit}}_{\omega _j}/k)\times
I_{\wt F^{S_{\cit}}_{\omega _j}}\\
\rresp{\Gal(\wt F^{S_{\cit}}_{\o\omega _j}/k)&=
\Gal(\wt F^{nr;S_{\cit}}_{\o\omega _j}/k)\times
I_{\wt F^{S_{\cit}}_{\o\omega _j}}}.\end{align*}
$I_{\wt F^{S_{\cit}}_{\omega _j}}$
\resp{$I_{\wt F^{S_{\cit}}_{\o\omega _j}}$} is the smallest normal subgroup (i.e. the subgroup of inner shifted automorphisms of Galois), of the subgroup 
$\Gal(\wt F^{S_{\cit}}_{\omega _j}/k)$
\resp{$\Gal(\wt F^{S_{\cit}}_{\o\omega _j}/k)$} of modular shifted automorphisms of Galois.
\vskip 11pt

\item If it is assumed that the global Weil group
$W^{ab}_{F^{S_{\cit}}_{\omega }}$
\resp{$W^{ab}_{F^{S_{\cit}}_{\o\omega }}$} is the Galois subgroup referring to pseudo-ramified extensions characterized by extension degrees $d=0\mod N$, then we have that:
\begin{align*}
W^{ab}_{F^{S_{\cit}}_{\omega }}
&\equiv \Gal(\dot{\wt F}^{S_{\cit}}_{\omega_\oplus }/k)
= \bigoplus_{j,m_j}\Gal(\dot{\wt F}^{S_{\cit}}_{\omega _{j,m_j}}/k)\\[11pt]
\rresp{W^{ab}_{F^{S_{\cit}}_{\o\omega }}
&\equiv \Gal(\dot{\wt F}^{S_{\cit}}_{\o\omega_\oplus }/k)
= \bigoplus_{j,m_j}\Gal(\dot{\wt F}^{S_{\cit}}_{\o\omega _{j,m_j}}/k)},\end{align*}
where
$\dot{\wt F}^{S_{\cit}}_{\omega _{j,m_j}}$
\resp{$\dot{\wt F}^{S_{\cit}}_{\o\omega _{j,m_j}}$} denote these shifted pseudo-ramified extensions with degrees $d=0\mod N$.

This leads to the product of the {\bf shifted global Weil groups} $W^{ab}_{F^{S_{\cit}}_{\o\omega _\oplus}}
\times W^{ab}_{F^{S_{\cit}}_{\omega_\oplus }}$
\[
W^{ab}_{F^{S_{\cit}}_{\o\omega }}
\times W^{ab}_{F^{S_{\cit}}_{\omega }}
= \Gal(\dot{\wt F}^{S_{\cit}}_{\o\omega_\oplus }/k)
\times \Gal(\dot{\wt F}^{S_{\cit}}_{\omega_\oplus }/k)
\subset \Gal({\wt F}^{S_{\cit}}_{\o\omega _\oplus }/k)
\times \Gal({\wt F}^{S_{\cit}}_{\omega _\oplus}/k)\;.\]
\Ei

\subsection{From abelian class field theory to its nonabelian equivalence}

The set of \lr pseudo-ramified extensions
$\wt F_{\omega _{j,m_j}}$
\resp{$\wt F_{\o\omega _{j,m_j}}$}, $1\le j\le r$, generates a one-dimensional complex affine semigroup 
$\SS^1_L$
\resp{$\SS^1_R$} in such a way the $n$-dimensional equivalent of their product 
$\SS^1_R\times \SS^1_L$ is a  complex bilinear algebraic semigroup
$G^{(2n)}(\wt F_{\o\omega }\times \wt F_\omega )$, isomorphic to the bilinear algebraic semigroup of matrices
\[ \GL_n(\wt F_{\o\omega }\times\wt F_\omega )
=T^t_n(\wt F_{\o\omega })\times  T_n(\wt F_\omega )\]
where:
\Bi
\item  $\wt F_\omega =\{\wt F_{\omega _1},\dots,\wt F_{\omega _{j,m_j}},\dots,\wt F_{\omega _{r,m_r}}\}$
\resp{$\wt F_{\o\omega} =\{\wt F_{\o\omega _1},\dots,\wt F_{\o\omega _{j,m_j}},\dots,\wt F_{\o\omega _{r,m_r}}\}$} denotes the set of complex pseudo-ramified finite extensions;
\item $T_n(\wt F_\omega )$ is the (semi)group of upper triangular matrices with entries in $\wt F_\omega $;
\item $T^t_n(\wt F_{\o\omega })$ is the (semi)group of lower triangular matrices with entries in $\wt F_{\o\omega} $.
\Ei
\vskip 11pt

\subsection{The algebraic general bilinear semigroup}

\Bi
\item Let $\wt B_{F_{\omega }}$
\resp{$\wt B_{F_{\omega }}$} be a \lr {\bbf division semialgebra of complex dimension $n$} over the set $\wt F_\omega $ \resp{$\wt F_{\o\omega }$} of the pseudo-ramified extensions
$\wt F_{\omega _{j,m_j}}$
\resp{$\wt F_{\o\omega _{j,m_j}}$} of $k$.

Then, $\wt B_{F_{\omega }}$
\resp{$\wt B_{F_{\o\omega }}$}, which is a \lr vector semispace of complex dimension $n$ over $\wt F_\omega $ \resp{$\wt F_{\o\omega} $}, is isomorphic to the algebra of Borel upper \resp{lower} triangular matrices:
\[
\wt B_{F_{\omega }}\approx T_n(\wt F_{\omega })
\rresp{\wt B_{F_{\o\omega }}\approx T^t_n(\wt F_{\o\omega })}\;.\]

\item This allows to define the {\bbf bilinear general semigroup} $\GL_n(\wt F_{\o\omega }\times \wt F_\omega )$ by:
\[
\wt B_{F_{\o\omega }}\times
\wt B_{F_{\o\omega }}
\simeq T_n^t(\wt F_{\o\omega })
\times T_n(\wt F_{\omega })\equiv \GL_n(\wt F_{\o\omega }\times \wt F_\omega )\]
such that its representation space is given by the tensor product of a right $\wt B_{F_{\o\omega }}$-semimodule $\wt M_R$ by a left $\wt B_{F_\omega }$-semimodule $\wt M_L$.
\vskip 11pt

\item Taking into account the definition of 
$\wt F_{\omega _{\oplus}}$
\resp{$\wt F_{\o\omega _{\oplus}}$}, the 
{\bbf $\wt B_{F_{\omega _{\oplus}}}$-semimodule $\wt M_{L_\oplus}$}
\resp{ $\wt B_{F_{\o\omega _{\oplus}}}$-semimodule $\wt M_{R_\oplus}$} decomposes according to:
\[ \wt M_{L_\oplus}=\bigoplus_{j}\bigoplus_{m_j}
 \wt M_{\omega _{j,m_j}}
 \rresp{\wt M_{R_\oplus}=\bigoplus_{j}\bigoplus_{m_j}
 \wt M_{\o\omega _{j,m_j}}}\]
 where 
\begin{align*}
 \wt M_{\omega _{j,m_j}}&\simeq t_n(\wt F_{\omega _{j,m_j}})\subset T_n(\wt F_\omega )\\
\rresp{\wt M_{\o\omega _{j,m_j}}&\simeq t^t_n(\wt F_{\o\omega _{j,m_j}})\subset T^t_n(\wt F_{\o\omega} )}
\end{align*}
is the representation subspace of 
$T_n(\wt F_\omega )$
\resp{$T^t_n(\wt F_{\o\omega} )$} restricted to the extension
$\wt F_{\omega _{j,m_j}}$
\resp{$\wt F_{\o\omega _{j,m_j}}$} and corresponds to the $m_j$-th representative of the $j$-th conjugacy class of
$\wt M_L$
\resp{$\wt M_R$}.
\vskip 11pt

\item Let $t_n(\wt F_{\omega _{j,m_j}})$
\resp{$t^t_n(\wt F_{\o\omega _{j,m_j}})$} be an element of 
$T_n(\wt F_\omega )$
\resp{$T^t_n(\wt F_{\o\omega })$} having the Levi decomposition:
\begin{align*}
t_n(\wt F_{\omega _{j,m_j}})&=d_n(\wt F_{\omega _{j,m_j}})\ u_n(\wt F_{\omega _{j,m_j}})\\
\rresp{t^t_n(\wt F_{\o\omega _{j,m_j}})&=u^t_n(\wt F_{\o\omega _{j,m_j}})\ d_n(\wt F_{\o\omega _{j,m_j}})}\end{align*}
where $d_n(\centerdot)$ is a diagonal matrix of order $n$ and where $u_n(\centerdot)$ is an upper unitriangular matrix.

So, any matrix $g_n(\wt F_{\o\omega _{j,m_j}}\times
\wt F_{\omega _{j,m_j}})\in \GL_n(\wt F_{\o\omega }\times \wt F_\omega )$ satisfies the {\bbf bilinear Gauss decomposition}:
\[
g_n(\wt F_{\o\omega _{j,m_j}}\times \wt F_{\omega _{j,m_j}})=
[(d_n(\wt F_{\omega _{j,m_j}})\times d_n(\wt F_{\o\omega _{j,m_j}}))]
[(u^t_n(\wt F_{\o\omega _{j,m_j}})\ u_n(\wt F_{\omega _{j,m_j}}))]\;.\]
\Ei

\subsection{Pseudo-ramified lattices}

Let $\Os_{\wt F_\omega }$ 
\resp{$\Os_{\wt F_{\o\omega }}$ } be the maximal order of $\wt F_\omega$  \resp{$\wt F_{\o\omega }$}.  Then,
$\Lambda _{\omega }=\Os_{\wt B_{F_{\omega }}}$
\resp{$\Lambda _{\o\omega }=\Os_{\wt B_{F_{\o\omega }}}$}
in the division semialgebra 
$\wt B_{F_{\omega }}$
\resp{$\wt B_{F_{\o\omega }}$} is a {\bbf pseudo-ramified $\zit/N\  \zit$-lattice} in the \lr 
$\wt B_{F_{\omega }}$-semimodule $\wt M_L$
\resp{$\wt B_{F_{\o\omega }}$-semimodule $\wt M_R$}.  So, we can fix the isomorphisms:
\[
\Lambda _\omega \simeq T_n(\Os_{\wt F_{\omega }})
\qquad \text{and} \qquad
\Lambda _{\o\omega} \simeq T^t_n(\Os_{\wt F_{\o\omega }})
\]
leading to $\Lambda _{\o\omega }\otimes\Lambda _\omega 
\simeq \GL_n(\Os_{\wt F_{\o\omega }}\times\Os_{\wt F_\omega })$.  And, if we take into account the decomposition of $\wt F_\omega $ and $\wt F_{\o\omega }$ into their pseudo-ramified extensions, we have that the sublattice
$\Lambda _{\omega _{j,m_j}}$
 \resp{$\Lambda _{\o\omega _{j,m_j}}$} into the 
$t_n(\wt F_{\omega _{j,m_j}})$-subsemimodule $\wt M_{\omega _{j,m_j}}$
\resp{$t^t_n(\wt F_{\o\omega _{j,m_j}})$-subsemimodule $\wt M_{\o\omega _{j,m_j}}$} verifies:
\[
\Lambda _{\omega _{j,m_j}}\simeq t_n(\Os_{\wt F_{\omega_{j,m_j}}})
\rresp{\Lambda _{\o\omega _{j,m_j}}\simeq t^t_n(\Os_{\wt F_{\o\omega_{j,m_j}}})}\]
and
\[
\Lambda _{\o\omega _{j,m_j}}\otimes
\Lambda _{\omega _{j,m_j}}\simeq
g_n ( \Os_{\wt F_{\o\omega _{,m_j}}}\times
\Os_{\wt F_{\omega _{,m_j}}})\in
\GL_n ( \Os_{\wt F_{\o\omega }}\times
\Os_{\wt F_{\omega }})\;.\]

\subsection{Proposition}

{\em Assume that we have fixed the isomorphism 
$\Lambda _{\o\omega }\otimes \Lambda _{\omega }\simeq \GL_n((\zit/N\  \zit)^2)$.

Then, the representation space 
$\Repsp (\GL_n((\zit/N\  \zit)^2)$ of
$\GL_n((\zit/N\  \zit)^2)$ decomposes according to:
\[
\Repsp (\GL_n((\zit/N\  \zit)^2)
= \bigoplus_j \bigoplus_{m_j}
(\Lambda _{\o\omega _{j,m_j}}\otimes
\Lambda _{\omega _{j,m_j}})\]
where the direct sums bear over the places of $F_\omega $ and $F_{\o\omega}$ having multiplicities $m^{(j)}=\sup(m_j)$.
}
\vskip 11pt

\subsection{Proposition}

{\em \Bena
\item The {\bbf pseudo-ramified Hecke bialgebra}  $\Hs\RL(n)$, generated by all the pseudo-ramified Hecke bioperators $T_R(n;t)\otimes T_L(n;t)$ has a representation in the arithmetic subgroup of matrices $\GL_n(\zit/N\  \zit)^2)$.

\item The {\bbf $j$-th coset representative of $T_R(n;t)\otimes T_L(n;t)$} is given by:
\[ U_{j_R}\times U_{j_L}
= [d_n(\Os_{\wt F_{\o\omega _{j,m_j}}})\centerdot d_n(\Os_{\wt F_{\omega _{j,m_j}}})]\times
[u^t_n(\Os_{\wt F_{\o\omega _{j,m_j}}})\centerdot u_n(\Os_{\wt F_{\omega _{j,m_j}}})]\;.\]

\item The Hecke bialgebra $\Hs\RL(n)$ generates the {\bbf endomorphisms of the pseudo-ramified $\wt B_{F_{\o\omega }}\otimes \wt B_{F_\omega }$-bisemimodule} 
$\wt M_{R_{\oplus}} \otimes \wt M_{L_{\oplus}}$ decomposing it according to the bisublattices
$(\Lambda _{\o\omega _{j,m_j}}\otimes\Lambda _{\omega _{j,m_j}})$ of
$\Lambda _{\o\omega }\otimes\Lambda _{\omega }$:
\[
\wt M_{R_\oplus}\otimes \wt M_{L_\oplus}
=\bigoplus_j\bigoplus_{m_j}
(\wt M_{\o\omega _{j,m_j}}\otimes
\wt M_{\omega _{j,m_j}})\]
where $(\wt M_{\o\omega _{j,m_j}}\otimes
\wt M_{\omega _{j,m_j}})$ is a 
$(\wt B_{F_{\o\omega _j}}\times \wt B_{F_{\omega _j}})$-bisubsemimodule representative.
\Ee}
\vskip 11pt

\begin{proof} These assertions were proved in \cite{Pie3}.
\end{proof}
\vskip 11pt

\subsection{Corollary}

{\em There exists an injective morphism:
\[
J_{\Lambda \to M}: \quad
\Lambda _{\o\omega }\otimes
\Lambda _{\omega }\To \wt M_R\otimes \wt M_L\]
from the bilattice $\Lambda _{\o\omega }\otimes \Lambda _{\omega }$ into the $\GL_n(\wt F_{\o\omega }\times \wt F_\omega )$-bisemimodule $\wt M_R\otimes \wt M_L$.
}
\vskip 11pt

\subsection{Toroidal compactification}

\Bi
\item Let $\GL_n(\wt F_R\times \wt F_L)$ be the bilinear algebraic semigroup over the product of symmetric splitting semifields $\wt F_R$ and $\wt F_L$.

Let $Y_{S\RL}=\GL_n(\wt F_R\times \wt F_L)
\big/\GL_n((\zit/N\  \zit)^2)$ be the non compact pseudo-ramified lattice bisemispace.
\vskip 11pt

\item The {\bbf Borel-Serre toroidal compactification of $Y_{S\RL}$} is a toroidal projective emergent isomorphism of compactification given by:
\[ \gamma ^c\RL : \qquad Y_{S\RL} \To \o Y_{S^T\RL}\]
where:
\Bi
\item $ \o Y_{S^T\RL}=\GL_n(F^T_R\times F^T_L) \big/
\GL_n((\zit/N\ \zit)^2)$;

\item $F^T_R$ and $F^T_L$ are toroidal compactifications of $\wt F_R$ and $\wt F_L$ respectively;
\Ei
such that:
\Bi
\item $Y_{S\RL}$ may be viewed as the interior of
$\o Y_{S^T\RL}$ in the sense that the isomorphism
 $\gamma ^c\RL$ is an inclusion isomorphism
$Y_{S\RL}\hookrightarrow 
\o Y_{S^T\RL}$ given by a homotopy equivalence;

\item $\o Y_{S^T\RL}$ is a $\GL_n(F^T_{\o\omega }\times F^T_\omega )$-bisemimodule $M^T_R\otimes M^T_L$ over the sets
$F^T_\omega =\linebreak \{F^T_{\omega _1},\dots,F^T_{\omega _{r,mr}}\}$ and
$F^T_{\o\omega} =\{F^T_{\o\omega _1},\dots,F^T_{\o\omega _{r,mr}}\}$ of toroidal completions.
\Ei

By this way, $\gamma ^c\RL$ sends all equivalent representatives of conjugacy classes of\linebreak $\GL_n(\wt F_R\times \wt F_L)$ into their toroidal compactified equivalents which are products of $n$-dimensional complex semitori  $T^{2n}_R[j,m_j]\times T^{2n}_L[j,m_j]$.
\vskip 11pt

\item On the other hand, let 
$F^T_{\omega ^1}=\{F^T_{\omega ^1_1},\dots,F^T_{\omega _{j,mr}}\}$
\resp{$F^T_{\o\omega ^1}=\{F^T_{\o\omega ^1_1},\dots,F^T_{\o\omega _{j,mr}}\}$} denote the set of irreducible toroidal completions.

The {\bbf bilinear complex parabolic semigroup $P_n(F^T_{\o\omega ^1_1}\times F^T_{\omega ^1})$} is the smallest normal bilinear subsemigroup of
$\GL_n(F^T_{\o\omega }\times F^T_\omega )$, representing the $n$-fold product of the global inertia subgroup $I_{\wt F_{\o\omega }}\times I_{\wt F_\omega }$.

The {\bbf double coset decomposition} of $\GL_n(F^T_R\times F^T_L)$ gives rise to the compactified bisemispace \cite{Vog}:
\[ 
S^{P_n}_{\GL_n}=P_n(F^T_{\o\omega ^1}\times F^T_{\omega ^1})\setminus \GL_n(F^T_R\times F^T_L)/\GL_n((\zit/N\  \zit)^2)\;.\]
\Ei

\subsection{Proposition}

{\em As a consequence of the double coset decomposition of the compactified bisemivariety $\o S^{P_n}_{K_n}$, the {\bbf modular conjugacy classes of $\GL_n(F^T_{\o\omega }\times F^T_\omega )$} with respect to the bilinear parabolic semigroup $P_n(F^T_{\o\omega ^1}\times F^T_{\omega ^1})$ correspond to the cosets of the compactified pseudo-ramified lattice bisemispace
$\o Y_{S^T\RL}=\GL_n(F^T_R\times F^T_L)/\GL_n((\zit/N\  \zit)^2)$.
}
\vskip 11pt

\begin{proof} As the bilinear parabolic semigroup
$P_n(F^T_{\o\omega ^1}\times F^T_{\omega^1} )$
is compact and as the cosets of
$\o Y_{S^T\RL}$ correspond to the set of lattices of
$((F^T_{\o\omega })^n\times(F^T_{\omega })^n)$, we have that:
\[ P_n(F^T_{\o\omega ^1}\times F^T_{\omega^1} )/
\GL_n(F^T_R\times F^T_L)\approx
\GL_n(F^T_R\times F^T_L)/\GL_n((\zit/N\  \zit)^2)
\]
implying that the modular conjugacy classes of
$\GL_n(F^T_{\o\omega }\times F^T_{\omega} )$ are the cosets of the bilinear quotient semigroup
$P_n(F^T_{\o\omega ^1}\times F^T_{\omega^1} )/\GL_n(F^T_R\times F^T_L)$.
\end{proof}
\vskip 11pt

\subsection{Reducible Galois cohomologies}

\Bi
\item Let $n=n_1+\dots+n_s$ be a partition of $n$ \cite{Rod}, \cite{Zel} and let
\begin{align*}
\o Y^{2n=2n_1+\dots+2n_s}_{S\RL}
&= \GL_{n=n_1+\dots+n_s}(F_R\times F_L)/\GL_n((\zit/N\  \zit)^2)\\
&=\Repsp ( \GL_{n=n_1+\dots+n_s}(F_{\o\omega }\times F_\omega ))
=\mathop{\boxplus}\limits_{n_\ell }\Repsp ( \GL_{n_\ell }(F_{\o\omega }\times F_\omega ))
\end{align*}
be the {\bbf reducible compactified representation space of}
$\GL_n(\wt F_{\o\omega }\times \wt F_\omega )$ decomposing according to the irreducible representation spaces
$\Repsp(\GL_{n_\ell }(F_{\o\omega }\times F_\omega ))$ of
$\GL_{n }(F_{\o\omega }\times F_\omega )$ given with respect to modular conjugacy classes ``$j$''.
\vskip 11pt

\item {\bbf The bilinear cohomology of} $\o Y^{2n=2n_1+\dots+2n_s}_{S\RL}$, 
introduced in \cite{Pie3}, (section 3.2),
decomposes according to
\[ H^*(\o Y^{2n=2n_1+\dots+2n_s}_{S\RL},
M^{2n}_{R}\otimes M^{2n}_{L})
=\bigoplus^{2n_s}_{2n_\ell =2n_1} H^{2n_\ell  }
(\o Y^{2n}_{S\RL},M^{2n_\ell }_R\otimes M^{2n_\ell }_L)\]
where $M^{2n_\ell }_L$ 
\resp{$M^{2n_\ell }_R$} is a \lr 
$T_{n_\ell }(F_\omega )$-subsemimodule
\resp{$T^t_{n_\ell }(F_{\o\omega })$-subsemimodule}
of real dimension $2n_\ell $.

As $M^{2n}\otimes M^{2n}\equiv G^{(2n)}(F_{\o\omega }\times F_\omega )$ is a smooth abstract bisemivariety, the bilinear cohomology on its algebraic equivalent
$\wt M^{2n}\otimes \wt M^{2n}\equiv G^{(2n)}(\wt F_{\o\omega }\times \wt F_\omega )$ will be similarly decomposed into
\[H^*(G^{(2n)}(\wt F_{\o\omega }\times \wt F_\omega ))=
\bigoplus_{2n_\ell}H^{2n_\ell}(G^{(2n)}(\wt F_{\o\omega }\times \wt F_\omega ),
\wt M_R^{2n_\ell}\otimes \wt M_L^{2n_\ell})\;.\]
\vskip 11pt

\item The cohomology of the reducible toroidal bisemivariety $\o S^{P_n}_{\GL_n}$ also decomposes according to:
\[
H^*(\o S^{P_n=n_1+\dots+n_s}_{\GL_n=n_1+\dots+n_s},
M^{2n}_{T_{R }}\otimes M^{2n}_{T_{L }})
=\bigoplus_{2n_\ell } H^{2n_\ell  }(\o S^{P_n}_{\GL_n},M^{2n_\ell }_{T_R}\otimes M^{2n_\ell }_{T_L})
\]
where $M^{2n_\ell }_{T_L}$
\resp{$M^{2n_\ell }_{T_R}$} is a \lr compactified 
$T_{n_\ell }(F^T_\omega )$
\resp{\linebreak $T^t_{n_\ell }(F^T_{\o\omega })$}-subsemimodule of dimension $2n_\ell $.
\vskip 11pt

\item However, {\bbf the coefficients of the cohomology} are generally considered in (bisemi)\-sheaves of rings over bilinear complete algebraic semigroups
$(M^{2n_\ell }_R\otimes M^{2n_\ell }_L)$.

In this purpose, {\bbf a (semi)sheaf 
$\widehat M^{2n_\ell }_L$
\resp{$\widehat M^{2n_\ell }_R$} of $C^\infty  $-differentiable functions on 
$M^{2n_\ell }_L$
\resp{$M^{2n_\ell }_R$}} will be envisaged and a (bisemi)sheaf\linebreak
$(\widehat M^{2n_\ell }_R\otimes \widehat M^{2n_\ell }_L)$ of $C^\infty $-differentiable bifunctions (i.e. products of cofunctions by functions) on
$(M^{2n_\ell }_R\otimes M^{2n_\ell }_L)$ will be considered as coefficients of the cohomology
$H^{2n_\ell  }(\o Y^{2n}_{S\RL}, \widehat M^{2n_\ell }_R\otimes \widehat M^{2n_\ell }_L)$.
\vskip 11pt


\Ei

\subsection{Algebraic semicycles of the Chow (semi)groups}

Let $Y_L$ \resp{$Y_R$} denote a \lr algebraic semigroup
$G^{(2n)}(F_\omega )$ \resp{$G^{(2n)}(F_{\o\omega })$} of complex dimension $n$ isomorphic to a \lr smoth semischeme.  Then, the
algebraic semicycle 
${\rm CY}^{2n_\ell }(Y_L)$
\resp{${\rm CY}^{2n_\ell }(Y_R)$}
of dimension $2n_\ell $  on 
$Y_L$
\resp{$Y_R$} is such that:
\begin{align*}
&{\rm CY}^{2n_\ell }(Y_L) 
\subset \bZ^{2n_\ell }(Y_L)
\subset {\rm CH}^{2n_\ell }(Y_L) 
\\
\rresp{&{\rm CY}^{2n_\ell }(Y_R) 
\subset \bZ^{2n_\ell }(Y_R)
\subset {\rm CH}^{2n_\ell }(Y_R) 
}
\end{align*}
where:
\Bi
\item $\bZ^{2n_\ell }(Y_L)$ is the {\bbf semigroup of algebraic semicycles} ${\rm CY}^{2n_\ell }(Y_L)$ of codimension $2n_\ell $;

\item ${\rm CH}^{2n_\ell }(Y_L)=\dfrac{\bZ^{2n_\ell }(Y_L)}{\bZ^{2n_\ell }_{\rm rat}(Y_L)}$ is the {\bbf $2n_\ell $-th Chow semigroup of $Y_L$} with $\bZ^{2n_\ell }_{\rm rat}(Y_L)$ the semigroup of algebraic semicycles of codimension $2n_\ell $ rationally equivalent to zero \cite{Mur}, \cite{Vis}.
\Ei
It is evident that ${\rm CY}^{2n_\ell }(Y_L)$
\resp{${\rm CY}^{2n_\ell }(Y_R)$} decomposes according to the equivalence classes ``$j$'' having representatives $m_j$ such that:
\[
{\rm CY}^{2n_\ell }(Y_R)\times {\rm CY}^{2n_\ell }(Y_L)=
\bigoplus_j
\bigoplus_{m_j}
({\rm CY}^{2n_\ell }(Y_R[j,m_j])\times {\rm CY}^{2n_\ell }(Y_L[j,m_j])\;.\]
\vskip 11pt

\subsection{\SV\ motivic presheaf \cite{Frie}}

\Bi
\item Let $\Sigma ^{2n_\ell }_L$
\resp{$\Sigma ^{2n_\ell }_R$} denote a \lr complex topological $2n_\ell $-simplex and let
$\Sigma ^\centerdot_L$
\resp{$\Sigma ^\centerdot_R$} denote a cosimplicial object from the collection of the 
$\Sigma ^{2n_\ell }_L$
\resp{$\Sigma ^{2n_\ell }_R$} in the category 
$Sm_L(k)$
\resp{$Sm_R(k)$} of \lr  (semi)schemes over $k $.
\vskip 11pt

\item A {\bbf \SV\ motivic \lr presheaf} of the left (resp.\linebreak right)  (semi)scheme 
$X^{\rm sv}_L$
\resp{$X^{\rm sv}_R$} of complex dimension $\ell $ on 
$Sm_L(k )$
\resp{$Sm_R(k )$} and denoted 
$\u c_*(X^{\rm sv}_L)$
\resp{$\u c_*(X^{\rm sv}_R)$}
is a functor from 
$X^{\rm sv}_L$
\resp{$X^{\rm sv}_R$} to the \lr chain complex associated to the abelian semigroup\linebreak
$\mathop{\sqcup}\limits_{i_\ell } \Hom_{Sm_L(k )}(\dot\Sigma _L,SP^{i_\ell }(X^{\rm sv}_L))$
\resp{$\mathop{\sqcup}\limits_{i_\ell } \Hom_{Sm_R(k )}(\dot\Sigma _R,SP^{i_\ell }(X^{\rm sv}_R))$} where
$SP^{i_\ell }(X^{\rm sv}_L)$
\resp{$SP^{i_\ell }(X^{\rm sv}_R)$}
denotes the $i_\ell$-th symmetric product of
$X^{\rm sv}_L$
\resp{$X^{\rm sv}_R$}.
\vskip 11pt

\item On the other hand, let {\bbf 
$Z_L(2n_\ell )$
\resp{$Z_R(2n_\ell )$}} denote the \lr
{\bbf \SV\ submotive of dimension $2n_\ell =i_\ell \times2\ell $}
as developed in chapter 1 of \cite{Pie3}.  
$Z_L(2n_\ell )$
\resp{$Z_R(2n_\ell )$} can be checked to correspond to a \lr element of the $2n_\ell $-th semigroup
$\bZ^{2n_\ell }(Y _{L})$
\resp{$\bZ^{2n_\ell }(Y _{R })$} of \lr algebraic semicycles over $Y_L$ \resp{$Y_R$} of dimension $2n$.
\vskip 11pt

\item Similarly, let $Z^T_L(2n_\ell )$
\resp{$Z^T_R(2n_\ell )$} be the resulting toroidal compactified \SV\ submotive obtained from $Z_L(2n_\ell )$
\resp{$Z_R(2n_\ell )$}.
\Ei

\subsection{Proposition}

{\em
The cohomologies $H^{2n_\ell }(Y\RL,{\rm CY}^{2n_\ell }(Y_R) \times {\rm CY}^{2n_\ell }(Y_L))$ and
$H^{2n_\ell } (\u c_*(X^{\rm sv}\RL),
Z\RL({2n_\ell }))$
are {\bfseries bilinear pure motivic cohomologies}, with
$Y\RL=Y_R\times Y_L$ and 
$ \u c_*(X^{\rm sv}\RL)=
\u c_*(X^{\rm sv}_R)\times
\u c_*(X^{\rm sv}_L)$.
}
\vskip 11pt

\begin{proof}
Let $Z\RL(2n_\ell  )\equiv Z_R(2n_\ell )\otimes Z_L(2n_\ell )$ denote the bilinear products \cite{Pie5}, right by left, of \SV\  submotives of complex codimension $n_\ell $.

If we have the isomorphisms:
\begin{align*}
i_{M-X}&: \quad 
H^{2n_\ell } 
(Y\RL,
{\rm CY}^{2n_\ell }(Y_R)\times
{\rm CY}^{2n_\ell }(Y_L))\\
& \qquad \qquad
\overset{\sim}{\To}
H^{2n_\ell } (\u c_*(X^{\rm sv}\RL),
Z\RL({2n_\ell  })))
\end{align*}
resulting from sections 1.15 and 1.16., and \cite{Pie3},



it is evident that the cohomologies
\[
H^{2n_\ell  }(Y\RL,{\rm CY}^{2n_\ell }(Y_R)
\times {\rm CY}^{2n_\ell }(Y_L))=
\Hom_{C_{M\RL}}(Y\RL,{\rm CY}^{2n_\ell }(Y_R)
\times {\rm CY}^{2n_\ell }(Y_L))\]
and $H^{2n_\ell } (\u c_*(X^{\rm sv}_{T\RL},
Z\RL({2n_\ell  }))$
are ``pure'' motivic, noticing that $C_{M\RL}$  is the category of smooth bisemischemes isomorphic to $\GL_{n_\ell}(F_{\o\omega }\times F_\omega )$-bisemimodules $M^{2n_\ell}_R\otimes M^{2n_\ell}_L$ .
\end{proof}

%% file: LanglandsIII_2.tex
\section{Bilinear cohomologies of mixed (bisemi)motives}

The objective consists now in introducing a left and a right triangulated category 
$DM_L(k)$ and
$DM_R(k)$ of mixed (semi)motives \cite{Del2}, \cite{Del3}, \cite{Jan2}  and in developing a corresponding suitable bilinear mixed motivic cohomology.
\vskip 11pt

\subsection{Definition: Correspondences on \SV\ (semi)\-motives}

The \SV\ \lr semimotive
$\u c_*(X^{\rm sv}_L)$
\resp{$\u c_*(X^{\rm sv}_R)$},
also noted
$M(X^{\rm sv}_L)$
\resp{$M(X^{\rm sv}_R)$}
 has the property to be a \lr {\bbf presheaf with transfers} \cite{Mor}, \cite{Frie}.  That is to say that there exist {\bbf \lr correspondences}, noted 
$\Corr(SP^{f_\ell }(X^{\rm sv}_L),X^{2n_\ell -2\fl}_L)$
\resp{$\Corr(SP^{f_\ell }(X^{\rm sv}_R),X^{2n_\ell -2\fl}_R)$}, on the set of irreducible subvarieties of 
$X^{2n_\ell }_L$
\resp{$X^{2n_\ell }_R$}.

Left \resp{right} correspondences are here introduced by:
\begin{align*}
\Corr ( SP^{f_\ell } ( X^{\rm sv}_L),X^{2n_\ell -2\fl}_L )
&: \;
SP^{i_\ell }(X^{\rm sv}_L)\to
X^{2n_\ell }_L=SP^{f_\ell }(X^{\rm sv}_L)\times X^{2n_\ell -2\fl}_L\\[11pt]
\rresp{\Corr ( SP^{f_\ell }(X^{\rm sv}_R),X^{2n_\ell -2\fl}_R )
&: \;
SP^{i_\ell }(X^{\rm sv}_R)\to
X^{2n_\ell }_R=SP^{f_\ell }(X^{\rm sv}_R)\times X^{2n_\ell -2\fl}_R},
\end{align*}
for the integers
\Bi
\item $\fl\le n_\ell \le n$;
\item $i_\ell \times \ell =n_\ell =(\fl)+(n_\ell -\fl)$;
\item $f_\ell \le i_\ell $;
\Ei
such that:
\Bi
\item the $i_\ell $-th sub(semi)motive $SP^{i_\ell }(X^{\rm sv}_L)$ of dimension $2n_\ell =i_\ell \times 2\ell $ be sent by the left correspondence $\Corr(\centerdot,\centerdot)$ to the product $X^{2n_\ell }_L$ of closed irreducible sub(semi)motives
$SP^{f_\ell }(X^{\rm sv}_L)$ by $X^{2n_\ell -2\fl}_L$ where $X^{2n_\ell -2\fl}_L$ is a smooth presheaf of complex dimension $n_\ell -\fl$;

\item there exists a projection from
$X^{2n_\ell }_L$
\resp{$X^{2n_\ell }_R$} to an irreducible component of
$SP^{f_\ell }(X^{\rm sv}_L)$
\resp{$SP^{f_\ell }(X^{\rm sv}_R)$}.
\Ei

\subsection[Definition: Fibre of the tangent bundle $\Tan(SP^{f_\ell }(X^{\rm sv}_{L,R}))$]{\bbf Definition: Fibre of the tangent bundle $\Tan(SP^{f_\ell }(X^{\rm sv}_{L,R}))$}

Let $\TAN[SP^{f_\ell }(X^{\rm sv}_L)]$
\resp{$\TAN[SP^{f_\ell }(X^{\rm sv}_R)]$} be the \lr tangent vector bundle given by the triple:
\begin{align*}
 &\Tan[SP^{f_\ell }(X^{\rm sv}_L)](\Delta ^{2\fl}_L,\ppr_L,SP^{f_\ell }(X^{\rm sv}_L))\\
\rresp{&\Tan[SP^{f_\ell }(X^{\rm sv}_R)](\Delta ^{2\fl}_R,\ppr_R,SP^{f_\ell }(X^{\rm sv}_R))}\end{align*}
where:
\Bi
\item $\Delta ^{2\fl}_L$
\resp{$\Delta ^{2\fl}_R$} is the total space obtained from the base space 
$SP^{f_\ell }(X^{\rm sv}_L)$
\resp{$SP^{f_\ell }(X^{\rm sv}_R)$} under the action of the upper \resp{lower} linear trigonal group
$T_{\fl}(\cit)$
\resp{$T^t_{\fl}(\cit)$} $\subset \GL_{\fl}(\cit\times \cit)$ such that
\begin{align*}
\Delta ^{2\fl}_L &= SP^{f_\ell } (X^{\rm sv}_L)\times \AdFRepsp ( T_{\fl}( \cit ))\\[11pt]
\rresp{\Delta ^{2\fl}_R &= SP^{f_\ell } (X^{\rm sv}_R)\times \AdFRepsp ( T_{\fl}( \cit ))}
\end{align*}
be defined with respect to the {\bbf \lr fibre 
$\AdFRepsp (T_{\fl}(\cit))$
\resp{$\AdFRepsp (T^t_{\fl}(\cit))$} which is given by the adjoint functional representation space of 
$T_{\fl}(\cit)$
\resp{$T^t_{\fl}(\cit)$};}

\item $\ppr_L$
\resp{$\ppr_R$} is the evident projection:
\[ \ppr_L: \quad \Delta ^{2\fl}_L\To SP^{f_\ell }(X^{\rm sv}_L)
\rresp{\ppr_R: \quad \Delta ^{2\fl}_R\To SP^{f_\ell }(X^{\rm sv}_R)}\;.\]
\Ei

\subsection{Definition: Shifted correspondences}
Taking into account the \lr tangent bundle
$\TAN[SP^{f_\ell }(X^{\rm sv}_L)]$
(resp.\linebreak {$\TAN[SP^{f_\ell }(X^{\rm sv}_R)]$})
as introduced in definition 2.2. and the \lr correspondences
$\Corr(SP^{f_\ell }(X^{\rm sv}_L),X^{2n_\ell -2\fl}_L)$
\resp{$\Corr(SP^{f_\ell }(X^{\rm sv}_R),X^{2n_\ell -2\fl}_R)$} on \SV\ semimotives, {\bbf shifted \lr correspondences can be defined by the homomorphism}:
\begin{align*}
\CORR^S_L: \quad
\Corr(SP^{f_\ell }(X^{\rm sv}_L),X^{2n_\ell -2\fl}_L)
&\To
\Corr^S_L(\Delta ^{2\fl}_L,X^{2n_\ell -2\fl}_L)\\[11pt]
\rresp{\CORR^S_L: \quad
\Corr(SP^{f_\ell }(X^{\rm sv}_R),X^{2n_\ell -2\fl}_R)
&\To
\Corr^S_R(\Delta ^{2\fl}_L,X^{2n_\ell -2\fl}_R)}
\end{align*}
where the \lr smooth presheaf 
$SP^{f_\ell }(X^{\rm sv}_L)$
\resp{$SP^{f_\ell }(X^{\rm sv}_R)$} has been sent to the corresponding smooth presheaf
\begin{align*}
 \Delta ^{2\fl}_L &= SP^{f_\ell }(X^{\rm sv}_L)
 \times \AdFRepsp (T_{\fl}(\cit ))
\\[11pt]
\rresp{ 2\Delta ^{\fl}_R &= SP^{f_\ell }(X^{\rm sv}_R)
 \times \AdFRepsp (T^t_{\fl}(\cit ))}
 \end{align*}
 by means of the inverse projection map
 $\ppr^{-1}_L$
 \resp{$\ppr^{-1}_R$} of
$\Tan[SP^{f_\ell }(X^{\rm sv}_L)]$
\resp{$\Tan[SP^{f_\ell }(X^{\rm sv}_R)]$}.
\vskip 11pt

\subsection{Triangulated category of mixed (semi)motives}

\Bi
\item Let the \SV\ \lr pure (semi)motive
$M(X^{\rm sv}_L)$
\resp{\linebreak $M(X^{\rm sv}_R)$}, provided with \lr shifted correspondences
$\Corr^S(\centerdot,\centerdot)$, be noted
$M_{DM_L}(X^{\rm sv}_L)$
\resp{$M_{DM_R}(X^{\rm sv}_R)$}: it is then a {\bbf \lr mixed (semi)motive} of the triangulated category
$DM_L(k )$
\resp{$DM_R(k)$} of \lr geometric (semi)motives. Indeed, the isomorphism:
\begin{multline*}
M_{\Corr_L}: \qquad
M(X^{\rm sv}_L)=\mathop{\sqcup}\limits_{i_\ell } \Hom_{Sm_L(k )}
(\Sigma ^\centerdot_L,SP^{i_\ell }(X^{\rm sv}_L))\\
\To M_{DM_L}(X^{\rm sv}_L)
= \mathop{\sqcup}\limits_{i_\ell ,f_\ell } \Hom_{Sm_L(k)}
 ( SP^{i_\ell }(X^{\rm sv}_L),X^{2n_\ell }_L[2\fl] ),
 \end{multline*}
 where $X^{2n_\ell }_L[\fl]=\Delta ^{2\fl}_L\times X^{2n_\ell -2\fl}$, maps the \SV\ pure (semi)motive 
$M(X^{\rm sv}_L)$ to the \SV\ mixed (semi)motive
$M_{DM_L}(X^{\rm sv}_L)$ by means of the left shifted correspondence
$\Corr^S(\Delta ^{2\fl}_L,X^{2n_\ell -2\fl}_L)$, taking into account that $\Delta ^{2\fl}_L$ {\bbf is a sub(semi)motive shifted in $2\fl$-dimensions}.
\vskip 11pt

\item Noticing that a triangulated category is an additive category graded by a translation functor and a set of distinguished triangles \cite{Ver}, we have that the isomorphism $M_{\Corr_L}$ can be viewed as belonging to the translation functor from the category of \SV\ pure (semi)motives to the {\bbf triangulated category $DM_L(k )$ of mixed (semi)motives} \cite{Hub}, \cite{C-F}.  

And, the derived category 
$D(M(X^{\rm sv}_L))$
\resp{$D(M(X^{\rm sv}_R))$} of pure \lr (semi)motives
$M(X^{\rm sv}_L)$
\resp{$M(X^{\rm sv}_R)$} with transfers is included into the corresponding triangulated category 
$DM_L(k)$
\resp{$DM_R(k)$}, {\bbf a derived category} resulting from a corresponding triangulated category with a condition of null homotopy on the automorphisms of translations \cite{F-S-V}.
\vskip 11pt

\item {\bf Remark} finally that
a triangulated category of mixed (semi)motives can also be defined from the toroidal  pure (semi)motives
$\u c_*(X^{\rm sv}_{T_L})$
\resp{$\u c_*(X^{\rm sv}_{T_R})$}: so, {\bbf the \SV\ \lr mixed (semi)motives} 
$M_{DM_L}(X^{\rm sv}_{T_L})$
\resp{$M_{DM_R}(X^{\rm sv}_{T_R})$}
belong to the \lr  derived category
$D(M(X^{\rm sv}_{T_L}))$
\resp{$D(M(X^{\rm sv}_{T_R}))$}.
\Ei

\subsection{Lemma}
{\em
Let
\begin{align*}
\Delta ^{2\fl}_{L} &= SP^{f_\ell } ( X^{\rm sv}_{L})
\times \AdFRepsp ( T_{\fl} ( \cit ))\\[11pt]
\rresp{\Delta ^{2\fl}_{R} &= SP^{f_\ell } ( X^{\rm sv}_{R})
\times \AdFRepsp ( T^t_{\fl} ( \cit ))}
\end{align*}
be a \lr $2\fl$-(semi)scheme in the category 
$Sm_L(k )$
\resp{$Sm_R(k )$} of \lr smooth (semi)schemes over
$k$.

Then, we have that the vector bisemispace 
$\Aut(\Tan_{e}(SP^{f_\ell } (X^{\rm sv}_R)\times
SP^{f_\ell } (X^{\rm sv}_L)))$ of the endomorphisms of the tangent bisemispace of
$SP^{f_\ell } ( X^{\rm sv}_{R} ) \times
SP(X^{\rm sv}_{L})$ is precisely a $(\fl\times\fl)$-subbisemischeme of
$\Delta ^{\fl}_{R}\times \Delta ^{\fl}_{L}$ such that:
\begin{multline*}
\Delta ^{2\fl}_{R}\times \Delta ^{2\fl}_{L}
\simeq \FRepsp(\GL_{\fl}((F_{\o\omega }\otimes\cit)\times
(F_{\omega }\otimes\cit )))\\
= \AdFRepsp ( \GL_{\fl}(\cit\times\cit ) )\times
\FRepsp ( \GL_{\fl}(F_{\o\omega }\times F_\omega  ) )
\end{multline*}
where:
\Bi
\item $SP^{f_\ell }(X^{\rm sv}_{R})\times
SP^{f_\ell }(X^{\rm sv}_{L})
\simeq \FRepsp (\GL_{\fl}(F_{\o\omega }\times _\omega ))$
is the functional representation space of $\GL_{\fl}(F_{\o\omega }\times F_\omega )$;

\item $\AdFRepsp(\GL_{\fl}(\cit\times\cit))$ is the bilinear fibre of the tangent bibundle\linebreak $\TAN[\SP^{f_\ell }(X^{\rm sv}_{R})\times
\SP^{f_\ell }(X^{\rm sv}_{L})]$ introduced in definition 2.2.\Ei
}
\vskip 11pt

\begin{proof}
As we are concerned with mixed bimotives of the product
$DM_L(k)\times DM_R(k)$ of triangulated categories, where a triangulated category is an additive category graded by a translation functor, we have that the total space
$(\Delta ^{2\fl}_{T_R}\times
\Delta ^{2\fl}_{T_L})$ of the tangent bibundle
$\TAN[SP^{f_\ell }(X^{\rm sv}_{T_R}))]\times
\TAN[SP^{f_\ell }(X^{\rm sv}_{T_L}))]$, introduced in definition 2.2, is the tangent bisemispace of
$(SP^{f_\ell }(X^{\rm sv}_{T_R})\times
SP^{f_\ell }(X^{\rm sv}_{T_L}))$ generated under the action of the Lie algebra of
$\GL_{\fl}(\cit\times\cit)$.

Let $\Tan_e(SP^{f_\ell }(X^{\rm sv}_{T_R})\times
SP^{f_\ell }(X^{\rm sv}_{T_L}))$ denote this tangent bisemispace at the identity element ``$e$'' in order to define differentials on it.  Then,
$\Aut(\Tan_e(SP^{f_\ell }(X^{\rm sv}_{T_R})\times
SP^{f_\ell }(X^{\rm sv}_{T_L})))$ is an open subset of the bilinear vector semispace of endomorphisms of
$\Tan_e(SP^{f_\ell }(X^{\rm sv}_{T_R})\times
SP^{f_\ell }(X^{\rm sv}_{T_L}))$ \cite{F-H}.  So, we have that:
\begin{multline*}
\Aut ( \Tan_e ( SP^{f_\ell } ( X^{\rm sv}_{T_R})\times
SP^{f_\ell }(X^{\rm sv}_{T_L}) ) )
\subset \Delta ^{2\fl}_{T_R}\times \Delta ^{2\fl}_{T_L}\\
\begin{aligned}
&\simeq \AdFRepsp ( \GL_{\fl}(\cit\times\cit ) )\times
 \FRepsp ( \GL_{\fl}(F_{\o\omega }\times F_\omega ) )\\
&=\FRepsp ( \GL_{\fl}((F_{\o\omega }\otimes\cit)\times (F_\omega \otimes\cit ) )))\;.\qedhere
\end{aligned}
\end{multline*}
\end{proof}

\subsection[Weil algebra of the adjoint representation of $\GL_{\fl}(\cit\times\cit)$]{\bbf Weil algebra of the adjoint representation of $\GL_{\fl}(\cit\times\cit)$}

Let $\TAN[SP^{f_\ell }(X^{\rm sv}_{R})\times
SP^{f_\ell }(X^{\rm sv}_{L})]$ be the tangent bibundle having as bilinear fibre the adjoint functional representation space of $\GL_{\fl}(\cit\times\cit)$ given by:
\[ \AdFRepsp ( \GL_{\fl} ( \cit\times\cit ))
\simeq (\Delta ^{2\fl}_{R}\times\Delta ^{2\fl}_{L})
\big/SP^{f_\ell }(X^{\rm sv}_{R})\times
SP^{f_\ell }(X^{\rm sv}_{L})\]
and denoted $\Fs^{2\fl}\RL ( \TAN)$.

The Lie algebra of $\Fs^{2\fl}\RL ( \TAN)$ is denoted 
${\rm Lie}(\Fs^{2\fl}\RL ( \TAN))$.

Let $A({\rm Lie}(\Fs^{2\fl}\RL ( \TAN)))$ be the exterior algebra of products, right by left, of differential forms of all degrees on 
${\rm Lie}(\Fs^{2\fl}\RL ( \TAN))$ and let
$S({\rm Lie}(\Fs^{2\fl}\RL ( \TAN)))$ denote the symmetric bialgebra corresponding to the symmetric multilinear forms on 
${\rm Lie}(\Fs^{2\fl}\RL ( \TAN))$.  

Then, the Weil bilinear algebra of the Lie algebra
${\rm Lie}(\Fs^{2\fl}\RL ( \TAN))$ is the graded bialgebra
\cite{G-H-V}, \cite{Hum},
\[
W ( {\rm Lie}(\Fs^{2\fl}\RL ( \TAN ) ) )
=A ( {\rm Lie}(\Fs^{2\fl}\RL ( \TAN ) ))\times
S ( {\rm Lie}(\Fs^{2\fl}\RL ( \TAN ) ))\;.\]

\subsection{Definition: Connection on the tangent bisemispace}

Let $\Lambda (SP^{f_\ell }(X^{\rm sv}_{R})\times
SP^{f_\ell }(X^{\rm sv}_{L}))$ denote the graded differential algebra of differential forms of
$SP^{f_\ell }(X^{\rm sv}_{R})\times
SP^{f_\ell }(X^{\rm sv}_{L})$ and let
$\Lambda (\Delta ^{2\fl}_{R}\times
\Delta ^{2\fl}_{L})$ denote the graded differential algebra of differential forms of $(\Delta ^{2\fl}_{R}\times
\Delta ^{2\fl}_{L})$.

A connection on the fibered tangent bisemispace
$(\Delta ^{2\fl}_{R}\times
\Delta ^{2\fl}_{L})$ consists in a bilinear mapping
$f^{\TAN}\RL$ of 
$A^1({\rm Lie} ( \Fs^{2\fl}\RL ( \TAN ) ))$ in the subspace of bielements of degree one of the bialgebra
$\Lambda (\Delta ^{2\fl}_{R}\times
\Delta ^{2\fl}_{L})$.
\vskip 11pt

\subsection{Proposition}

{\em
Let $I_s({\rm Lie} ( \Fs^{2\fl}\RL ( \TAN ) )) $ denote the subalgebra of invariant elements of\linebreak 
$S({\rm Lie} ( \Fs^{2\fl}\RL ( \TAN ) )) $} \cite{Car} {\em which is the algebra of symmetric multilinear forms\linebreak
$V({\rm Lie} ( \Fs^{2\fl}\RL ( \TAN ) )) $ on
${\rm Lie} ( \Fs^{2\fl}\RL ( \TAN ) ) $.

Then, there is a homomorphism:
\[ h^{\TAN}\RL: \quad
V({\rm Lie} ( \Fs^{2\fl}\RL ( \TAN ) ))\To
H_{2\fl} (
\Lambda (SP^{f_\ell }(X^{\rm sv}_{R})
\times SP^{f_\ell }(X^{\rm sv}_{L})))\;, \]
corresponding to the Chern-Weil homomorphism, such that a connection, associated to the homomorphism
$I_s({\rm Lie} ( \Fs^{2\fl}\RL ( \TAN ) )) 
\To \Lambda (SP^{f_\ell }(X^{\rm sv}_{T_R})
\times SP^{f_\ell }(X^{\rm sv}_{T_L}))$, is equivalent to the existence of a bilinear ``{\bf contracting}'' fibre
$\Fs^{2\fl}\RL ( \TAN )$ in the tangent bibundle
$\TAN[SP^{f_\ell }(X^{\rm sv}_{R})
\times SP^{f_\ell }(X^{\rm sv}_{L})]$ which implies that:
\[ H_{2\fl,}(\Delta ^{2\fl}\RL,\Fs\RL^{2\fl}(\TAN ))
\simeq \AdFRepsp (\GL_{\fl}(\cit\times\cit ))\]
and thus that:
\[ H^{2\fl} [ X^{2n_\ell }_{R}[2\fl ]\times X^{2n_\ell }_{L}[2\fl ],
\Delta ^{2\fl}_{R}\times\Delta ^{2\fl}_{L}]
\simeq \FRepsp (\GL_{\fl}(F_{\o\omega }\times\cit)
\otimes (F_{\omega }\times\cit )))\]
where $X^{2n_\ell }_{L}[2\fl]$ denotes a (semi)scheme of dimension $2n_\ell $ shifted in $2\fl$ dimensions according to
$X^{2n_\ell }_{L}[2\fl]=
\Delta ^{2\fl}_{L}\times
X ^{2n_\ell -2\fl}_{L}$.}
\vskip 11pt

\begin{proof}
\Bena
\item The connection $f^{\TAN}\RL$ on the fibered tangent bisemispace $(\Delta ^{2\fl}_{R}\times\Delta ^{2\fl}_{L})$ can be extended to a homomorphism:
\[ f^{\TAN '}\RL: \quad
A({\rm Lie} ( \Fs^{2\fl}\RL ( \TAN ) ))
\To
\Lambda (\Delta ^{2\fl}_{R}\times\Delta ^{2\fl}_{L})\;.\]

\item According to H. Cartan \cite{Car}, the knowledge of
$(SP^{f_\ell }(X^{\rm sv}_{R})
\times SP^{f_\ell }(X^{\rm sv}_{L}))$ together with the connection $f^{\TAN}\RL$ is sufficient to know
\[ H^{2\fl}(X^{2n_\ell }_{R}[2\fl]\times X^{2n_\ell }_{L}[2\fl],
\Delta ^{2\fl}_{R}\times\Delta ^{2\fl}_{L}) \;.\]

\item Thus, the existence of a connection $f^{\TAN}\RL$, associated to the knowledge of 
$\Lambda (\Delta ^{2\fl}_{R}\times\Delta ^{2\fl}_{L})$ via the homomorphism $f^{\TAN'}\RL$, is equivalent to the existence of a bilinear fibre $\Fs^{2\fl}\RL ( \TAN )$ on
$(SP^{f_\ell }(X^{\rm sv}_{R})
\times SP^{f_\ell }(X^{\rm sv}_{L}))$.

\item If this bilinear fibre is contracting, we have that the homology of this bilinear fibre is given by:
\[ H_{2\fl}(\Delta ^{2\fl}_{R}\times\Delta ^{2\fl}_{L},\Fs\RL^{2\fl}(\TAN ))
\simeq \AdFRepsp (\GL_{\fl}(\cit\times\cit ))\;.\]

\item And thus, the bilinear cohomology with coefficients in
$(\Delta ^{2\fl}_{R}\times\Delta ^{2\fl}_{L})$ must be developed according to:
\begin{multline*}
H^{2\fl}(X^{2n_\ell }_{R}[2\fl]\times X^{2n_\ell }_{L}[2\fl],
\Delta ^{2\fl}_{R}\times\Delta ^{2\fl}_{L})\\
= H_{2\fl}(\Delta ^{2\fl}_{R}\times\Delta ^{2\fl}_{L},\Fs\RL^{2\fl}(\TAN ))\hspace*{5cm}\\
\times
H^{2\fl}[X^{2n_\ell-2\fl }_{R}\times X^{2n_\ell-2\fl }_{L},SP^{f_\ell }(X^{\rm sv}_{T_R})\times SP(X^{\rm sv}_{T_L}))\\
\simeq \AdFRepsp (\GL_{\fl}(\cit\times\cit)
\times \FRepsp (\GL_{\fl}(F_{\o\omega }\times F_{\omega }))\end{multline*}
taking into account that \cite{Pie3}:
\begin{multline*}
H^{2\fl}(X^{2n_\ell-2\fl }_{R}\times X^{2n_\ell-2\fl }_{L},SP^{f_\ell }(X^{\rm sv}_{R})\times SP(X^{\rm sv}_{T_L})\\
\simeq
\FRepsp (\GL_{\fl}(F_{\o\omega }\times F_{\omega }))\;.\qedhere\end{multline*}
\Ee
\end{proof}

\subsection{Definition}
{\bbf The bilinear Lie algebra
$\Lie(G_\Ls^{(n)}(F_{\o\omega }\times F_{\omega }))$} of the Lie bilinear semigroup
$G_\Ls^{(n)}(F_{\o\omega }\times F_{\omega })$ associated with the bilinear semigroup
$G^{(n)}(F_{\o\omega }\times F_{\omega })$ can be introduced by noting that
$\Lie(G_\Ls^{(n)}(F_{\o\omega }\times F_{\omega }))$ naturally decomposes according to:
\[ \Lie(G_\Ls^{(n)}(F_{\o\omega }\times F_{\omega }))
= \Lie(T_\Ls^{(n)}(F_{\o\omega }))
\otimes \Lie(T_\Ls^{(n)}(F_{\omega })\]
where $\Lie(T_\Ls^{(n)}(F_{\omega }))$ is the linear Lie algebra of the Lie semigroup
$T_\Ls^{(n)}(F_{\omega })$ associated with the linear semigroup $T^{(n)}(F_{\omega })$ being the representation semispace of the group
$T_n(F_{\omega })$ of upper triangular matrices (see section 1.6).

The Lie algebra
$\Lie(G_\Ls^{(n)}(F_{\o\omega }\times F_{\omega }))$
corresponds to the bilinear tensor product of a vector semispace
$\Lie(T_\Ls^{(n)}(F_{\omega }))$ by is dual 
$\Lie(T_\Ls^{(n)}(F_{\o\omega }))$, also called a bilinear vector semispace \cite{Pie5}.

Each element of
$\Lie(T_\Ls^{(n)}(F_{\omega }))$
\resp{$\Lie(T_\Ls^{(n)}(F_{\o\omega }))$}
defines a one-parameter semigroup of automorphisms of
$T_\Ls^{(n)}(F_{\omega })$
\resp{$T_\Ls^{(n)}(F_{\o\omega })$}, which are the right translations by a one-parameter subgroup of
$T_\Ls^{(n)}(F_{\omega })$
\resp{$T_\Ls^{(n)}(F_{\o\omega })$}.

More exactly, the bilinear Lie algebra
$\Lie(G_\Ls^{(n)}(F_{\o\omega }\times F_{\omega }))$ is defined by the two conditions:
\Bean
\item $\Lie(G_\Ls^{(n)}(F_{\o\omega }\times F_{\omega }))$ is a bilinear vector semispace over the product
$F_{\o\omega }\times F_{\omega }$ of sets of completions;

\item to each pair
$(\tau ^t_{F_{\o\omega }},\tau _{F_{\omega }})$, with
$\tau ^t_{F_{\o\omega }}\in\Lie(T_\Ls^{(n)}(F_{\o\omega }))$ and
$\tau _{F_{\omega }}\in\Lie(T_\Ls^{(n)}(F_{\omega }))$,
corresponds an element of 
$\Lie(G_\Ls^{(n)}(F_{\o\omega }\times F_{\omega }))$,
noted $[ \tau ^t_{F_{\o\omega }},\tau _{F_{\omega }}]$
\Bi
\item which is linear with respect to
$\tau ^t_{F_{\o\omega }}$ and to $\tau _{F_{\omega }}$;
\item whose value is given by
$[ \tau ^t_{F_{\o\omega }},\tau _{F_{\omega }}]
= \tau ^t_{F_{\o\omega }}\centerdot\tau _{F_{\omega }}
- \tau _{F_{\omega }}\centerdot\tau^t_{F_{\o\omega }}$;
\item which verifies the Jacobi identity.
\Ei
\Ee

\subsection{Proposition}

{\em
The bilinear cohomology with coefficients in the tangent bisemispace
$(\Delta ^{2\fl}_{R}
\times\Delta ^{2\fl}_{L})$, noted
$H^{2\fl}(X^{2n_\ell }_{R}[2\fl]
\times X^{2n_\ell }_{L}[2\fl],
\Delta ^{2\fl}_{R}
\times\Delta ^{2\fl}_{L})$, is in one-to-one correspondence with the Lie algebra of the general bilinear semigroup
$\GL_{\fl}(F_{\o\omega }
\times F_{\omega })$:
\[
H^{2\fl}(X^{2n_\ell }_{R}[2\fl]
\times X^{2n_\ell }_{L}[2\fl],
\Delta ^{2\fl}_{R}
\times\Delta ^{2\fl}_{L})
\simeq \Lie (\GL_{\fl}(F_{\o\omega }
\times F_{\omega }))\;.\]
}
\vskip 11pt

\begin{proof}
Indeed, according to proposition 2.8, we have that
\begin{multline*}
H^{2\fl}(X^{2n_\ell }_{R}[2\fl]
\times X^{2n_\ell }_{L}[2\fl],
\Delta ^{2\fl}_{R}
\times\Delta ^{2\fl}_{L})\\
\begin{aligned}
&\simeq \FRepsp ( \GL_{\fl}(F_{\o\omega }\otimes \cit)\times (F_{\omega }\otimes \cit )))\\
&= \AdFRepsp (\GL_{\fl}(\cit\times \cit )) 
\times \FRepsp ( \GL_{\fl}(F_{\o\omega }
\times F_{\omega }))
\end{aligned}\end{multline*}
from which it clearly appears that:
\[
\Lie(\GL_{\fl} ( F_{\o\omega }\times F_{\omega })
=\FRepsp(\GL_{\fl} ( F_{\o\omega }\otimes \cit )
\times ( F_{\omega }\otimes \cit ))
\]
since the bilinear fibre
$\Fs^{2\fl}\RL(\TAN)$ of the tangent bibundle
$\TAN[SP^{f_\ell }(X^{\rm sv}_{R})\times
SP^{f_\ell }(X^{\rm sv}_{L})]$ is precisely the adjoint functional representation space
$\AdFRepsp (\GL_{\fl}(\cit\times \cit))$ of\linebreak
$\GL_{\fl}(\cit\times \cit)$.\end{proof}
\vskip 11pt

\subsection{Proposition}

{\em
Let 
$H^{2n_\ell }({\u c}_*(X^{\rm sv}\RL),Z\RL(2n_\ell ))$
be the bilinear cohomology of the  \SV\ pure bisemimotives
${\u c}_*(X^{\rm sv}_{T\RL})$ with coefficients in the product, right by left, of \SV\  sub(bisemi)motives of complex codimension $n_\ell $.

Then, the cohomology of the corresponding \SV\ mixed bisemimotive\linebreak
$M_{DM_R}(X^{\rm sv}_{R})\times
M_{DM_L}(X^{\rm sv}_{L})$, noted
 $M_{DM\RL}(X^{\rm sv}\RL)$, can be reached throughout the following endomorphism:
 \begin{multline*}
 \HH D_{2\fl}: \quad
 H^{2n_\ell }({\u c}_*(X^{\rm sv}\RL),Z\RL(2 n_\ell ))\\
 \To
 H^{2n_\ell-2\fl}(M_{DM\RL}(X^{\rm sv}\RL
),Z\RL(2n_\ell [2\fl] )\end{multline*}
where $Z\RL(2n_\ell [2\fl] )$ is the product, right by left, of \SV\  mixed subbisemimotives of complex codimension $n_\ell $ shifted in $\fl$ complex dimensions and written
$(X^{2n_\ell }_{R} [2\fl]\times X^{2n_\ell }_{L} [2\fl])$ in section 2.4, such that the cohomology of the \SV\ mixed bisemimotive decomposes according to:
\begin{multline*}
 H^{2n_\ell-2\fl }(M_{DM\RL} ( X^{\rm sv}\RL
),Z\RL ( 2n_\ell [2\fl])\\
= H_{2\fl}(\Delta ^{2\fl}_{R}\times\Delta ^{2\fl}_{L},
\Fs^{2\fl}\RL ( \TAN ))\times
H^{2n_\ell }({\u c}_*(X^{\rm sv}\RL),Z\RL(2 n_\ell ))\;.\end{multline*}
}
\vskip 11pt

\begin{proof}
Taking into account that
$X^{2n_\ell }_{L}[2\fl] = \Delta ^{2\fl}_{L}\times X^{2n_\ell -2\fl}_{L}$
\resp{$X^{2n_\ell }_{R}[2\fl] = \Delta ^{2\fl}_{R}\times X^{2n_\ell -2\fl}_{R}$} according to section 2.4, we have that the cohomology of
${\u c}_*(X^{\rm sv}\RL)$, submitted to a translation functor acted by the tangent bibundle 
$\TAN[SP^{f_\ell }(X^{\rm sv}_{R})
\times SP^{f_\ell }(X^{\rm sv}_{L})]$, is transformed according to:
\begin{multline*}
 \HH D_{2\fl}: \quad
 H^{2n_\ell }({\u c}_*(X^{\rm sv}\RL),Z\RL(2n_\ell ))\\
 \To
 \begin{aligned}[t]
&  H^{2\fl}(M_{DM\RL}(X^{\rm sv}\RL
),\Delta ^{2\fl}_{R}
\times\Delta ^{2\fl}_{L}
)\\
& \qquad 
\oplus H^{2n_\ell -2\fl}({\u c}_*(X^{\rm sv}\RL),
X^{2n_\ell -2\fl}_{R}
\times X^{2n_\ell -2\fl}_{L})\\
&=( H_{2\fl}(\Delta ^{2\fl}_{R}
\times \Delta ^{2\fl}_{L},\Fs^{2\fl}\RL(\TAN) ))\\
& \qquad \times
\bigl[ H^{2\fl}({\u c}_*(X^{\rm sv}\RL),
SP^{f_\ell }( X^{\rm sv}_{R})\times SP^{f_\ell }(X^{\rm sv}_{L} )) \bigr.\\
& \qquad \bigl. \oplus 
H^{2n_\ell -2\fl}({\u c}_*(X^{\rm sv}\RL),
X^{2n_\ell -2\fl}_{R}\times X^{2n_\ell -2\fl}_{L} )\bigr] \\
& = H_{2\fl}(\Delta ^{2\fl}_{R}
\times\Delta ^{2\fl}_{L}, \Fs^{2\fl}\RL ( \TAN ))\\
& \qquad \times H^{2n_\ell }({\u c}_*(X^{\rm sv}\RL),Z\RL (2 n_\ell ))
\end{aligned}\end{multline*}
such that $H_{2\fl}(\Delta ^{2\fl}_{R}
\times\Delta ^{2\fl}_{R}, \Fs^{2\fl}\RL ( \TAN ))
\times H^{2n_\ell }({\u c}_*(X^{\rm sv}\RL),Z\RL (2 n_\ell ))$ be noted\linebreak
$H^{2n_\ell-2\fl}(M_{DM\RL}(X^{\rm sv}\RL
),Z\RL(2 n_\ell[2\fl] )$.

These equalities essentially result from proposition 2.8.

It then results that the cohomology
$H^{2n_\ell-2\fl}(M_{DM\RL}(X^{\rm sv}\RL
),Z\RL(2n_\ell [2\fl] )$
of the \SV\ mixed bisemimotive decomposes into the
$(2\fl)$-homology with coefficients in the bilinear fibre
$\Fs^{\fl}\RL(\TAN)$ acting on the $(2n_\ell )$-coho\-mology of the \SV\ pure bisemimotive ${\u c}_*(X^{\rm sv}\RL)$.\end{proof}

\subsection{Proposition}

{\em
The bilinear cohomology of the \SV\ mixed bisemimotive $M_{DM\RL}(X^{\rm sv}\RL)$ is in bijection with the functional representation space of the bilinear general semigroup
$\GL_{n_\ell }(F_{\o\omega }\times F_\omega )$ shifted in
$\fl$-complex dimensions:
\begin{multline*}
H^{2n_\ell-2\fl}(M_{DM\RL}(X^{\rm sv}\RL
),Z\RL(2n_\ell [2\fl]] )\\
\simeq \FRepsp(\GL_{n_\ell [\fl]}((F_{\o\omega }\otimes\cit)\times
((F_{\omega }\otimes\cit))\end{multline*}
where $\FRepsp(\GL_{n_\ell [\fl]}((F_{\o\omega }\otimes\cit)\times
((F_{\omega }\otimes\cit)$ is a condensed notation for\linebreak
$\AdFRepsp(\GL_{\fl}(\cit\times\cit)\times
\FRepsp(\GL_{n_\ell }(F_{\o\omega }\times F_\omega ))$.
}
\vskip 11pt

\begin{proof}
According to proposition 2.11, we have that:
\begin{multline*}
H^{2n_\ell-2\fl}(M_{DM\RL}(X^{\rm sv}\RL
),Z\RL(2 n_\ell[2\fl] )\\
 \begin{aligned}[t]
&= (H_{2\fl}(\Delta ^{2\fl}_R
\times \Delta ^{2\fl}_L,\Fs^{2\fl}\RL(\TAN) )\\
& \qquad \times
\bigl[ H^{2\fl}({\u c}_*(X^{\rm sv}\RL),
SP^{f_\ell }(X^{\rm sv}_{R})\times SP^{f_\ell }(X^{\rm sv}_{L} )) \bigr.\\
& \qquad \bigl. \oplus
H^{2n_\ell -2\fl}({\u c}_*(X^{\rm sv}\RL),
X^{2n_\ell -2\fl}_{R}\times X^{2n_\ell -2\fl}_{L} ))\bigr]\;.
\end{aligned}\end{multline*}
And, propositions 2.8 and 2.10 give the following isomorphisms:
\Bi
\item $H_{2\fl}(\Delta ^{2\fl}_{R}
\times \Delta ^{2\fl}_{L},\Fs^{2\fl}\RL(\TAN) )\simeq
\AdFRepsp(\GL_{\fl}(\cit\times\cit))$;

\item $H^{2\fl}({\u c}_*(X^{\rm sv}\RL),
SP^{f_\ell }(X^{\rm sv}_{R})\times SP^{f_\ell }(X^{\rm sv}_{L} ))\simeq \FRepsp (\GL_{\fl}(F_{\o\omega }\times F_\omega ))$;

\item $H^{2n_\ell -2\fl}({\u c}_*(X^{\rm sv}\RL),
X^{2n_\ell -2\fl}_{R}\times X^{2n_\ell -2\fl}_{L} )
 \simeq \FRepsp (\GL_{n_\ell -\fl}(F_{\o\omega }\times F_\omega ))$;
\Ei
leading to:
\begin{multline*}
H^{2n_\ell-2\fl}(M_{DM\RL}(X^{\rm sv}\RL
),Z\RL(2n_\ell [2\fl] )\\
 \begin{aligned}[t]
 &\simeq \AdFRepsp(\GL_{\fl}(\cit\times \cit ))
\\ & \qquad
\times
\bigl[ \FRepsp(\GL_{\fl}(F_{\o\omega }\times F_\omega ))
 \oplus \FRepsp(\GL_{n_\ell -\fl}(F_{\o\omega }\times F_\omega ))\bigr]\\
&= \AdFRepsp(\GL_{\fl}(\cit\times \cit ))\times
 \FRepsp(\GL_{n_\ell }(F_{\o\omega }\times F_\omega ))\hspace*{2cm}\qedhere
\end{aligned}\end{multline*}
\end{proof}

\subsection[Higher Chow semigroups]{Higher Chow semigroups \cite{Blo}, \cite{Gil1}}

\Bi
\item According to proposition 1.17, the \lr  \SV\ subsemimotive 
$Z_L(2n_\ell )$
\resp{$Z_R(2n_\ell )$} of complex dimension $n_\ell $ can be isomorphic to the \lr semicycle 
${\rm CY}^{2n_\ell }(Y_L)$
\resp{${\rm CY}^{2n_\ell }(Y_R)$} belonging to the
$2n_\ell $-th Chow semigroup 
${\rm CH}^{2n_\ell }(Y_L)$
\resp{${\rm CH}^{2n_\ell }(Y_R)$}:
\begin{align*}
Z_L(2n_\ell )&\simeq {\rm CY}^{2n_\ell }(Y_L)\in {\rm CH}^{2n_\ell }(Y_L)\\[11pt]
\rresp{Z_R(2n_\ell )&\simeq {\rm CY}^{2n_\ell }(Y_R)\in {\rm CH}^{2n_\ell }(Y_R)}\;.\end{align*}

\item Similarly, the \lr \SV\  mixed submotive 
$Z_L(2n_\ell [2\fl])$
\resp{$Z_R(2n_\ell [2\fl])$} of complex dimension $n_\ell $, shifted in $\fl$ complex dimensions, can  be isomorphic to the \lr cycle 
${\rm CY}^{2n_\ell }(Y_L,[2\fl])$
\resp{${\rm CY}^{2n_\ell }(Y_L,[2\fl])$} of complex dimension $n_\ell $, shifted in $\fl$ complex dimensions, belonging to the $2n_\ell $-th higher Chow semigroup
${\rm CH}^{2n_\ell }(Y_L,[2\fl])$
\resp{${\rm CH}^{2n_\ell }(Y_R,[2\fl])$}:
\begin{align*}
Z_L(2n_\ell [2\fl]) &\simeq {\rm CY}^{2n_\ell }(Y_L,[2\fl])\in {\rm CH}^{2n_\ell }(Y_L,[2\fl])\\[11pt]
\rresp{Z_R(2n_\ell [2\fl]) &\simeq {\rm CY}^{2n_\ell }(Y_R,[2\fl])\in {\rm CH}^{2n_\ell }(Y_R,[2\fl])}.
\end{align*}
\Ei

\subsection{Representation of the general bilinear shifted semigroup}

According to section 1.15,  the product, right by left, of  cycles
${\rm CY}^{2n_\ell }(Y_R)\times {\rm CY}^{2n_\ell }(Y_L)$ leads to:
\begin{align*}
Z_R(2n_\ell )\times Z_L(2n_\ell )
&\simeq {\rm CY}^{2n_\ell }(Y_R)\times {\rm CY}^{2n_\ell }(Y_L)
\\[11pt]
&\simeq \FRepsp (\GL_{n_\ell }(F_{\o\omega }\times F_\omega ))\;.\end{align*}

By the same way, the product, right by left, of the \SV\  mixed subsemimotives of complex dimension $n_\ell $ shifted in $\fl$ complex dimensions  gives rise to the bijections:
\begin{align*}
Z\RL(2n_\ell [2\fl])
&\simeq {\rm CY}^{2n_\ell } ( Y_R[2\fl] ) \times {\rm CY}^{2n_\ell }(Y_L[2\fl] )\\
&\simeq \FRepsp (\GL_{n_\ell [\fl]}
((F_{\o\omega }\otimes \cit)\times (F_\omega \otimes \cit))\;.\end{align*}
\vskip 11pt

\subsection{Proposition}

{\em
Taking into account the isomorphism 
$i_{M-X}$ between the bilinear cohomology of the  \SV\ pure bisemimotive ${\u c}_*(X^{\rm sv}\RL)$ and the bilinear coho\-mology of  $Y\RL$ introduced in proposition 1.17, as well as the endomorphism $\HH D_{2\fl}$ between the bilinear cohomology of
${\u c}_*(X^{\rm sv}\RL)$ and the corresponding cohomology of the \SV\ mixed bisemimotive
$M_{DM\RL}(X^{\rm sv}\RL)$, we are led to the following commutative diagram:
\[ \begin{CD}
H^{2n_\ell }({\u c}_*(X^{\rm sv}\RL), Z\RL(2 n_\ell ))
@>{i_{M-X}}>> H^{2n_\ell }(Y\RL,{\rm CY}^{2n_\ell }(Y_R)\times {\rm CY}^{2n_\ell }(Y_L))\\
@V{\HH D_{2\fl}}VV
@V{\HH D^{X-M}_{2\fl}}VV\\
\begin{aligned}[t]
&H^{2n_\ell -2\fl}(M_{DM\RL} (X^{\rm sv}\RL
),\\
& \qquad \qquad
Z\RL(2n_\ell[2\fl])\end{aligned}
@>{i^{\rm sc}_{M-X}}>>
\begin{aligned}[t]
& H^{2n_\ell -2\fl}
( Y\RL[2\fl],\\
& \quad
{\rm CY}^{2n_\ell }(Y_R[2\fl] )\times
{\rm CY}^{2n_\ell }(Y_L[2\fl] )\end{aligned}
\end{CD}\]
where $Y\RL[2\fl]$ is the  bisemigroup $Y\RL$ shifted in $\fl$ complex dimensions on its right and left parts.
}


\subsection{Bilinear mixed cohomology}

The introduction in this chapter of the bilinear cohomology of mixed bisemimotives naturally leads to precise what must be a general bilinear mixed (or shifted) cohomology referring to the introduction of a general bilinear cohomology in section 3.2 of \cite{Pie3} and taking into account the isomorphisms:
\begin{align*}
Z\RL(2n_\ell[2\fl])
&\simeq {\rm CY}^{2n_\ell }(Y\RL[2\fl] )\\
&\simeq \FRepsp(\GL_{n_\ell[\fl]}(F_{\o\omega} \otimes\cit)\times(F_{\omega }\otimes\cit))
\end{align*}
and the bilinear mixed homology
\[ H_{2\fl}(\Delta ^{2\fl}\RL,\Fs^{2\fl}\RL(\TAN))\simeq 
\AdFRepsp(\GL_{\fl}(\cit\times \cit ))\]
associated with the tangent bibundle $\TAN[\SP^{f_\ell}(X^{\rm sv}\RL)]$.

\subsection{Proposition}
{\em
{\bbf A general bilinear mixed cohomology theory is a contravariant bifunctor\/}:
\begin{multline*}
\HH^{2i-2k}: \quad
\{\text{smooth abstract shifted bisemivarieties\ } G^{(n)}((F^+_{\o v}\otimes\rit)\times (F^+_v\otimes\rit))\}\\
\To \begin{aligned}[t]
&\{\text{graded functional representation spaces of the complete shifted}\\ 
 &\text{bilinear semigroups\ }\GL_{2i[2k]}((F^+_{\o v}\otimes\rit)\times (F^+_v\otimes\rit))\}\ , 0\le k\le i\;,\end{aligned}
\end{multline*}
given by
\[ H^{2i-2k}(G^{(n)}((F^+_{\o v}\otimes\rit)\times (F^+_v\otimes\rit)),
\FRepsp(\GL_{2i[2k]}
((F^+_{\o v}\otimes\rit)\times (F^+_v\otimes\rit)))\]
where
\begin{multline*}
\FRepsp ( \GL_{2i[2k]}
 ( ( F^+_{\o v}\otimes\rit ) \times ( F^+_v\otimes\rit ) ) ) \\
=\AdFRepsp ( \GL_{2k} ( \rit\times\rit )) \times
\Repsp ( \GL_{2i} ( F^+_{\o v}\times F^+_v ) )\;.
\end{multline*}
This general bilinear mixed cohomology is characterized by:
\Be
\item {\bf isomorphic embeddings\/}
\begin{align*}
&G^{(n)}((F^+_{\o v}\otimes\rit)\times (F^+_{v}\otimes\rit))\\
&\hspace{2cm} \overset{\sim}{\scalebox{1.5}{$\hookrightarrow$}}\quad
G^{(n)}((F_{\o \omega }\otimes\cit)\times (F_{\omega }\otimes\cit))\\
&\FRepsp(\GL_{2i[2k]}
((F^+_{\o v}\otimes\rit)\times (F^+_v\otimes\rit)))\\
&\hspace{2cm} \overset{\sim}{\scalebox{1.5}{$\hookrightarrow$}}\quad
\FRepsp(\GL_{2i[2k]}((F_{\o \omega }\otimes\cit)\times (F_\omega \otimes\cit)))
\end{align*}
of ``real'' shifted bisemivarieties into their complex equivalents.

\item {\bbf mixed (or shifted) bisemicycle maps\/}:
\begin{multline*}
\gamma ^{i[k]}_{G^{(n)}_{\o v\times v}}:
\quad \bZ^{i[k]}(G^{(n)}((F^+_{\o v}\otimes \rit)\times(F^+_v\otimes\rit)))\\
\To
H^{2i-2k} ( G^{(n)} ( ( F^+_{\o v}\otimes \rit ) \times ( F^+_v\otimes\rit )),\FREPSP (
\GL_{2i[2k]} ( ( F^+_{\o v}\otimes \rit )\times ( F^+_v\otimes\rit ))))
\end{multline*}
where $\bZ^{i[k]}$ denotes the bilinear semigroup of mixed bisemicycles of codimension $i$ shifted in $k$ dimensions.

\item {\bf Hodge mixed (or shifted) bisemicycles\/}:
\[H^{2i-2k}(G^{(n)}((F_{\o \omega }\otimes \cit)\times(F_\omega \otimes\cit)),\FREPSP(
\GL_{2i[2k]}(F^+_{\o v}\otimes \rit)\times(F^+_v\otimes\rit)))\]
from the abstract ``complex'' shifted bisemivariety
$G^{(n)}((F_{\o \omega }\otimes\cit)\times (F_{\omega }\otimes\cit))$ to the functional representation space of the ``real"" shifted general bilinear semigroup
$\GL_{2i[2k]}(F^+_{\o v}\otimes \rit)\times(F^+_v\otimes\rit))$.

There is the shifted bifiltration $F^{p[r]}\RL$ given by:
\begin{multline*}
F^{p[r]}\RL: \quad H^{2i-2k}(G^{(n)}(\pt\times\pt),-)\\
= \bigoplus_{\substack{i=p+q\\k=r-s}} H^{(2p-2r)+(2q-2s)}(
G^{(n)}((F_{\o \omega }\otimes \cit)\times(F_\omega  \otimes\cit)),\\
\FREPSP(
\GL_{2p[2r]+2q[2s]}(F^+_{\o v}\otimes \rit)\times(F^+_v\otimes\rit)))\;.
\end{multline*}

\item {\bf a Künneth standard conjecture\/}:

implying that the projectors on $H^{2i-2k}(G^{(n)}(\pt\times\pt),-)$ (see {\em\cite{Pie3}\/}) are induced by mixed (or shifted) compactified bisemicycles
$\CY^{i[k]}(G^{(n)}(F^+_{\o v}\otimes\rit)\times(F^+_v\otimes\rit))\subset
\bZ^{i[k]}(G^{(n)}(\pt\times\pt))$ decomposing into rational mixed (or shifted) subbisemicycles according to the conjucacy class representatives of
$\GL_{2i[2k]}((F^+_{\o v}\otimes\rit)\times(F^+_v\otimes\rit))$.

\item {\bf a Künneth biisomorphism\/}:
\begin{multline*}
H^{2i-2k} ( G^{(n)} ( F^+_{\o v}\otimes\rit ) ,\FREPSP ( \GL_{2i[2k]} ( F^+_{\o v}
\otimes\rit ) ) ) \\
\otimes_{F^+_{\o v}\times F^+_v}
H^{2i-2k} ( G^{(n)} ( F^+_{v}\otimes\rit ) ,
\FREPSP ( \GL_{2i[2k]} ( F^+_{v}\otimes\rit ) ) )\\
\To 
\begin{aligned}[t]
& H^{2i-2k} ( G^{(n)} ( F^+_{\o v}\otimes\rit ) \times  ( F^+_{v}\otimes \rit ) ,\\
& \qquad \FREPSP (\GL_{2i[2k]} ( F^+_{\o v}\otimes\rit) ) ) \times ( F^+_{v}\otimes \rit ) ) )
\end{aligned}\end{multline*}
in such a way that
\begin{multline*}
H^{2p-2r}(G^{(n)}(F^+_{\o v}\otimes\rit)\times (F^+_{v}\otimes \rit),
\FREPSP(\GL_{2p[2r]} ( F^+_{\o v}\otimes\rit)
\times ({F^+_{v}\otimes \rit} ) ) ) \\
\otimes
H^{(2i-2k)-(2p-2r)}(G^{(n)}(F^+_{\o v}\otimes\rit)\times (F^+_{v}\otimes \rit),\qquad \qquad\\
\FREPSP(\GL_{(2i[2k]-(2p[2r])}(F^+_{\o v}\otimes\rit)))\times ({F^+_{v}\otimes \rit}))\\
\To \begin{aligned}[t]
& H^{0}(G^{(n)}(F^+_{\o v}\otimes\rit)\times (F^+_{v}\otimes \rit),\\
& \qquad \FREPSP(\GL_{1}(F^+_{\o v}\otimes\rit)))\times ({F^+_{v}\times \rit}))
\end{aligned}
\end{multline*}
is the bilinear version of the mixed intersection cohomology according to section 3.2 of {\em \cite{Pie3}\/}.
\Ee
}

\bpr The introduction of the general bilinear mixed cohomology follows from the introduction of general bilinear cohomology in section 3.2 of \cite{Pie3} to which we refer.\epr

%% file: LanglandsIII_3.tex
\section[Bilinear $K$-homology associated with an elliptic bioperator]{\bbf Bilinear $K$-homology associated with an elliptic bioperator}

\subsection{Modular conjugacy classes of \SV\ pure subsemimotives}

Let $Z_L(2n_\ell )\equiv X^{2n_\ell }_{L}$
\resp{$Z_R(2n_\ell )\equiv X^{2n_\ell }_{R}$} be a \SV\ \lr pure subsemimotive of complex dimension $n_\ell $, i.e. a \lr subpresheaf with transfers (or correspondences).

Referring to sections 2.1 and 2.14, we have that:
\begin{align*}
X^{2n_\ell }_{R}\times X^{2n_\ell }_{L}
&\simeq \FRepsp (\GL_{n_\ell }(F_{\o\omega }\times F_\omega ))\\
&\simeq {\rm CY}^{2n_\ell }(Y_R)\times {\rm CY}^{2n_\ell }(Y_L)\\
&= \{ {\rm CY}^{2n_\ell }(Y_R(j,m_j))\times {\rm CY}^{2n_\ell }(Y_L(j,m_j))\}_{j,m_j}
\end{align*}
in such a way that the product, right by left, of $2n_\ell $-dimensional semicycles decomposes according to the set of conjugacy class representatives of
$\GL_{n_\ell }(F_{\o\omega }\times F_\omega )$.

It follows that the product, right by left, $X^{2n_\ell }_{R}\times X^{2n_\ell }_{L}$ of \SV\ subsemimotives of complex dimension $n_\ell $ also decomposes according to the conjugacy class representatives of 
$\GL_{n_\ell }(F_{\o\omega }\times F_\omega )$:
\[ 
X^{2n_\ell }_{R}\times X^{2n_\ell }_{L}
=  \{ X^{2n_\ell }_{R}(j,m_j)\times X^{2n_\ell }_{L}(j,m_j)\}_{j,m_j}\;.\]

\subsection{Modular conjugacy classes of \SV\ mixed subsemimotives}
Let
\[ Z_L(2n_\ell [2\fl]) \equiv X^{2n_\ell }_{L}[2\fl]\;, 
\rresp{Z_R(2n_\ell [2\fl]) \equiv X^{2n_\ell }_{R}[2\fl]}
\]
be the \lr \SV\ mixed subsemimotive of complex dimension $n_\ell $ shifted in $2\fl$ dimensions.

Then, as in section 3.1, we have that:
\begin{align*}
X^{2n_\ell }_{R}[2\fl]\times X^{2n_\ell }_{L}[2\fl]
&\simeq \FRepsp (\GL_{n_\ell[\fl] }
(F_{\o\omega }\otimes \cit )\times
(F_{\omega }\otimes \cit ))\\
&\simeq {\rm CY}^{2n_\ell }(Y_R[2\fl])\times {\rm CY}^{2n_\ell }(Y_L[2\fl])\\
&= \{ {\rm CY}^{2n_\ell }(Y_R,[2\fl],(j,m_j))\times {\rm CY}^{2n_\ell }(Y_L,[2\fl],(j,m_j))\}_{j,m_j}
\end{align*}
where ${\rm CY}^{2n_\ell }(Y_R,[2\fl],(j,m_j))$ is the $m_j$-th representative of the $j$-th conjugacy class of the $2n_\ell $-th (semi)cycle
${\rm CY}^{2n_\ell }$ shifted in $2\fl$ dimensions of the semigroup $Y_R$.
\vskip 11pt

\subsection[Definition: Differential bioperator $D^{2\fl}_R\otimes D^{2\fl}_L$]{\bbf Definition: Differential bioperator $D^{2\fl}_R\otimes D^{2\fl}_L$}

Let $D^{2\fl}_R\otimes D^{2\fl}_L$ be the product  of a right linear differential (elliptic) operator
$D^{2\fl}_R$ acting on $2\fl$ variables by its left equivalent.  This bioperator is defined by its biaction
\[ D^{2\fl}_R\otimes D^{2\fl}_L: \qquad
X^{2n_\ell }_{R}\times X^{2n_\ell }_{L}\To
X^{2n_\ell }_{R}[2\fl]\times X^{2n_\ell }_{L}[2\fl]
\]
from the \SV\ pure subbisemimotive
$ X^{2n_\ell }_{R}\times X^{2n_\ell }_{L}$ to the corresponding mixed subbisemimotive
$X^{2n_\ell }_{R}[2\fl]\times X^{2n_\ell }_{L}[2\fl]$ shifted in $(\fl)$ complex dimensions.

In fact, $(D^{2\fl}_R\otimes D^{2\fl}_L)$ acts on the set of smooth bisections
$\{X^{2n_\ell }_{R}(j,m_j)\times X^{2n_\ell }_{L}(j,m_j)\}_{j,m_j}$ of $X^{2n_\ell }_{R}\times X^{2n_\ell }_{L}$.

$D^{2\fl}_L$ has the form 
$D^{2\fl}_L(X^{n_\ell }_{L}(j,m_j))
= \sum\limits_a \dots  \sum\limits_t\ 
\partial^a_1\ \dots\ \partial^t_{2\fl}
(X^{2n_\ell }_{L}(j,m_j))$, where\linebreak $\partial_{2\fl}=i\ \dfrac d{dx_{2\fl}}$ is the differential operator with respect to the $2\fl$-th variable $x_{2\fl}$.
\vskip 11pt

\subsection[Definition: Symbol of the bioperator $D^{2\fl}_R\otimes D^{2\fl}_L$]{\bbf Definition: Symbol of the bioperator $D^{2\fl}_R\otimes D^{2\fl}_L$}

Referring to section 2.4 and lemma 2.5, 
$X^{2n_\ell }_{R}[2\fl]\times X^{2n_\ell }_{L}[2\fl]$ develops according to:
\[X^{2n_\ell }_{R}[2\fl]\times X^{2n_\ell }_{L}[2\fl]=
(\Delta ^{2\fl}_{R}\times \Delta ^{2\fl}_{L})\times
(X ^{2n_\ell -2\fl}_{R}\times X ^{2n_\ell -2\fl}_{L})\]
where 
$(\Delta ^{2\fl}_{R}\times \Delta ^{2\fl}_{L})$ is the total space of the tangent bibundle
$\TAN[SP^{f_\ell }(X^{\rm sv}_{R}]\times\linebreak \TAN[SP^{f_\ell }(X^{\rm sv}_{L}]$ and develops according to:
\begin{align*}
\Delta ^{2\fl}_{R}\times \Delta ^{2\fl}_{L}
&\simeq \AdFRepsp(\GL_{\fl}(\cit\times\cit ))\times \FRepsp ( \GL_{\fl} ( F_{\o\omega }\times F_\omega ))\\
&=\FRepsp ( \GL_{\fl} (( F_{\o\omega }\times\cit )\times ( F_{\o\omega }\times \cit )))\;.\end{align*}
Then, referring to the classical definition \cite{A-S} of the symbol $\sigma (D)$ of a differential linear operator $D$, we can admit that the symbol
$\sigma (D^{2\fl}_R\otimes D^{2\fl}_L)$ of the bioperator
$(D^{2\fl}_R\otimes D^{2\fl}_L)$ can be introduced by \cite{Ma1}, \cite{Ma2}, \cite{L-T}:
\[
\sigma (D^{2\fl}_R\otimes D^{2\fl}_L)
=\FRepsp (P_{\fl}((F_{\o\omega }\otimes \cit)\times
(F_{\omega }\otimes \cit)))\;, \]
i.e. by the unitary functional representation space of
$\GL_{\fl}((F_{\o\omega }\otimes \cit)\times
(F_{\omega }\otimes \cit))$ given by the functional representation space of the shifted bilinear parabolic semigroup\linebreak
$P_{\fl}((F_{\o\omega }\otimes \cit)\times
(F_{\omega }\otimes \cit))$ \cite{Pie3}.
\vskip 11pt

\subsection{Definition}

The differential bioperator $(D^{2\fl}_R\otimes D^{2\fl}_L)$ is elliptic if its symbol
$\sigma (D^{2\fl}_R\otimes D^{2\fl}_L)$ is invertible.
\vskip 11pt

In connection with the work of G. Kasparov \cite{Kas} who constructed a general $K_*\ K^*$ functor on the categories of compact operators and Hilbert modules, {\bbf we shall introduce a bilinear $K_*\ K^*$ functor on the categories of elliptic bioperators and 
products, right by left, of \SV\  pure semimotives allowing to set up a bilinear version of the index theorem\/}
\cite{B-F-M}, \cite{A-H}, \cite{Mil}, \cite{Jan1}.
\vskip 11pt

\subsection[Chern character of the pure bimotive ${\u c}_*(X^{\rm sv}\RL)$]{\bbf Chern character of the pure bimotive ${\u c}_*(X^{\rm sv}\RL)$}

Let $H^{2n_\ell }({\u c}_*(X^{\rm sv}\RL),
X ^{2n_\ell }_{R}\times X ^{2n_\ell }_{L})$ be the bilinear cohomology of the  \SV\ pure bisemimotive ${\u c}_*(X^{\rm sv}\RL)$ and let
\[
H^*({\u c}_*(X^{\rm sv}\RL))
=\bigoplus_{n_\ell }H^{2n_\ell }
({\u c}_*(X^{\rm sv}\RL),X ^{2n_\ell }_{R}\times X ^{2n_\ell}_{L})\]
denote the total bilinear cohomology of
${\u c}_*(X^{\rm sv}\RL)$.

Taking into account the definition of a pure bisemimotive
${\u c}_*(X^{\rm sv}\RL)$ as being a functor from
$X^{\rm sv}\RL$ of complex dimension $\ell $ to the chain bicomplex associated to the product, right by left, of abelian semigroups
\[ \mathop{\sqcup}\limits_{i_\ell }\Hom_{Sm_L(k)\times
Sm_R(k)}
(\dot\Sigma  _{R}\times
\dot\Sigma _{L},
SP^{i_\ell }(X^{\rm sv}_{R})\times
SP^{i_\ell }(X^{\rm sv}_{L}))\;, \]
we can introduce the product, right by left, of abelian semigroups generated by the complex vector bundles \cite{Laf} over $X^{\rm sv}_{R}\times
X^{\rm sv}_{L}$ and noted 
$K^*(X^{\rm sv}\RL)$

$K^*(X^{\rm sv}\RL)$ is then the $K$-cohomology associated to the pure bisemimotive
${\u c}_*(X^{\rm sv}\RL)$.
\vskip 11pt

{\bbf The total Chern character \cite{Gil2} in the bilinear $K$-cohomology of the pure bimotive}
${\u c}_*(X^{\rm sv}\RL)$ is then given by the homomorphism:
\[ ch^*({\u c}_*(X^{\rm sv}\RL)): \qquad
K^*(X^{\rm sv}\RL))\To
H^*({\u c}_*(X^{\rm sv}\RL))\]
and defined, as classically, according to:
\[ch^*({\u c}_*(X^{\rm sv}\RL))
=\sum^{n_\ell /\ell }_{i_\ell =1}
e^{\gamma _{i_\ell }}\centerdot
e^{\gamma _{i_\ell }}\;, \qquad i_\ell \centerdot \ell =n_\ell \le n\]
where the $\gamma _{i_\ell }$ result from the factorization \cite{Hir}:
\[ (1+c_1\ x+\dots +c_{i_\ell }\ x^{i_\ell }+\dots+
c_{n_\ell /\ell }\ x^{n_\ell /\ell })
=\prod_{i_\ell }(1+\gamma _{i_\ell }\ x)\]
with the $c_{i_\ell }$ being the Chern classes.
\vskip 11pt

\subsection[Bilinear $K$-homology]{Bilinear $K$-homology \cite{B-D-F}, \cite{Blo}, \cite{Gil1}}

\Bi
\item Let $H_{2\fl}
(\Delta ^{2\fl}_{R}\times\Delta ^{2\fl}_{L}),\Fs^{*}\RL(\TAN))
\simeq \AdFRepsp(\GL_{\fl}(\cit\times\cit))$ be the
  homology with coefficients in the bilinear fiber $\Fs^{2\fl}\RL(\TAN)$ of the tangent bibundle
$\TAN[SP^{f_\ell }(X^{\rm sv}_{R})\times 
SP^{f_\ell }(X^{\rm sv}_{L})]$ and let
\begin{align*}
H_*(\Delta ^*_{R}\times\Delta ^*_{L},\Fs^{*}\RL(\TAN)
&= \bigoplus_{\fl}
H_{2\fl}(\Delta ^{2\fl}_{R}\times\Delta ^{2\fl}_{L}),\Fs^{2\fl}\RL(\TAN))\\
&\simeq \bigoplus_{\fl}
\AdFRepsp(\GL_{\fl}(\cit\times\cit))
\end{align*}
be the total bilinear homology with coefficients in the set of bilinear fibres\linebreak
$\Fs^{2\fl}\RL(\TAN)$ such that, for $i_\ell \times\ell =n_\ell =(\fl)+(n_\ell -\fl)$, $f_\ell \le i_\ell $ and $\fl\le n_\ell \le n$.

\item Then, a {\bbf bilinear $K$-homology}, noted
$K_*(SP^{FL}(X^{\rm sv}\RL))$, can be introduced as being the product, right by left, of abelian semigroups generated by the set of tangent bibundles
$\TAN[SP^{f_\ell }(X^{\rm sv}_{R})\times 
SP^{f_\ell }(X^{\rm sv}_{L})]$, for all $\fl\le n_\ell \le n$, on the product
$SP^{f_\ell }(X^{\rm sv}_{R})\times 
SP^{f_\ell }(X^{\rm sv}_{L})$ of smooth presheaves.

$K_*(SP^{FL}(X^{\rm sv}\RL))$ is the $K$-homology associated to the pure bisemimotive 
${\u c}_*(X^{\rm sv}\RL)$.

\item {\bbf The Chern character in this bilinear $K$-homology} is thus given by the homomorphism:
\[
ch_*({\u c}_*(X^{\rm sv}\RL))
: \quad 
K_*(SP^{FL}(X^{\rm sv}\RL))\To
H_*(\Delta ^*_{R}\times\Delta ^*_{L},\Fs^{*}\RL(\TAN))\;.\]
Taking into account that
$H_*(\Delta ^*_{R}\times\Delta ^*_{L},\Fs^{*}\RL(\TAN))$ is the homology of\linebreak
$(\Delta ^*_{R}\times\Delta ^*_{L})$ with coefficients in the set of bilinear fibres
$\Fs^{2\fl}\RL(\TAN)$ which are ``contracting'', the total Chern character in this bilinear $K$-homology will be defined by
\[ch_*({\u c}_*(X^{\rm sv}\RL))
= \sum_{f_\ell}e^{-\gamma _{f_\ell}}\centerdot
e^{-\gamma _{f_\ell}}\]
such that the $\gamma _{f_\ell}$ are obtained from a formal factorisation $\sum\limits_{f_\ell} c_{-f_\ell}\ x^{f_\ell}
=\linebreak \prod\limits_{f_\ell}(1-\gamma _{f_\ell}\ x)$ where the
$c_{-f_\ell}\in H^{-2\fl}(\Fs^{\fl}_{\left\{ \substack{R \\ L}\right.} (\TAN),\zit)$ are Chern classes associated with the homology.
\Ei

\subsection{Proposition}

{\em The total Chern character
$ch^*(M_{DM\RL}(X^{\rm sv}\RL))$ of the \SV\ mixed bisemimotive
$M_{DM\RL}(X^{\rm sv}\RL)$ in the mixed bilinear $K$-homology-$K$-cohomology is given by {\bbf the homomorphism:}
\begin{multline*}
ch^*(M_{DM\RL}(X^{\rm sv}\RL)): \quad
K_*(SP^{FL}(X^{\rm sv}\RL)) \times
K^*(X^{\rm sv}\RL)\\
\To
H_*(\Delta ^*_{R}\times\Delta ^*_{L},
\Fs^*\RL(\TAN))\times H^*({\u c}_*(X^{\rm sv}\RL))
\end{multline*}
such that:
\[
ch^*(M_{DM\RL}(X^{\rm sv}\RL))=
ch_*({\u c}_*(X^{\rm sv}\RL))
\times ch^*({\u c}_*(X^{\rm sv}\RL))\]
{\bbf corresponds to a bilinear version of the index theorem.}
}
\vskip 11pt

\begin{proof}
Taking into account that:
\[ 
H^*({\u c}_*(X^{\rm sv}\RL))
= \bigoplus_{n_\ell } H^{2n_\ell }
({\u c}_*(X^{\rm sv}\RL),X^{2n_\ell }_{R}\times X^{2n_\ell }_{L})\]
and that:
\[ H_*(\Delta ^*_{R}\times\Delta ^*_{L},
\Fs^*\RL(\TAN))
=\bigoplus_{\fl} H_{2\fl}
(\Delta ^{2\fl}_{R}\times\Delta ^{2\fl}_{L},
\Fs^{2\fl}\RL(\TAN))\]
according to sections 3.6 and 3.7, as well as the decomposition of the cohomology of the \SV\ mixed bimotive into:
\begin{multline*}
H^{2n_\ell -2\fl}(M_{DM\RL}(X^{\rm sv}\RL),
X^{2n_\ell }_{R}[2\fl]\times
X^{2n_\ell }_{L}[2\fl])\\
\To
H_{2\fl}
(\Delta ^{2\fl}_{R}\times\Delta ^{2\fl}_{L},
\Fs^{2\fl}\RL(\TAN))
\times H^{2n_\ell}(
{\u c}_*(X^{\rm sv}\RL),X^{2n_\ell }_{R}\times X^{2n_\ell }_{L})\;, \end{multline*}
we have that:
\begin{multline*}
\qquad H_*
(\Delta ^*_{R}\times\Delta ^*_{L})
\times H^*(
{\u c}_*(X^{\rm sv}\RL))\\
= \bigoplus_{n_\ell } \bigoplus_{\fl }
H^{2n_\ell -2\fl}(M_{DM\RL}(X^{\rm sv}\RL),
X^{2n_\ell }_{R}[2\fl]\times
X^{2n_\ell }_{L}[2\fl])
\end{multline*}
is the total bilinear cohomology of
$M_{DM\RL}(X^{\rm sv}\RL)$, noted
$H^*(M_{DM\RL}(X^{\rm sv}\RL)$.

Similar arguments can be used to prove that
\[ K^*(M_{DM\RL}(X^{\rm sv}\RL))=
 K_*(SP^{FL}(X^{\rm sv}\RL))\times
 K^*(X^{\rm sv}\RL)\]
 is the mixed bilinear $K$-homology-$K$-cohomology associated with the \SV\ mixed bisemimotive
$M_{DM\RL}(X^{\rm sv}\RL)$.

And, thus, it follows that:
\begin{multline*}
\qquad \qquad ch^*(M_{DM\RL}(X^{\rm sv}\RL)): \qquad 
K^*(M_{DM\RL}(X^{\rm sv}\RL))\\
\To
H^*(M_{DM\RL}(X^{\rm sv}\RL))\qquad \qquad
\end{multline*}
is the total Chern character of the \SV\ mixed bisemimotive.
\end{proof}
\vskip 11pt

\subsection{Corollary}

{\em
Let $(D^{2\fl}_R\otimes D^{2\fl}_L)$ be a differential bioperator defined by its biaction:
\[
D^{2\fl}_R\otimes D^{2\fl}_L: \qquad
X^{2n_\ell }_{R}\times X^{2n_\ell }_{L}\To
X^{2n_\ell }_{R}[2\fl]\times X^{2n_\ell }_{L}[2\fl]
\]
from the \SV\ pure subbisemimotive
$X^{2n_\ell }_{R}\times X^{2n_\ell }_{L}$ to the corresponding mixed subbimotive
$X^{2n_\ell }_{R}[2\fl]\times X^{2n_\ell }_{L}[2\fl]$.

Let 
\[
ch_*(D^{2\fl}_R\otimes D^{2\fl}_L): \quad
K_*(SP^{f_\ell}(X^{\rm sv}\RL))
\To
H_{2\fl}
(\Delta ^{2\fl}_{R}\times\Delta ^{2\fl}_{L},
\Fs^{2\fl}\RL(\TAN))
\]
be an element of the Chern character
$ch_*({\u c}_*(X^{\rm sv}\RL))$ associated with the biaction of\linebreak 
$(D^{2\fl}_R\otimes D^{2\fl}_L)$
on $(X^{2n_\ell }_{R}\times X^{2n_\ell }_{L})$.

Let 
\[
ch^*(X^{2n_\ell }_{R}\times X^{2n_\ell }_{L}): \quad
K^*(SP^{f_\ell}(X^{\rm sv}\RL))
\To
H^{2n_\ell }(
{\u c}_*(X^{\rm sv}\RL),X^{2n_\ell }_{R}\times X^{2n_\ell }_{L})
\]
denote an element of
$ch^*(
{\u c}_*(X^{\rm sv}\RL))$.

Then, $ch_*(D^{2\fl}_R\otimes D^{2\fl}_L)\times
ch^*(X^{2n_\ell }_{R}\times X^{2n_\ell }_{L})$ will allow to define an index
$\Ind(D^{2\fl}_R\otimes D^{2\fl}_L)$ of the elliptic bioperator which is
 different from the classical Atiyah-Singer index
$\gamma (D^{2\fl}_R\otimes D^{2\fl}_L)$ except if $i_\ell =f_\ell$.
}
\vskip 11pt

\begin{proof}
Referring to sections 3.6 and 3.7 where
\begin{align*}
ch^*({\u c}_*(X^{\rm sv}\RL))
&=\sum^{n_\ell /\ell }_{i_\ell =1}
e^{\gamma _{i_\ell }}\centerdot
e^{\gamma _{i_\ell }}\;, && i_\ell \centerdot \ell =n_\ell \le n\;, \\
\text{and} \qquad
ch_*({\u c}_*(X^{\rm sv}\RL))
&=\sum_{f_\ell}
e^{-\gamma _{f_\ell}}\centerdot
e^{-\gamma _{f_\ell}}\;,
\end{align*}
are introduced, we define
$\Ind(D^{2\fl}_R\otimes D^{2\fl}_L)$ by:
\begin{align*}
\Ind(D^{2\fl}_R\otimes D^{2\fl}_L)
&=
ch_*(D^{2\fl}_R\otimes D^{2\fl}_L)
\times ch^*(X ^{2n_\ell }_{R}\times X ^{2n_\ell }_{L})\\
&= e^{2\gamma _{i_\ell }}\centerdot e^{-2\gamma _{f_\ell }}-\delta _{i_\ell ,f_\ell}\;, \end{align*}
where $\delta _{i_\ell ,f_\ell} \begin{aligned}[t]
=0 \quad \text{if} \quad i_\ell \neq f_\ell\;,\\
=1 \quad \text{if} \quad i_\ell = f_\ell\;.
\end{aligned}$

Thus, $\Ind(D^{2\fl}_R\otimes D^{2\fl}_L)=0$ if and only if
$i_\ell = f_\ell$.  In that case, if 
$(D^{2\fl}_R\otimes D^{2\fl}_L)$ is of finite rank, 
$\Ind(D^{2\fl}_R\otimes D^{2\fl}_L)=0$ and could correspond to the classical Atiyah-Singer index, defined by
\[ \gamma (D^{2\fl}_R\otimes D^{2\fl}_L)
=\dim {\rm Ker} (D^{2\fl}_R\otimes D^{2\fl}_L)
-\dim {\rm coKer} (D^{2\fl}_R\otimes D^{2\fl}_L)\;.\qedhere
\]
\end{proof}

\subsection{Corollary}

{\em The equivalent of the classical index theorem}
\cite{Gil1} {\em for a \SV\ bisemimotive asserts that, if
\[ \{D^{2\fl}_R\otimes D^{2\fl}_L\}_{\fl}: \qquad
{\u c}_*(X^{\rm sv}\RL)
\To M_{DM\RL}(X^{\rm sv}\RL)\]
is a proper morphism from the \SV\ pure bisemimotive
${\u c}_*(X^{\rm sv}\RL)$ to the \SV\ mixed bisemimotive
$M_{DM\RL}(X^{\rm sv}\RL)$ under the action of the set of bioperators
$\{D^{2\fl}_R\otimes D^{2\fl}_L\}_{\fl}$, we have that:
\begin{align*}
ch_*({\u c}_*(X^{\rm sv}\RL))\times
ch^*({\u c}_*(X^{\rm sv}\RL))
&= ch^*(M_{DM\RL}(X^{\rm sv}\RL))\\
&= \Ind\{(D^{2\fl}_R\otimes D^{2\fl}_L)\}_{\fl}\;.\end{align*}
}
\vskip 11pt

\begin{proof}
More specifically, the index theorem would assert that:
\[
{\rm Im}\{D^{2\fl}_R\otimes D^{2\fl}_L\}_{\fl}
ch^*({\u c}_*(X^{\rm sv}\RL))
= ch^*(M_{DM\RL}(X^{\rm sv}\RL))\]
where ${\rm Im}\{D^{2\fl}_R\otimes D^{2\fl}_L\}_{\fl}$ is the image of the morphism generated by\linebreak
$\{D^{2\fl}_R\otimes D^{2\fl}_L\}_{\fl}$.

Now, if we take into account the considered notations, it appears that:
\[
\qquad \qquad{\rm Im}\{D^{2\fl}_R\otimes D^{2\fl}_L\}_{\fl}
(ch^*({\u c}_*(X^{\rm sv}\RL)))
= ch_*({\u c}_*(X^{\rm sv}\RL))\times
ch^*({\u c}_*(X^{\rm sv}\RL))
\]
since ${\rm Im}\{D^{2\fl}_R\otimes D^{2\fl}_L\}_{\fl}$ is generated by the set of bioperators
$\{D^{2\fl}_R\otimes D^{2\fl}_L\}_{\fl}$ to which the Chern character $ch_*({\u c}_*(X^{\rm sv}\RL))$ in this bilinear $K$-homology corresponds.

Furthermore, we have that:
\begin{align*}
\Ind\{(D^{2\fl}_R\otimes D^{2\fl}_L)\}_{\fl}
&= ch_*({\u c}_*(X^{\rm sv}\RL))
\times ch^*({\u c}_*(X^{\rm sv}\RL))\\
&= ch_*(M_{DM\RL}(X^{\rm sv}\RL))\;.\qedhere
\end{align*}
\end{proof}

%% file: LanglandsIII_4.tex
\section{The toroidal spectral representation of an elliptic bioperator}

Chapters 2 and 3 were essentially devoted to pure and mixed bimotives of \SV, while chapter 4 and 5 will  more particularly concern the functional representation spaces of bilinear algebraic semigroups.
\vskip 11pt

\subsection[The shifted compactified bisemispace
$\o S^{P_{\nfl}}_{\GL_{\nfl}}$]{\bbf The shifted compactified bisemispace
$\o S^{P_{\nfl}}_{\GL_{\nfl}}$}

\Bi
\item Let
\[\o S^{P_{n_\ell }}_{\GL_{n_\ell }}
=P_{n_\ell }(F^T_{\o\omega ^1_1}\times F^T_{\omega ^1_1})\setminus 
\GL_{n_\ell }(F^T_R\times F^T_L)\big/
\GL_{n_\ell }((\zit/N\ \zit)^2)\]
be the toroidal compactified bisemispace representing the double coset decomposition of the algebraic bilinear semigroup
$\GL_{n_\ell }(F^T_R\times F^T_L)$ as introduced in section 1.12.
\vskip 11pt

\item The corresponding {\bbf double coset decomposition of the bilinear general semigroup shifted in $(\fl\times \fl)$ complex dimensions}
$\GL_{\nfl }(F^T_R\otimes \cit)
\times (F^T_L\otimes \cit))$, as developed in proposition 2.12, is given by:
\begin{multline*}
\qquad  \o S^{P_{\nfl}}_{\GL_{\nfl}}
= P_{\nfl}((F^T_{\o\omega ^1}\otimes \cit)\times
(F^T_{\omega ^1}\otimes \cit))\\
\setminus \GL_{\nfl}((F^T_R\otimes\cit)\times
(F_L^T\otimes\cit)) \big/
\GL_{\nfl }((\zit/N\ \zit)^2\otimes\cit^2)\end{multline*}
in such a way that:
\Bena
\item $\o S^{P_{\nfl}}_{\GL_{\nfl}}
=\Repsp(\GL_{\nfl}((F^T_{\o\omega }\otimes\cit)\times
(F^T_{\omega }\otimes\cit)))$ implies that
\begin{multline*}
 \FRepsp(\GL_{\nfl}((F^T_{\o\omega }\otimes\cit)\times(F^T_{\omega }\otimes\cit)))\\
= \AdFRepsp(\GL_{\fl }((\cit\otimes\cit)\times
 \FRepsp(\GL_{n_\ell}(F^T_{\o\omega }\otimes F^T_{\omega }))\end{multline*}
according to proposition 2.12.

\item The {\bbf shifted complex bilinear parabolic semigroup}
$P_{\nfl}((F^T_{\o\omega ^1}\otimes \cit)\times
(F^T_{\omega ^1}\otimes \cit)$ is generated from its unshifted equivalent
$P_{n_\ell }((F^T_{\o\omega ^1}\times F^T_{\omega ^1})$ by the shift homomorphism:
\[
SH_{P_{n_\ell }}: \qquad
P_{n_\ell }(F^T_{\o\omega ^1} \times F^T_{\omega ^1})
\To P_{\nfl}((F^T_{\o\omega ^1}\otimes \cit)\times
(F^T_{\omega ^1}\otimes \cit))\;.\]
A similar shift homomorphism can be introduced in order to generate\linebreak
$\GL_{\nfl}((F^T_{\o\omega }\otimes \cit)\times
(F^T_{\omega }\otimes \cit))$.

\item The bilinear arithmetic subgroup
$\GL_{n_\ell }((\zit/N\ \zit)^2)$, generating a
$(\zit/N\ \zit)^2$-bilattice in
$\o S^{P_{n_\ell }}_{\GL_{n_\ell }}$, is transformed by the shift homomorphism:
\[
SH_{\GL_{n_\ell }}: \qquad
\GL_{n_\ell }((\zit/N\ \zit)^2)
\To \GL_{\nfl}((\zit/N\ \zit)^2)\otimes \cit^2)\]
into a {\bbf shifted bilinear arithmetic subgroup}
$\GL_{\nfl}((\zit/N\ \zit)^2)\otimes \cit^2)$ in such a way that:
\Bi
\item the (functional) representation space of
$\GL_{\nfl}((\zit/N\ \zit)^2)\otimes \cit^2)$ corresponds to the Lie algebra of
$\GL_{n_\ell }((\zit/N\ \zit)^2))$:
\[ \Lie (
\GL_{n_\ell}(\zit/N\ \zit)^2)
\approx
\FRepsp(\GL_{\nfl}( ( \zit/N \ \zit )^2)\otimes \cit^2)\]
by considerations similar as given in proposition 2.10.

\item a shifted pseudoramified Hecke bialgebra
$\Hs\RL(\nfl)$, generated by the shifted pseudoramified Hecke bioperators
$T_R(\nfl;t)\otimes T_L(\nfl;t)$, has a representation in
$\GL_{\nfl}((\zit/N\ \zit)^2)\otimes \cit^2)$ as developed in the next proposition.
\Ei
\Ee
\Ei

\subsection{Proposition}

{\em Let $\GL_{\nfl}((\zit/N\ \zit)^2\otimes \cit^2)$ be the shifted bilinear arithmetic subgroup generated from
$\GL_{n_\ell }((\zit/N\ \zit)^2)$.

Then, the pseudoramified Hecke bialgebra, generated by all the shifted pseudoramified Hecke bioperators $T_R(\2nfl;t)\otimes T_L(\2nfl;t)$, is a shifted pseudoramified bialgebra of Hecke noted $\Hs\RL(\2nfl)$.
}
\vskip 11pt

\begin{proof}
\Bena
\item Referring to sections 1.4 and 1.8, a shifted maximal order of 
$F^T_\omega $
\resp{$F^T_{\o\omega }$} will be given by
$(\Os  _{F^T_{\omega }}\otimes\cit)$
\resp{$(\Os  _{F^T_{\o\omega }}\otimes\cit)$}.

Then, a lattice of dimension $2n_\ell $, noted 
$\Lambda ^{2n_\ell }_{\omega }$
\resp{$\Lambda ^{2n_\ell }_{\o\omega }$}
shifted in $2\fl$ dimensions will be introduced by:
\[\Lambda ^{\2nfl}_{\omega }=\Lambda ^{2n_\ell }_{\omega } \otimes_{[\fl]}\cit\;,
\rresp{\Lambda ^{\2nfl}_{\o\omega }=\Lambda ^{2n_\ell }_{\o\omega } \otimes_{[\fl]}\cit},\]
where the tensor product $\otimes_{[\fl]}$ bears on the $\fl$ shifted complex dimensions, and will be defined by the isomorphism:
\begin{align*}
\Lambda ^{2n_\ell }_{\omega } \otimes_{[\fl]}\cit
&\simeq T_{\nfl}((\zit/N\ \zit)\otimes\cit)\\
&=T_{\nfl}(\Os  _{F^T_{\omega }}\otimes\cit)\\
\rresp{\Lambda ^{2n_\ell }_{\o\omega } \otimes_{[\fl]}\cit
&\simeq T^t_{\nfl}((\zit/N\ \zit)\otimes\cit)\\
&=T^t_{\nfl}(\Os  _{F^T_{\o\omega }}\otimes\cit)}
\end{align*}
leading to:
\[ \Lambda ^{\2nfl}_{\o\omega }\otimes \Lambda ^{\2nfl}_{\omega }
\simeq \GL_{\nfl}((\zit/N \ \zit)^2)\otimes\cit^2)\;.\]

\item According to proposition 1.10, the $(j,m_j)$-th coset representative $U_{j,m_{j_R}}(\2nfl)\times U_{j,m_{j_L}}(\2nfl)$ of the shifted pseudoramified Hecke bioperator $T_R(\2nfl;t)\otimes T_L(\2nfl;t)$ will be given by:
\begin{multline*}
U_{j,m_{j_R}}(\2nfl)\times U_{j,m_{j_L}}(\2nfl)\\
= [d_{\nfl}((\zit/N\ \zit)\otimes\cit)\times
d_{\nfl}((\zit/N\ \zit)\otimes\cit)]\\
\times [u^t_{\nfl}((\zit/N\ \zit)\otimes\cit)\times
u_{\nfl}((\zit/N\ \zit)\otimes\cit)]
\end{multline*}
taking into account the Gauss bilinear decomposition.
\qedhere
\Ee
\end{proof}
\vskip 11pt

\subsection{Proposition}

{\em The differential bioperator 
$(D^{2\fl}_R\otimes
D^{2\fl}_L)\in \Ds_R\otimes\Ds_L$ maps the bisemisheaf 
$(\widehat M^{2n_\ell }_{T_R}\otimes
\widehat M^{2n_\ell }_{T_L})$ 
 into the corresponding perverse bisemisheaf
$(\widehat M^{2n_\ell }_{T_R}[2\fl]\otimes \widehat M^{2n_\ell }_{T_L}[2\fl])$ according to:
\[ D^{2\fl}_R\otimes D^{2\fl}_L: \qquad
(\widehat M^{2n_\ell }_{T_R}\otimes \widehat M^{2n_\ell }_{T_L})\To
(\widehat M^{2n_\ell }_{T_R}[2\fl]\otimes \widehat M^{2n_\ell }_{T_L}[2\fl])\]
such that
$(\widehat M^{2n_\ell }_{T_R}[2\fl]\otimes \widehat M^{2n_\ell }_{T_L}[2\fl])$ is the tensor product of perverse (semi)sheaves which are $\Ds_R\otimes \Ds_L$-bisemimodules.
}
\vskip 11pt

\begin{proof} \Bena
\item The action of the differential bioperator
$(D^{2\fl}_R\otimes D^{2\fl}_L)$ on
$\widehat M^{2n_\ell }_{T_R}\otimes \widehat M^{2n_\ell }_{T_L}$ corresponds to the shift homomorphism:
\[
SH_{G_{n_\ell }}: \qquad \GL_{n_\ell }(
F^T_{\o\omega }\times F^T_{\omega })
\To \GL_{\nfl}((F^T_{\o\omega }\otimes\cit)\times (F^T_\omega \otimes\cit))\;, \]
as introduced in section 4.1, since 
$M^{2n_\ell }_{T_R}\otimes M^{2n_\ell }_{T_L}$ is the representation space of\linebreak $\GL_{n_\ell }(F^T_{\o\omega }\times F^T_\omega )$ (see section 1.14 and \cite{Pie3}).

\item $\widehat M^{2n_\ell }_{T_R}[2\fl]\otimes \widehat M^{2n_\ell }_{T_L}[2\fl]$ is the tensor product of perverse sheaves because it is an object of the derived category of
$(\o S^{P_{n_\ell }}_{\GL_{n_\ell }})$ (see sections 2.4, 1.14 and \cite{B-B-D}).\qedhere
\Ee
\end{proof}
\vskip 11pt



\subsection{Proposition}

{\em The bilinear cohomology of the shifted compactified bisemispace
$\o S^{P_{\Nfl}}_{\GL{\Nfl}}$ is\linebreak
$H^{2n_\ell -2\fl}(\o S^{P_{\Nfl}}_{\GL{\Nfl}},
\widehat M^{2n_\ell }_{T_R}[2\fl]\otimes
\widehat M^{2n_\ell }_{T_L}[2\fl])$ and is isomorphic to the
bilinear cohomology of the \SV\ mixed bimotive
$M_{DM\RL}(
X^{\rm sv}\RL)$, noted
\[ H^{2n_\ell -2\fl}(
M_{DM\RL} (
X^{\rm sv}\RL ),
X^{2n_\ell }_{R}[2\fl]\times
X^{2n_\ell }_{L}[2\fl] )\;.\]
}
\vskip 11pt

\begin{proof}
The isomorphism:
\begin{multline*}
H^{2n_\ell -2\fl}(
M_{DM\RL} 
(X^{\rm sv}\RL),
X^{2n_\ell }_{R}[2\fl]\times
X^{2n_\ell }_{L}[2\fl] )\\
\simeq
H^{2n_\ell -2\fl}(
\o S^{P_{\Nfl}}_{\GL_{\Nfl}}
,
\widehat M^{2n_\ell }_{R}[2\fl]\otimes
\widehat M^{2n_\ell }_{L}[2\fl] )
\end{multline*}
results from the isomorphisms
\begin{multline*}
\qquad \qquad H^{2n_\ell -2\fl}(
M_{DM\RL} 
(X^{\rm sv}\RL),
X^{2n_\ell }_{R}[2\fl]\times
X^{2n_\ell }_{L}[2\fl] )\\
\simeq
\FRepsp ( \GL_{\nfl} (( F^T_{\o\omega }\otimes\cit )\times
 ( F^T_\omega \otimes\cit ))) \qquad 
\end{multline*}
and
\begin{multline*}
\qquad \qquad H^{2n_\ell -2\fl}(
\o S^{P_{\Nfl}}_{\GL{\Nfl}},
\widehat M^{2n_\ell }_{R}[2\fl]\otimes
\widehat M^{2n_\ell }_{L}[2\fl] )\\
\simeq
\FRepsp ( \GL_{\nfl} (( F^T_{\o\omega }\otimes\cit )\times
( F^T_\omega \otimes\cit ))) \;.\qquad \qedhere
\end{multline*}
\end{proof}
\vskip 11pt

\subsection{Proposition}

{\em
The bilinear cohomology
\begin{multline*}
\qquad \qquad H^{2n_\ell -2\fl}(
\o S^{P_{\Nfl}}_{\GL_{\Nfl}}
,
\widehat M^{2n_\ell }_{T_R}[2\fl]\otimes
\widehat M^{2n_\ell }_{T_L}[2\fl] )\\
\simeq
{\rm CY}^{2n_\ell }(Y_R,[2\fl])\times {\rm CY}^{2n_\ell }(Y_L,[2\fl])\\
= \{ ( 
{\rm CY}^{2n_\ell }(Y_R,[2\fl],(j,m_j))\times 
({\rm CY}^{2n_\ell }(Y_L,[2\fl],(j,m_j))\}_{j,m_j}
\end{multline*}
is in bijection with the decomposition in equivalence classes ``$j$, having multiplicities $m^{(j)}$, of the products, right by left, of the right  cycles
${\rm CY}^{2n_\ell }(Y_R,[2\fl],(j,m_j))$
shifted in $2\fl$ dimensions by their left equivalents
${\rm CY}^{2n_\ell }(Y_L,[2\fl],(j,m_j))$ such that:
\begin{align*}
&{\rm CY}^{2n_\ell }(Y_R,[2\fl],(j,m_j))\in {\rm CH}^{2n_\ell }(Y_R,2\fl)\\
\rresp{&{\rm CY}^{2n_\ell }(Y_L,[2\fl],(j,m_j))\in {\rm CH}^{2n_\ell }(Y_L,2\fl)}
\end{align*}
where ${\rm CH}^{2n_\ell }(Y_R,2\fl)$
\resp{${\rm CH}^{2n_\ell }(Y_L,2\fl)$}
is the $2n_\ell $-th higher Chow semigroup (see section 2.13).
}
\vskip 11pt

\begin{proof} This results from section 2.13, the isomorphisms of proposition 4.4 and proposition 4.6.
\end{proof}

\subsection{Proposition}

{\em
The decomposition of the product, right by left,
${\rm CY}^{2n_\ell }(Y_R,[2\fl])\times
{\rm CY}^{2n_\ell }(Y_L,[2\fl])$ of  cycles of codimension $2n_\ell $ shifted in $2\fl$  dimensions into equivalence class representatives corresponds to the decomposition of
$\GL_{\nfl}((F^T_{\o\omega }\otimes\cit)\times
(F^T_\omega \otimes\cit))$, into the set of products, right by left, of conjugacy class representatives
$g^{2n_\ell }_{T\RL}([2\fl],(j,m_j))$ shifted in $2\fl$ dimensions:
\begin{multline*}
\qquad \qquad {\rm CY}^{2n_\ell }(Y_R,[2\fl])\times
{\rm CY}^{2n_\ell }(Y_L,[2\fl])\\
=\{ (
{\rm CY}^{2n_\ell } (Y_R,[2\fl],(j,m_j) )\times
{\rm CY}^{2n_\ell } (Y_L,[2\fl],(j,m_j) )
\}_{j,m_j}\\
\simeq \{ \phi (g^{2n_\ell }_{T_R}([2\fl],(j,m_j) ) )\times
 \phi (g^{2n_\ell }_{T_L}([2\fl],(j,m_j) ) )
 \}_{j,m_j}
 \end{multline*}
 in such a way that each cofunction
 $\phi (g^{2n_\ell }_{T_R}([2\fl],(j,m_j) ) )$
 \resp{function $\phi (g^{2n_\ell }_{T_L}([2\fl],(j,m_j) ) )$} on
$g^{2n_\ell }_{T_R}([2\fl],(j,m_j) ) $
\resp{$g^{2n_\ell }_{T_L}([2\fl],(j,m_j) ) $}
be a $n_\ell $-dimensional complex semitorus
$T^{2n_\ell }_R([2\fl],(j,m_j))$
\resp{$T^{2n_\ell }_L([2\fl],(j,m_j))$}
shifted in $\fl$ complex dimensions and localized in the lower \resp{upper} half space (toroidal case only).
}
\vskip 11pt

\begin{proof}
We have that:
\begin{multline*}
\qquad \qquad \Repsp(\GL_{\nfl} ( ( F^T_{\o\omega }\otimes\cit )\times
( F^T_{\omega }\otimes\cit ) ) \\
\begin{aligned}
&= \{ 
g^{2n_\ell }_{T_R}([2\fl],(j,m_j) ) )\times
g^{2n_\ell }_{T_L}([2\fl],(j,m_j) )\}_{j,m_j}\\
&\equiv \{g^{2n_\ell }_{T\RL}([2\fl],(j,m_j) )\}_{j,m_j}
\;.\end{aligned}
\end{multline*}
Indeed, $\GL_{n_\ell }(F^T_{\o\omega }\times F^T_\omega )$ decomposes into conjugacy class representatives
$g^{(2n_\ell )}_{T\RL}(j,m_j)\equiv
g^{(2n_\ell )}_{T_R}(j,m_j)\times
g^{(2n_\ell )}_{T_L}(j,m_j)$ consisting in products, right by left, of $n_\ell $-dimensional complex semitori \cite{Pie3}.

Then,  the bilinear complete semigroup
$\GL_{\nfl} ( ( F^T_{\o\omega }\otimes\cit )\times
( F^T_{\omega }\otimes\cit ) $, shifted in $\fl$ complex dimensions from
$\GL_{n_\ell } (  F^T_{\o\omega }\times
 F^T_{\omega })$, has for conjugacy class representatives $g^{(2n_\ell )}_{T\RL}([2\fl],(j,m_j))$ which are the conjugacy class representatives
$ g^{(2n_\ell )}_{T\RL}(j,m_j)$ of
$\GL_{n_\ell }(F^T_{\o\omega }\times F^T_\omega )$ shifted in $2\fl$ dimensions.  And, thus,
$\phi (g^{2n_\ell }_{T_R}([2\fl],(j,m_j) ) )\times
 \phi (g^{2n_\ell }_{T_L}([2\fl],(j,m_j) ) )$
 consists of the product, right by left, of the analytical representatives of $n_\ell $-dimensional complex semitori shifted in $\fl$ dimensions.
\end{proof}
\vskip 11pt

\subsection{Proposition}

{\em Let $z_{n_\ell }=\sum\limits^{2{n_\ell }}_{\alpha =1} z_\alpha \ \vec e_\alpha $,
 $z_{\fl}=\sum\limits^{2{\fl}}_{\beta =1} z_\beta  \ \vec e_\beta  $ and
 $(z_{n_\ell -\fl)}=\sum\limits^{2({n_\ell -\fl})}_{\gamma  =1} z_\gamma \ \vec e_\gamma $
 be respectively a vector of 
$\cit^{n_\ell }$,
$\cit^{\fl }$ and
$\cit^{n_\ell -\fl }$.

Then, every \lr $n_\ell $-dimensional complex semitorus
$T^{2n_\ell }_L([2\fl],(j,m_j))$
\resp{$T^{2n_\ell }_R([2\fl],(j,m_j))$} shifted in $\fl$ complex dimensions and localized in the upper \resp{lower} half space has the following analytic development:
\begin{align*}
T^{2n_\ell }_L([2\fl],(j,m_j))
&\simeq E_{2\fl}(2n_\ell ,j,m_j)\ \lambda ^{\half}(2n_\ell ,j,m_j)\ e^{2\pi ijz_{n_\ell }}\\
&\equiv \prod^{2(n_\ell -\fl)}_{c=1}\ \lambda_c ^{\half}(2n_\ell ,j,m_j)\ e^{2\pi ijz_{n_\ell -\fl }}\\
&\qquad \prod^{2\fl}_{d=1}\ \lambda_d ^{\half}(2n_\ell ,j,m_j) \ E_d(2n_\ell ,j,m_j)\ e^{2\pi ijz_{\fl }}\\
\rresp{%
T^{2n_\ell }_R([2\fl],(j,m_j))
&\simeq E_{2\fl}(2n_\ell ,j,m_j)\ \lambda ^{\half}(2n_\ell ,j,m_j)\ e^{-2\pi ijz_{n_\ell }}}
\end{align*}
where:
\Bi
\item $\lambda ^{\half}(2n_\ell ,j,m_j)= \prod^{2(n_\ell -\fl)}_{c=1}\ \lambda_c ^{\half}(2n_\ell ,j,m_j)\ \prod^{2\fl}_{d=1}
\ \lambda_d ^{\half}(2n_\ell ,j,m_j) \simeq (j\centerdot N)^{2n_\ell }$;

\item
$E_{2\fl}(2n_\ell ,j,m_j )= \prod^{2\fl}_{d=1}\  E_d(2n_\ell ,j,m_j)$ is the shift of the Hecke character\linebreak
$\lambda ^{\half}(2n_\ell ,j,m_j)$ in such a way that
$ E_d(2n_\ell ,j,m_j)$ be a generator of the Lie algebra component $d_{\fl}(\Os  _{F_{\omega _{j,m_j}}}) \in D_{\fl}(\Os  _{F_\omega } )$ (see proposition 1.10).
\Ei
}
\vskip 11pt

\begin{proof}
\Bena
\item According to propositions 4.3, 4.4 and 4.5, the cohomology\linebreak
$H^{2n_\ell -2\fl}(
\o S^{P_{\Nfl}}_{\GL_{\Nfl}}
,
\widehat M^{2n_\ell }_{T_R}[2\fl]\otimes
\widehat M^{2n_\ell }_{T_L}[2\fl] )$, associated with an endomorphism of
$\o S^{P_{\nfl}}_{\GL_{\nfl}}$ into itself, decomposes into conjugacy class functional representatives
$\phi (g^{2n_\ell }_{T\RL}([2\fl],(j,m_j) ) )$ which correspond to the cosets of 
$\GL_{\nfl} ( ( F^T_{\o\omega }\otimes\cit )\times
( F^T_{\omega }\otimes\cit ))\big/
\GL_{\nfl} ( ( \zit/N\ \zit)^2\otimes\cit^2 ) $.  So, the scalar
$(E^2_{2\fl}(2n_\ell ,j,m_j)\centerdot\lambda (2n_\ell ,j,m_j)$ will correspond to the eigenvalues of the
$(j,m_j)$-th coset representative of the Hecke shifted bioperator $T_R(\2nfl;t)\otimes T_L(\2nfl;t)$, since it has a representation into the Lie algebra of
$\GL_{n_\ell } ( ( \zit/N\ \zit)^2)$ according to section 4.1 and proposition 4.2, while the scalar $\lambda (2n_\ell ,j,m_j)$ will correspond to the eigenvalues of the $(j,m_j)$-th coset representative of the Hecke bioperator
$(T_R(2n_\ell ;t)\otimes T_L(2n_\ell ;t))$ by means of the equality:
\[
\lambda (2n_\ell ,j,m_j)
= \prod^{2(n_\ell -\fl)}_{c=1}\ \lambda_c (2n_\ell ,j,m_j)
\
\prod^{2\fl}_{d=1}\  \lambda_d (2n_\ell ,j,m_j)\;.\]
Remark that the eigenvalues of the $(j,m_j)$-coset representative of
$(T_R(\2nfl ;t)\otimes T_L(\2nfl ;t)$ are partitioned into unshifted eigenvalues
$\lambda_c (2n_\ell ,j,m_j)$ and into shifted eigenvalues
$(\lambda_d (2n_\ell ,j,m_j) \centerdot E^2_d(2n_\ell ,j,m_j))$
such that:
\begin{multline*}
E_{2\fl}^2(2n_\ell ,j,m_j)\centerdot
\lambda (2n_\ell ,j,m_j) \\
=\prod^{2(n_\ell -\fl)}_{c=1}\ \lambda_c (2n_\ell ,j,m_j)\ \prod^{2\fl}_{d=1}\ \lambda_d (2n_\ell ,j,m_j) \centerdot E^2_d(2n_\ell ,j,m_j)\;.\end{multline*}
$\lambda(2n_\ell ,j,m_j)\simeq j^{2n_\ell }\centerdot N^{2n_\ell }$ since
$\lambda(2n_\ell ,j,m_j)=\det (\alpha _{2n^2_\ell ,j^2}\times D_{j^2,m^2_j})_{\rm ss}$
where $D_{j^2,m^2_j}$ is the decomposition group element of the $(j,m_j)$-th bisublattice
$(\Lambda _{\o\omega _{j,m_j}}\otimes\Lambda _{\omega _{j,m_j}})$ and where $\alpha _{2n^2_\ell ,j^2}$ is the corresponding split Cartan subgroup \cite{Pie3}.

\item On the other hand, the $(j,m_j)$-th conjugacy class functional representative\linebreak
$\phi (g^{2n_\ell }_{T_L}([2\fl],(j,m_j) ) )$
\resp{$\phi (g^{2n_\ell }_{T_R}([2\fl],(j,m_j) ) )$} is generated by means of the global Frobenius substitution:
\[ e^{2\pi iz_{n_\ell }}\To e^{2\pi ijz_{n_\ell }}
\rresp{e^{-2\pi iz_{n_\ell }}\To e^{-2\pi ijz_{n_\ell }}}\]
from the $1$-th conjugacy class functional representative
\begin{align*}
\phi (g^{2n_\ell }_{T_L}([2\fl],1 ) ) & \simeq
E_{2\fl}(2n_\ell ,1)\ \lambda ^{\half}(2n_\ell ,1)\ e^{2\pi iz_{n_\ell }}\\
\rresp{\phi (g^{2n_\ell }_{T_R}([2\fl],1 ) ) & \simeq
E_{2\fl}(2n_\ell ,1)\ \lambda ^{\half}(2n_\ell ,1)\ e^{-2\pi iz_{n_\ell }}}
\end{align*}
which is a $n_\ell $-dimensional complex semitorus shifted in $\fl$ dimensions and localized in the upper \resp{lower} half space.\qedhere
\Ee
\end{proof}
\vskip 11pt

\subsection{Proposition}

{\em The cohomology 
$H^{2n_\ell -2\fl}(
\o S^{P_{\Nfl}}_{\GL_{\Nfl}}
,
\widehat M^{2n_\ell }_{T_R}[2\fl]\otimes
\widehat M^{2n_\ell }_{T_L}[2\fl] )$ has the analytic development:
\begin{multline*}
\qquad \qquad H^{2n_\ell -2\fl}(
\o S^{P_{\Nfl}}_{\GL_{\Nfl}}
,
\widehat M^{2n_\ell }_{T_{R_\oplus}}[2\fl]\otimes
\widehat M^{2n_\ell }_{T_{L_\oplus}}[2\fl] )\\
= \L[
\bigoplus^r_{j=1}\bigoplus_{m_j}
(E_{2\fl}(2n_\ell ,j,m_j)\ \lambda ^{\half}(2n_\ell ,j,m_j)\ e^{-2\pi ijz_{n_\ell }}\R]\\
\times
\L[
\bigoplus^r_{j=1}\bigoplus_{m_j}
(E_{2\fl}(2n_\ell ,j,m_j)\ \lambda ^{\half}(2n_\ell ,j,m_j)\ e^{+2\pi ijz_{n_\ell }}\R]
\end{multline*}
according to the conjugacy class representatives
$g^{2n_\ell }_{T\RL}([2\fl],(j,m_j) ) $ where:
\begin{align*}
\EIS_L(\2nfl,(j,m_j) ) &=
\bigoplus^r_{j=1}\bigoplus_{m_j}
E_{2\fl}(2n_\ell ,j,m_j)\ \lambda ^{\half}(2n_\ell ,j,m_j)\ e^{2\pi ijz_{n_\ell }}\\
\noalign{\quad \hfill $r\le \infty\;,$}
\rresp{%
\EIS_R(\2nfl,(j,m_j) ) &=
\bigoplus^r_{j=1}\bigoplus_{m_j}
E_{2\fl}(2n_\ell ,j,m_j)\ \lambda ^{\half}(2n_\ell ,j,m_j)\ e^{-2\pi ijz_{n_\ell }}}
\end{align*}
is the (truncated) Fourier development of a normalized $2n_\ell $-dimensional \lr shifted cusp form of weight $k=2$ restricted to the upper \resp{lower} half space}
\cite{Kub} {\em (see also} \cite{Pie3}), {\em chapter3, for the introduction of the equivalent unshifted cusp form).
}
\vskip 11pt

\begin{proof}[Sketch of the proof]
This directly results from the decomposition of\linebreak
$H^{2n_\ell -2\fl}(
\o S^{P_{\Nfl}}_{\GL_{\Nfl}}
,
\widehat M^{2n_\ell }_{T_{R_\oplus}}[2\fl]\otimes
\widehat M^{2n_\ell }_{T_{L_\oplus}}[2\fl] )$ into conjugacy class functional representatives
$\phi (g^{2n_\ell }_{T\RL}([2\fl],j,m_j ) )$ whose analytic representations are given in proposition~4.7.
\end{proof}
\vskip 11pt

\subsection{Theorem (Origin of the (bilinear) spectral theorem)}

{\em The analytic development of the cohomology
$H^{2n_\ell -2\fl}(
\o S^{P_{\Nfl}}_{\GL_{\Nfl}}
,
\widehat M^{2n_\ell }_{T_{R_\oplus}}[2\fl]\otimes\linebreak
\widehat M^{2n_\ell }_{T_{L_\oplus}}[2\fl] )$ gives rise to the eigen(bi)value equation:
\begin{multline*}
\qquad \qquad (D_R^{2\fl}\otimes D_L^{2\fl})
(\EIS_R ( 2n_\ell ,j^{\rm up}=j,m_j)\otimes
(\EIS_L (2 n_\ell ,j^{\rm up}=j,m_j))\\
= E^2_{2\fl}(2n_\ell ,j,m_j)
(\EIS_R ( 2n_\ell ,j^{\rm up}=j,m_j)\otimes
(\EIS_L ( 2n_\ell ,j^{\rm up}=j,m_j))
\qquad 
\end{multline*}
where:
\Bi
\item $(D_R^{2\fl}\otimes D_L^{2\fl})$ acts on the space of smooth (unshifted) bisections
$\phi (g^{2n_\ell }_{T\RL}(j,m_j)  )$ of
$(\widehat M^{2n_\ell }_{T_{R}}\otimes
\widehat M^{2n_\ell }_{T_{L}})$ such that
$\phi (g^{2n_\ell }_{T_L}(j,m_j) ) $
\resp{$\phi (g^{2n_\ell }_{T_R}(j,m_j) ) $}
be a $C^\infty$-function localized in the upper \resp{lower} half space (see definition 3.3);

\item {\bbf\em the eigenvalues 
$E^2_{2\fl}(2n_\ell ,j,m_j)$ are the shifts of the corresponding generalized Hecke (bi)characters $\lambda (2n_\ell ,j,m_j)$;}

\item {\bbf\em the eigenbivectors
$\EIS_R ( 2n_\ell ,j^{\rm up}=j,m_j)\otimes
(\EIS_L ( 2n_\ell ,j^{\rm up}=j,m_j)$ are (tensor) products of truncated Fourier developments at the ``$j$'' classes of normalized $2n_\ell $-dimensional cusp forms, $j$ varying from $1$ to $r$.}
\Ei

{\bbf \em The set of $r$-tuples:
$\{ \EIS_R ( 2n_\ell ,1,m_1)\otimes
(\EIS_L ( 2n_\ell ,1,m_1),\dots,
\EIS_R ( 2n_\ell ,j^{\rm up}=j,m_j)\otimes
(\EIS_L ( 2n_\ell ,j^{\rm up}=j,m_j),\dots,
\EIS_R ( 2n_\ell ,j^{\rm up}=r,m_r)\otimes
(\EIS_L ( 2n_\ell ,j^{\rm up}=r,m_r)\}$ is the toroidal spectral representation of the elliptic bioperator\linebreak
$(D_R^{2\fl}\otimes D_L^{2\fl})\in \Ds_R\otimes \Ds_L$.}

The spectral measure
$\mu _{\EIS_L}$
\resp{$\mu _{\EIS_R}$} on the spectrum
$\sigma (D_L^{2\fl})$
\resp{$\sigma (D_R^{2\fl})$} of
$D_L^{2\fl}$
\resp{$D_R^{2\fl}$}
can be assumed to be the Haar measure.
}\vskip 11pt

\begin{proof} \Bena
\item The $\GL_{n_\ell }(F_{\o\omega _\oplus}\times F_{\omega _\oplus})$-bisemimodule
$(M^{2n_\ell }_{R_\oplus}\otimes M^{2n_\ell }_{L_\oplus})$ decomposes into subbisemimodules under the endomorphism
\[ E_{D_R}\otimes E_{D_L}: \quad
M^{2n_\ell }_{R_\oplus}\otimes M^{2n_\ell }_{L_\oplus}
\To\bigoplus_j\bigoplus_{m_j}
(M^{2n_\ell }_{\o\omega _{j,m_j}}\otimes M^{2n_\ell }_{\omega _{j,m_j}})\]
generated under the action of the Hecke bialgebra $\Hs\RL(n)$ according to proposition 1.10.

\item There exists a toroidal isomorphism of compactification
\[\gamma \RL: \qquad
M^{2n_\ell }_{R_\oplus}\otimes M^{2n_\ell }_{L_\oplus}
\To
M^{2n_\ell }_{T_{R_\oplus}}\otimes M^{2n_\ell }_{T_{L_\oplus}}\]
sending
$(M^{2n_\ell }_{R_\oplus}\otimes M^{2n_\ell }_{L_\oplus})$ into its toroidal equivalent
$M^{2n_\ell }_{T_{R_\oplus}}\otimes M^{2n_\ell }_{T_{L_\oplus}}$ according to section 1.12, such that:
\[ (M^{2n_\ell }_{T_{R_\oplus}}\otimes M^{2n_\ell }_{T_{L_\oplus}})
=\bigoplus_j\bigoplus_{m_j}
(M^{2n_\ell }_{T_{\o\omega _{j,m_j}}}\otimes M^{2n_\ell }_{T_{\omega _{j,m_j}}})\]
has an analytic development given by 
$\EIS_R(2n_\ell ,j,m_j)\otimes
\EIS_L(2n_\ell ,j,m_j)$ (see proposition 4.8 and \cite{Pie3}.

\item
\Bi
\item The elliptic bioperator
$(D^{2\fl}_R\otimes
D^{2\fl}_L)$ maps
$(\widehat M^{2n_\ell }_{R_\oplus}\otimes \widehat M^{2n_\ell }_{L_\oplus})$ into its shifted equivalent
$\widehat M^{2n_\ell }_{R_\oplus}[2\fl]\otimes \widehat M^{2n_\ell }_{L_\oplus}[2\fl]$ according to:
\[ D^{2\fl}_R\otimes
D^{2\fl}_L: \qquad
\widehat M^{2n_\ell }_{R_\oplus}\otimes \widehat M^{2n_\ell }_{L_\oplus}\To
\widehat M^{2n_\ell }_{R_\oplus}[2\fl]\otimes \widehat M^{2n_\ell }_{L_\oplus}[2\fl]\;.\]

\item $\widehat M^{2n_\ell }_{R_\oplus}[2\fl]\otimes \widehat M^{2n_\ell }_{L_\oplus}[2\fl]$ is transformed by the unitary action of
$(\gamma \RL\circ(E_{D_R}\otimes E_{D_L}))$ into:
\begin{multline*}
\gamma \RL\circ(E_{D_R}\otimes E_{D_L}): \quad
\widehat M^{2n_\ell }_{R_\oplus}[2\fl]\otimes \widehat M^{2n_\ell }_{L_\oplus}[2 \fl]\\
\hspace*{4cm}\To
\bigoplus_j\bigoplus_{m_j}
(\widehat M^{2n_\ell }_{T_{\o\omega _{j,m_j}}}[2\fl]\otimes \widehat M^{2n_\ell }_{T_{\omega _{j,m_j}}}[2\fl])\\
\equiv
\bigoplus_j\bigoplus_{m_j}
\phi (g^{2n_\ell }_{T\RL}([2\fl],(j,m_j) ) )\qquad \quad
\end{multline*}
in such a way that the eigenbivalue equation
\begin{multline*}
 \qquad D^{2\fl}_R\otimes
D^{2\fl}_L(\EIS_R(2n_\ell ,j^{\rm up}=j,m_j)\otimes \EIS_L(2n_\ell ,j^{\rm up}=j,m_j))\\
=E^2_{2\fl}(2n_\ell ,j,m_j)_{\rm eig}
((\EIS_R(2n_\ell ,j^{\rm up}=j,m_j)\otimes ( \EIS_L(2n_\ell ,j^{\rm up}=j,m_j))
\end{multline*}
corresponds to the map:
\[
D^{2\fl}_R\otimes
D^{2\fl}_L: \qquad
\widehat M^{2n_\ell }_{T_{R_\oplus}}\otimes \widehat M^{2n_\ell }_{T_{L_\oplus}}
\To
\widehat M^{2n_\ell }_{T_{R_\oplus}}[2\fl]\otimes \widehat M^{2n_\ell }_{T_{L_\oplus}}[2\fl ]\;.\]
\Ei

\item Indeed, according to proposition 4.6, we have that:
\[
M^{2n_\ell }_{T_R}[2\fl]\otimes   M^{2n_\ell }_{T_L}[2\fl]
\simeq \Repsp ( \GL_{\nfl} ( ( F^T_{\o\omega }\otimes\cit )\times
(F^T_\omega \otimes\cit )))\;.\]
But, $\GL_{\nfl} ( ( F^T_{\o\omega }\otimes\cit )\times
(F^T_\omega \otimes\cit ))$ is a complete solvable bilinear semigroup implying the chain of embedded normal bilinear subsemigroups:
\begin{multline*}
g^{2n_\ell }_{T\RL}([2\fl],1)\subset \dots\\
 \subset 
\dots\subset \bigoplus_{j=1}^j\bigoplus_{m_j} g^{2n_\ell }_{T\RL}([2\fl],(j,m_j))\subset \dots\\
\dots\subset \bigoplus_{j=1}^r\bigoplus_{m_j} g^{2n_\ell }_{T\RL}([2\fl],(j,m_j))\;.\end{multline*}

The analytic representation of the $j$-th normal bilinear subsemigroup:\linebreak
$\bigoplus_{j=1}^j\bigoplus_{m_j} g^{2n_\ell }_{T\RL}([2\fl],(j,m_j))$ is precisely the product
$\EIS_R(\2nfl,j^{\rm up}=j,m_j)\otimes
\EIS_L(\2nfl,j^{\rm up}=j,m_j)$ of Fourier truncated  series at ``$j$'' classes of normalized shifted $2n_\ell $-dimensional cusp forms.

So, $\EIS_R(\2nfl,j^{\rm up}=j,m_j)\otimes
\EIS_L(\2nfl,j^{\rm up}=j,m_j)$ is the $j$-th analytic representative of
$( M^{2n_\ell }_{T_{R_\oplus}}[2\fl]\otimes   M^{2n_\ell }_{T_{L_\oplus}}[2\fl ])$ and develops as follows:
\begin{multline*}
 \EIS_R(\2nfl,j^{\rm up}=j,m_j)\otimes
\EIS_L(\2nfl,j^{\rm up}=j,m_j)\\
\begin{aligned}
&=\L( \bigoplus_{j=1}^j\bigoplus_{m_j}
(E_{2\fl}(2n_\ell ,j,m_j)\ \lambda ^{\half}(2n_\ell ,j,m_j)\
e^{-2\pi ijz_{n_\ell }})\R)\\
& \qquad \qquad \otimes\L( \bigoplus_{j=1}^j\bigoplus_{m_j}
(E_{2\fl}(2n_\ell ,j,m_j)\ \lambda ^{\half}(2n_\ell ,j,m_j)\
e^{2\pi ijz_{n_\ell }})\R)\\
&= E^2_{2\fl}(2n_\ell ,j,m_j)_{\rm eig}
(\EIS_R(2n_\ell ,j^{\rm up}=j,m_j)\otimes ( \EIS_L(2n_\ell ,j^{\rm up}=j,m_j))\end{aligned}
\end{multline*}
according to proposition 4.8, where:
\Bi
\item $E_{2\fl}(2n_\ell ,j,m_j)_{\rm eig}=
\bigoplus_{j=1}^j\bigoplus_{m_j}
E_{2\fl}(2n_\ell ,j,m_j)$

\item $\EIS_L(2n_\ell ,j^{\rm up}=j,m_j)
= \bigoplus_{j=1}^j\bigoplus_{m_j}
\lambda ^{\half}(2n_\ell ,j,m_j)\ e^{+2\pi ijz_{n_\ell }}$.
\Ei

\item Thus, we have that:
\Bi
\item $E^2_{2\fl}(2n_\ell ,j,m_j)_{\rm eig}$ is the $j$-th eigenbivalue of the elliptic bioperator
$D^{2\fl}_R\otimes
D^{2\fl}_L$;

\item $(\EIS_R(2n_\ell ,j^{\rm up}=j,m_j)\otimes ( \EIS_L(2n_\ell ,j^{\rm up}=j,m_j))$ is the corresponding $j$-th eigenbifunction which is also the eigenbifunction of the product of Hecke operators $T_R(2n;j)\otimes T_L(2n;j)$ according to proposition 1.10 and \cite{Pie3}.
\Ei

And the set of $r$-tuples:
\begin{multline*}
\{
\EIS_R(2n_\ell ,j^{\rm up}=1,m_1)\otimes\EIS_L(2n_\ell ,j^{\rm up}=1,m_1),\dots,\\
\dots,\EIS_R(2n_\ell ,j^{\rm up}=j,m_j)\otimes\EIS_L(2n_\ell ,j^{\rm up}=j,m_j),\dots,\\
\dots,\EIS_R(2n_\ell ,j^{\rm up}=r,m_r)\otimes\EIS_L(2n_\ell ,j^{\rm up}=r,m_r)\}\end{multline*}
is the toroidal spectral representation of
$D^{2\fl}_R\otimes
D^{2\fl}_L$.\qedhere
\Ee
\end{proof}
\vskip 11pt

\subsection{Corollary}

{\em
Let $\widehat M^{2n_\ell }_{T_{R_\oplus}} \otimes \widehat M^{2n_\ell }_{T_{L_\oplus}} $ be the (truncated) normalized cusp biform over the
$\GL_{n_\ell }(F_{\o\omega _\oplus}\otimes F_{\omega _\oplus})$-bisemimodule
\[M^{2n_\ell }_{{R_\oplus}}\otimes   M^{2n_\ell }_{{L_\oplus}}
= \bigoplus_{j=1}^j\bigoplus_{m_j}
(M^{2n_\ell }_{\o\omega _{j,m_j}}\otimes
M^{2n_\ell }_{\omega _{j,m_j}})\]
decomposing into the sum of subbisemimodules
$(M^{2n_\ell }_{\o\omega_{j,{m_j}} }\otimes
M^{2n_\ell }_{\omega_{j,m_j}})$ according to the conjugacy classes of the complete bilinear semigroup
$\GL_{n_\ell }(F_{\o\omega }\times F_\omega )$ having multiplicities $m^{(j)}=\sup(m_j)$, $m_j$ being an increasing integer superior or equal to $1$.

Then, there exist Haar bimeasures $\mu _{j_R}\times\mu _{j_L}$ on the spectrum $\sigma (D^{2\fl}_R\otimes
D^{2\fl}_L)$ of the elliptic bioperator
$(D^{2\fl}_R\otimes
D^{2\fl}_L)$ and an isomorphism
\begin{multline*}
 \qquad \gamma \RL\circ (E_{D_R}\otimes E_{D_L}): \qquad
\widehat M^{2n_\ell }_{{R_\oplus}}\otimes \widehat M^{2n_\ell }_{{L_\oplus}}\\
\To \bigoplus_{j=1}^j\bigoplus_{m_j} (
( \EIS_R(2n_\ell ,j,m_j )\otimes\EIS_L ( 2n_\ell ,j,m_j ))\qquad 
\end{multline*}
leading to the eigenbivalue equation(s):
\[
[(\gamma _R\circ D^{2\fl}_R)\otimes
(D^{2\fl}_L\circ\gamma _L)(
\widehat M^{2n_\ell }_{R_\oplus}\otimes
\widehat M^{2n_\ell }_{L_\oplus}]=
E^2_{2\fl}(2n_\ell ,j,m_j)_{\rm eig} 
(\widehat M^{2n_\ell }_{R_\oplus}\otimes \widehat M^{2n_\ell }_{L_\oplus})
\] whose spectral representation is given by the set of eigenbifunctions
$\{\EIS_R(2n_\ell ,j^{\rm up}=j,m_j)\otimes
\EIS_L(2n_\ell ,j^{\rm up}=j,m_j)\}^r_{j=1,m_j}$
having multiplicities $m^{(j)}=\sup(m_j)$.
}
\vskip 11pt

\subsection{Proposition}

{\em
{\bfseries A trace formula}}  \cite{Art} {\bbf \em corresponding to a shifted Plancherel formula}  and {\em associated with the $j$-th eigenbifunction 
$\EIS_R(2n_\ell ,j^{\rm up}=j,m_j)\otimes
\EIS_L(2n_\ell ,j^{\rm up}=j,m_j)$ of the eigenvalue equation:
\[
[(\gamma _R\circ D^{2\fl}_R)\otimes
(D^{2\fl}_L\circ\gamma _L)(
\widehat M^{2n_\ell }_{R_\oplus}\otimes
\widehat M^{2n_\ell }_{L_\oplus}]
=
E^2_{2\fl}(2n_\ell ,j,m_j)_{\rm eig} 
(\widehat M^{2n_\ell }_{R_\oplus}\otimes \widehat M^{2n_\ell }_{L_\oplus})
\] 
is given by the {\bfseries bilinear form}:
\begin{multline*}
 (\EIS_R(\2nfl ,j^{\rm up}=j,m_j),
\EIS_L(\2nfl ,j^{\rm up}=j,m_j))\\
= \bigoplus_{j=1}^j\bigoplus_{m_j}
(\lambda (2n_\ell ,j,m_j)\ E^2_{2\fl}
(2n_\ell ,j,m_j))\end{multline*}
from $\FRepsp ( \GL_{\nfl} (( F^T_{\o\omega _\oplus}\otimes\cit ) \times
(F^T_{\omega _\oplus}\otimes\cit)))$ to $\cit$.
}
\vskip 11pt

\begin{proof}
 This trace formula directly results from point 4) of the proof of proposition 4.9 and corresponds to the shifted Plancherel formula since the trace formula
\[
(\EIS_R(2n_\ell  ,j,m_j),
\EIS_L(2n_\ell  ,j,m_j))
= \bigoplus_{j=1}^r
\lambda  (2n_\ell ,j,m_j)\]
from $\FRepsp ( \GL_{n_\ell } ( F^T_{\o\omega _\oplus}\otimes F^T_{\omega _\oplus}))$ to $\cit$ is the
{\bbf Plancherel formula} associated with the bilinear semigroup
$ \GL_{n_\ell } ( F^T_{\o\omega } \times
F^T_{\omega} )$.
\end{proof}
\vskip 11pt

\subsection{Proposition}

{\em
The product, right by left, 
$\EIS_R(\2nfl,j,m_j)\otimes \EIS_L(\2nfl,j,m_j)$ of the truncated Fourier development of the shifted $2n_\ell$ dimensional cusp biform, constitutes:
\Bena
\item a supercuspidal representation of the shifted complete bilinear semigroup\linebreak
$\GL_{\nfl}((F^T_{\o\omega }\otimes\cit)\times
(F^T_{\omega }\otimes\cit))$;

\item a shifted supercuspidal representation of the complete bilinear semigroup\linebreak
$\GL_{n_\ell }((F^T_{\o\omega }\times F^T_{\omega})$.
\Ee
}
\vskip 11pt

\begin{proof}
\Bena
\item According to proposition 4.8,
$\EIS_R(\2nfl,j,m_j)$ \resp{$\EIS_L(\2nfl,j,m_j)$} is the truncated Fourier development of a normalized $2n_\ell $-dimensional \rl shifted cusp form of weight $k=2$.  Consequently,
$\EIS\RL(\2nfl,j,m_j)=
\EIS_R(\2nfl,j,m_j)\otimes \EIS_L(\2nfl,j,m_j)$ is a truncated cuspidal biform over 
$\cit^{n_\ell }\times_D\cit^{n_\ell }$.  
On the other hand, as we have the equality:
\[ \FRepsp ( \GL_{\nfl}(( F^T_{\o\omega _\oplus}\otimes\cit )\times
( F^T_{\omega _\oplus}\otimes\cit ) ) )=
\EIS\RL(\2nfl,j,m_j)\]
according to propositions 4.6 and 4.7, and as
$\FRepsp(\GL_{\nfl}((F^T_{\o\omega  }\otimes\cit)\times
(F^T_{\omega  }\otimes\cit)))$ is irreducible and the coefficients of
$\EIS\RL(\2nfl,j,m_j)$ have compact support in
$\GL_{\nfl}((F^T_{\o\omega}\otimes\cit)\times
(F^T_{\omega }\otimes\cit))$,
$\EIS\RL(\2nfl,j,m_j)$ constitutes a supercuspidal representation of
$\GL_{\nfl}((F^T_{\o\omega  }\otimes\cit)\times
(F^T_{\omega  }\otimes\cit))$.

\item Taking into account proposition 4.9, it clearly appears that $\EIS\RL(\2nfl,j,m_j)$ also constitutes a shifted supercuspidal representation of
$\GL_{n_\ell }(F^T_{\o\omega }\times F^T_{\omega })$.\qedhere
\Ee
\end{proof}
\vskip 11pt

\subsection{Holomorphic spectral representation}

This chapter has been essentially devoted to the toroidal spectral representation of an elliptic bioperator given by a set of $r$-tuples of products, right by left, of truncated Fourier developments of cusp forms.

Indeed, the aim of this paper, put in concrete form in chapter 5, deals with supercuspidal representations of shifted algebraic bilinear semigroups in the frame of geometric-shifted global bilinear correspondences of Langlands.

A cusp form, being a holomorphic function, the conclusions obtained for the toroidal spectral representation of an elliptic bioperator result in fact from its {\bbf ``holomorphic'' spectral representation} as developed succinctly in the next sections of this chapter.

As the toroidal spectral representation of an elliptic bioperator $(D^{2\fl}_R\otimes D^{2\fl}_L)$ is directly connected to the functional representation space
$\FRepsp ( \GL_{\nfl}(( F^T_{\o\omega }\otimes\cit )\times
( F^T_{\omega }\otimes\cit ) ) )$ of the shifted bilinear complete semigroup
$ \GL_{\nfl}(( F^T_{\o\omega  }\otimes\cit )\times
( F^T_{\omega  }\otimes\cit ) ) $, the associated holomorphic spectral representation will be proved to result from the functional representation space
$\FRepsp ( \GL_{\nfl}(( F_{\o\omega }\otimes\cit )\times
( F_{\omega }\otimes\cit ) ) )$
of the shifted bilinear complete semigroup
$ \GL_{\nfl}(( F_{\o\omega }\otimes\cit )\times
( F_{\omega }\otimes\cit ) ) )$.
\vskip 11pt

\subsection{Proposition}

{\em
The differential bioperator $(D^{2\fl}_R\otimes D^{2\fl}_L)\in\Ds_R\otimes\Ds_L$  maps the bisemisheaf 
$(\widehat M^{2n_\ell }_{R}\otimes \widehat M^{2n_\ell }_{L})$ on the
$ \GL_{n_\ell }( F_{\o\omega }\times F_{\omega } )$-bisemimodule $( M^{2n_\ell }_{R}\otimes  M^{2n_\ell }_{L})$ into the perverse bisemisheaf
$(\widehat M^{2n_\ell }_{R}[2\fl]\otimes \widehat M^{2n_\ell }_{L}[2\fl])$ according to:
\[ (D^{2\fl}_R\otimes D^{2\fl}_L)
: \qquad (\widehat M^{2n_\ell }_{R}\otimes \widehat M^{2n_\ell }_{L})
\To(\widehat M^{2n_\ell }_{R}[2\fl]\otimes \widehat M^{2n_\ell }_{L}[2\fl])\;.\]
}
\vskip 11pt

\begin{proof}
This is an adaptation of proposition 4.3.\end{proof}
\vskip 11pt

\subsection{Proposition}

{\em 
The bilinear cohomology
$H^{2n_\ell -2\fl}(M_{DM\RL})(X^{\rm sv}\RL),X^{2n_\ell }_R[2\fl]\times X^{2n_\ell }_L[2\fl])$ of the \SV\ mixed bimotive
$M_{DM\RL}(X^{\rm sv}\RL)$ is isomorphic to the decomposition in conjugacy classes of the product, right by left,
${\rm CY} ^{2n_\ell }(Y_R,[2\fl])\times
{\rm CY} ^{2n_\ell }(Y_L,[2\fl])$ of $2n_\ell $-dimensional cycles shifted in $2\fl$-dimensions:
\begin{multline*}
H^{2n_\ell -2\fl}(M_{DM\RL})(X^{\rm sv}\RL),X^{2n_\ell }_R[2\fl]\times X^{2n_\ell }_L[2\fl])\\
\begin{aligned}
&\simeq {\rm CY} ^{2n_\ell }(Y_R,[2\fl])\times
{\rm CY} ^{2n_\ell }(Y_L,[2\fl])\\
&= \{ {\rm CY} ^{2n_\ell }(Y_R,[2\fl],(j,m_j))\times
{\rm CY} ^{2n_\ell }(Y_L,[2\fl],(j,m_j))\}_{j,m_j}\;.\end{aligned}
\end{multline*}
}
\vskip 11pt

\begin{proof}
This results from the isomorphism
\begin{multline*}
\qquad H^{2n_\ell -2\fl}(M_{DM\RL})(X^{\rm sv}\RL),X^{2n_\ell }_R[2\fl]\times X^{2n_\ell }_L[2\fl])\\
\simeq
H^{2n_\ell -2\fl}
( \widehat M^{2n}_R[2\fl]\otimes \widehat M^{2n}_L[2\fl] ),
( \widehat M^{2n_\ell }_R[2\fl]\otimes \widehat M^{2n_\ell }_L[2\fl] )\qquad 
\end{multline*}
between cohomologies according to propositions 4.4 and 4.5 in such a way that
\Bean
\item $\widehat M^{2n_\ell }_R[2\fl]\otimes \widehat M^{2n_\ell }_L[2\fl] )$ is the bisemisheaf over the representation space\linebreak
$\Repsp ( \GL_{\nfl}(( F_{\o\omega }\otimes\cit )\times
( F_{\omega }\otimes\cit ) ) )$
of the bilinear complete semigroup\linebreak
$ \GL_{\nfl}(( F_{\o\omega }\otimes\cit )\times
( F_{\omega }\otimes\cit ) ) $;

\item ${\rm CY} ^{2n_\ell }(Y_R,[2\fl])\times
{\rm CY} ^{2n_\ell }(Y_L,[2\fl])$ is the $2n_\ell $-th bicycle isomorphic to\linebreak $ \GL_{\nfl}(( F_{\o\omega }\otimes\cit )\times
( F_{\omega }\otimes\cit ) ) $.\qedhere
\Ee
\end{proof}
\vskip 11pt

\subsection{\bbf Laurent polynomials on 
$ \GL_{n_\ell }(F_{\o\omega }\times
 F_{\omega } ) $}
 
 \Bi
 \item Let $g^{2n_\ell }_L(j,m_j)$
\resp{$g^{2n_\ell }_R(j,m_j)$} be the $(j,m_j)$-th \lr linear conjugacy class representative of
$\GL_{n_\ell }(F_{\o\omega }\times
 F_{\omega } ) $
and let 
$\psi (g^{2n_\ell }_L(j,m_j))$
\resp{$\psi (g^{2n_\ell }_R(j,m_j))$} be a differentiable function on it, into $\cit$, given simply by
\begin{align*}
\psi (g^{2n_\ell }_L(j,m_j))
&= \lambda ^{\half}(2n_\ell ,j,m_j)(y^j_1\times\dots\times y^j_{n_\ell })\;, &&   \\
&= \lambda ^{\half}(2n_\ell ,j,m_j)y^j\;, && y=y_1\times\dots\times y_{n_\ell }\;, \\
\rresp{
\psi (g^{2n_\ell }_R(j,m_j))
&= \lambda ^{\half}(2n_\ell ,j,m_j)(y^{*j}_1\times\dots\times y^{*j}_{n_\ell })\;, && y_{n_\ell }^*\text{\ being the}\\
\noalign{\mbox{}\hfill  conjugate complex of $y_{n_\ell }$}
&= \lambda ^{\half}(2n_\ell ,j,m_j)(y^*)^j\;,}
\end{align*}
where $y_1,\dots,y_{n_\ell}$ are functions of complex variables on unitary closed supports.

\item If the conjugacy class representatives of
$T_{n_\ell }(F_\omega )$
\resp{$T^t_{n_\ell }(F_{\o\omega })$} $\subset \GL_{n_\ell }(F_{\o\omega }\times F_\omega )$ are glued together, {\bbf a Laurent polynomial on the representation space of 
$T_{n_\ell }(F_{\omega _\oplus})$
\resp{$T^t_{n_\ell }(F_{\o\omega _\oplus})$}}
will be defined by:
\begin{align*}
\psi (\Repsp(T_{n_\ell }(F_{\omega _\oplus}))
&=\sum_{j=1}^r \sum_{m_j} \lambda ^{\half}
(2n_\ell ,j,m_j)\ y^j\;, && r\le\infty \\
\rresp{
\psi (\Repsp(T^t_{n_\ell }(F_{\o\omega _\oplus}))
&=\sum_{j=1}^r \sum_{m_j} \lambda ^{\half}
(2n_\ell ,j,m_j)\ (y^*)^j\;,}\end{align*}
where $\lambda ^{\half}
(2n_\ell ,j,m_j)$ is the square root of the product of the eigenvalues of the $(j,m_j)$-th coset representative of the Hecke bioperator as described in proposition 4.7.

\item And, a Laurent bipolynomial on the representation space 
$\GL_{n_\ell }(F_{\o\omega _\oplus}\times F_{\omega _\oplus})$ with respect to its conjugacy classes glued together will be given by:
\[ \psi ( \Repsp ( \GL_{n_\ell } ( F_{\o\omega _\oplus}\times F_{\omega _\oplus})))
= \sum_{j=1}^j\sum_{m_j} (\lambda ^{\half}(2n_\ell ,j,m_j) \times (y^*)^j)\times
 (\lambda ^{\half}(2n_\ell ,j,m_j) \times y^j)\;.\]
\Ei
\vskip 11pt

\subsection{Proposition}

{\em On the glued together conjugacy class representatives of 
$\GL_{n_\ell }(F_{\o\omega }\times F_\omega )$, the function
$\psi (\Repsp(T_{n_\ell }(F_{\omega _\oplus})))$
\resp{$\psi (\Repsp(T^t_{n_\ell }(F_{\o\omega _\oplus})))$}, defined in a neighbourhood of a point $y'_0$ \resp{$y_0^*{}'$} of $\cit^{n_\ell }$, is holomorphic at $y'_0 $ \resp{$y_0^*{}'$} if we have the following multiple power series developments:
\begin{align*}
\psi (\Repsp(T_{n_\ell }(F_{\omega _\oplus}))
&= \sum_{j=1}^r\sum_{m_j} \lambda ^{'\half}\ (2n_\ell ,j,m_j)\ (y'-y'_0)^j\\
\rresp{%
\psi (\Repsp(T^t_{n_\ell }(F_{\o\omega _\oplus}))
&= \sum_{j=1}^r\sum_{m_j} \lambda ^{'\half}\ (2n_\ell ,j,m_j)\ (y^{*'}-y^{*'}_0)^j}.
\end{align*}
And, the holomorphic bifunction
\[
\psi (\Repsp(T^t_{n_\ell }(F_{\o\omega _\oplus}))\otimes
\psi (\Repsp(T_{n_\ell }(F_{\omega _\oplus}))
= \sum_{j=1}^r\sum_{m_j} \lambda' \ (2n_\ell ,j,m_j)\ ( y^{*'}y'-y^{*'}_0y'_0)^j \]
at the bipoint $(y_0^*{}'y'_0)$ constitutes an irreducible holomorphic representation\linebreak $\Irr\hol ((\GL_{n_\ell }(F_{\o\omega }\times F_\omega )))$ of the bilinear  semigroup $\GL_{n_\ell }(F_{\o\omega }\times F_\omega )$.
}
\vskip 11pt

\begin{proof}[Sketch of proof]
This is a consequence of the introduction of Laurent polynomials on $\GL_{n_\ell }(F_{\o\omega }\times F_\omega )$ in section 4.16 and \cite{Pie3}.\end{proof}
\vskip 11pt

\subsection[Shifted holomorphic representation of
$\GL_{\nfl }((F_{\o\omega }\otimes\cit)\times (F_\omega \otimes\cit))$]{\bbf Shifted holomorphic representation of
$\GL_{\nfl }((F_{\o\omega }\otimes\cit)\times (F_\omega \otimes\cit))$}

\Bi

\item Similarly as it was done in proposition 4.7, a function
$\psi (\Repsp(T_{\nfl}(F_{\omega  _\oplus}\otimes\cit)))$
\resp{a cofunction $\psi (\Repsp(T^t_{\nfl}(F_{\o\omega _\oplus}\otimes\cit)))$} on the representation space of the linear complete semigroup
$T_{\nfl}(F_{\omega _\oplus}\otimes\cit))$
\resp{$T^t_{\nfl}(F_{\o\omega _\oplus}\otimes\cit))$}, shifted in $2\fl$ dimensions, will be introduced by:
\begin{align*}
&\psi (\Repsp(T_{\nfl}(F_{\omega _\oplus}\otimes\cit))))
\qquad \qquad \\ &\qquad \qquad = \sum_{j=1}^r\sum_{m_j} c_{2\fl}\ (2n_\ell ,j,m_j)\ \lambda ^{\half} (2n_\ell ,j,m_j)\ (y^{j})\\
\rresp{%
&\psi (\Repsp(T^t_{\nfl}(F_{\o\omega _\oplus}\otimes\cit))))
\qquad \qquad \\ &\qquad \qquad = \sum_{j=1}^r\sum_{m_j} c_{2\fl}\ (2n_\ell ,j,m_j)\ \lambda ^{\half} (2n_\ell ,j,m_j)\ (y^{*j})}\end{align*}
where:
 $c_{2\fl}\ (2n_\ell ,j,m_j)$ is the shift in $\fl$ complex dimensions of the Hecke character $\lambda ^{\half}(2n_\ell ,j,m_j)$.

\item And the bifunction
\begin{multline*}
\qquad \qquad \psi ( \Repsp ( T^t_{\nfl} ( F_{\o\omega _\oplus}\otimes\cit )))\otimes
\psi ( \Repsp ( T_{\nfl} ( F_{\omega _\oplus}\otimes\cit )))\\
=\sum_{j=1}^r\sum_{m_j} c^2_{2\fl}\ (2n_\ell ,j,m_j)\ \lambda  (2n_\ell ,j,m_j)\ (y^*y)^{j}
\end{multline*}
on the representation space of the bilinear complete semigroup $\GL_{\nfl }(F_{\o\omega_\oplus }\otimes\cit)\times (F_{\omega_\oplus} \otimes\cit))$, shifted in $\fl$ complex dimensions on its right and left parts, constitutes an irreducible shifted holomorphic representation
$\Irr\hol(\GL_{\nfl }(F_{\o\omega }\otimes\cit)\times (F_\omega \otimes\cit)))$ of
$\GL_{n_\ell  }(F_{\o\omega }\otimes\cit)\times (F_\omega \otimes\cit))$.
\Ei
\vskip 11pt

\subsection{Theorem (Holomorphic spectral theorem)}

{\em The analytic development of the cohomology
\begin{multline*}
H^{2n_\ell  -2\fl}
(\widehat M^{2n}_R[2\fl]\otimes \widehat M^{2n}_L[2\fl],
\widehat M^{2n_\ell }_R[2\fl]\otimes \widehat M^{2n_\ell }_L[2\fl])
\\
=\sum_{j=1}^r\sum_{m_j}
(c^2_{2\fl}(2n_\ell ,j,m_j)\
\lambda (2n_\ell ,j,m_j)\ (y^*y)^{j}
\end{multline*}
gives rise to the eigen(bi)value equation:
\begin{multline*}
(D^{2\fl}_R\otimes D^{2\fl}_L)
(\psi ( \Repsp ( T^t_{n_\ell } ( F_{\o\omega _\oplus}
,j^{\rm up}=j))
\otimes
\psi ( \Repsp ( T_{n_\ell } ( F_{\omega _\oplus},j^{\rm up}=j)))\\
= c_{2\fl}(2n_\ell ,j,m_j)_{\rm eig}\ 
(\psi ( \Repsp ( T^t_{n_\ell } ( F_{\o\omega _\oplus}
,j^{\rm up}=j))
\otimes
\psi ( \Repsp ( T_{n_\ell } ( F_{\omega _\oplus},j^{\rm up}=j)))
\end{multline*}
where:
\Bi
\item $(D^{2\fl}_R\otimes D^{2\fl}_L)$ acts on the space of smooth (unshifted) bisections
$\psi (g^{2n_\ell }\RL(j,m_j))$ of $(\widehat M^{2n_\ell }_R\otimes\widehat M^{2n_\ell }_L)$;

\item {\bbf \em the eigenbivectors
$(\psi ( \Repsp ( T^t_{n_\ell } ( F_{\o\omega _\oplus}
,j^{\rm up}=j))
\otimes
\psi ( \Repsp ( T_{n_\ell } ( F_{\omega _\oplus},j^{\rm up}=j)))$ are tensor products of truncated holomorphic functions at $j^{\rm up}=j$ terms in such a way that the $r$-tuple:}
\begin{multline*}
\bigl\{
(\psi ( \Repsp ( T^t_{n_\ell } ( F_{\o\omega _\oplus}
,j^{\rm up}=1))
\otimes
\psi ( \Repsp ( T_{n_\ell } ( F_{\omega _\oplus},j^{\rm up}=1)))\bigr.\\
\begin{aligned}
&,\dots,\\
&(\psi ( \Repsp ( T^t_{n_\ell } ( F_{\o\omega _\oplus}
,j^{\rm up}=j))
\otimes
\psi ( \Repsp ( T_{n_\ell } ( F_{\omega _\oplus},j^{\rm up}=j)))\\
&,\dots,\end{aligned}\\
\bigl.(\psi ( \Repsp ( T^t_{n_\ell } ( F_{\o\omega _\oplus}
,j^{\rm up}=r))
\otimes
\psi ( \Repsp ( T_{n_\ell } ( F_{\omega _\oplus},j^{\rm up}=r)))\bigr\}
\end{multline*}
{\bbf \em constitutes the holomorphic spectral representation of the elliptic bioperator
$(D^{2\fl}_R\otimes D^{2\fl}_L)$;}

\item {\bbf \em the eigenbivalues $c_{2\fl}(2n_\ell ,j,m_j)$ are shifts of the generalized Hecke\linebreak (bi)characters in the sense of section 4.18.}
\Ei
}
\vskip 11pt

\begin{proof}[Sketch of proof]
This theorem is an adaptation to the holomorphic case of the spectral theorem 4.9 having led to a toroidal spectral representation of the elliptic bioperator
$(D^{2\fl}_R\otimes D^{2\fl}_L)$.
\end{proof}
\vskip 11pt

\subsection{Proposition}

{\em
\Bena
\item The holomorphic spectral representation of theorem 4.19 is isomorphic to the toroidal spectral representation of theorem 4.9.

\item Every spectral representation on the representation space of the bilinear complete semigroup
$\GL_{n_\ell }(F_{\o\omega }\times F_\omega )$ is isomorphic to the holomorphic and toroidal spectral representations mentioned above.
\Ee
}
\vskip 11pt

\begin{proof}
Indeed, every spectral representation on a functional representation space\linebreak
$F(\Repsp(\GL_{n_\ell }(F_{\o\omega }\times F_\omega ))$ of the bilinear complete semigroup
$\GL_{n_\ell }(F_{\o\omega }\times F_\omega )$ has the structure of a $r$-tuple:

\begin{multline*}
\bigl\{
(F( \Repsp ( T^t_{n_\ell } ( F_{\o\omega _\oplus}
,j^{\rm up}=1))
\otimes
F ( \Repsp ( T_{n_\ell } ( F_{\omega _\oplus},j^{\rm up}=1)))\bigr.\\
\begin{aligned}
&,\dots,\\
&(F ( \Repsp ( T^t_{n_\ell } ( F_{\o\omega _\oplus}
,j^{\rm up}=j))
\otimes
F ( \Repsp ( T_{n_\ell } ( F_{\omega _\oplus},j^{\rm up}=j)))\\
&,\dots,\end{aligned}\\
\bigl.(F( \Repsp ( T^t_{n_\ell } ( F_{\o\omega _\oplus}
,j^{\rm up}=r))
\otimes
F ( \Repsp ( T_{n_\ell } ( F_{\omega _\oplus},j^{\rm up}=r)))\bigr\}\\
1\le j\le r\le\infty 
\end{multline*}
as indicated for the holomorphic and toroidal spectral representations of the elliptic bioperator
$(D^{2\fl}_R\otimes D^{2\fl}_L)$.
\end{proof}

%% file: LanglandsIII_5.tex
\section{Geometric-shifted global bilinear correspondences of Langlands}

\subsection{Lemma}

{\em
Let $W^{ab}_{F^{S_{\cit}}_{\o\omega }}\times
W^{ab}_{F^{S_{\cit}}_{\omega }}$ be the product of the shifted global Weil (semi)group
$W^{ab}_{F^{S_{\cit}}_{\o\omega }}$ by its equivalent
$W^{ab}_{F^{S_{\cit}}_{\omega }}$ as introduced in section 1.5.

Then, there exists an irreducible representation:
\[ 
\Irr W^{(2n_\ell )}_{F\RL}: \qquad
W^{ab}_{F^{S_{\cit}}_{\o\omega }}
\times W^{ab}_{F^{S_{\cit}}_{\omega }}
\To \GL_{\nfl}((
F_{\o\omega  }\otimes \cit)\times
(F_{\omega  s}\otimes \cit))\]
from 
$(W^{ab}_{F^{S_{\cit}}_{\o\omega }}\times
W^{ab}_{F^{S_{\cit}}_{\omega }})$ to the complex bilinear complete semigroup 
$\GL_{\nfl}((
F_{\o\omega  }\otimes \cit)\times
(F_{\omega  }\otimes \cit))$ shifted in $\fl$ complex dimensions in such a way that {\em \cite{Pie1}, \cite{Pie2}:}
\Bena
\item $G^{(\2nfl)}((
F_{\o\omega  }\otimes \cit)\times
(F_{\omega  }\otimes \cit))
\simeq \GL_{\nfl}((
F_{\o\omega  }\otimes \cit)\times
(F_{\omega }\otimes \cit))$ where\linebreak
$G^{(\2nfl)}((
F_{\o\omega  }\otimes \cit)\times
(F_{\omega  }\otimes \cit))$ is a condensed notation for the shifted bilinear complete semigroup
$M^{2n_\ell }_R[2\fl]\otimes M^{2n_\ell }_L[2\fl]$;

\item $\Irr\Rep^{(\2nfl)}_{W_{F\RL}}
(W^{ab}_{F^{S_{\cit}}_{\o\omega }}\times
W^{ab}_{F^{S_{\cit}}_{\omega }})
=G^{(\2nfl)}((
F_{\o\omega _\oplus}\otimes \cit)\times
(F_{\omega _\oplus}\otimes \cit))$ where\linebreak
$\Irr\Rep^{(\2nfl)}_{W_{F\RL}}
(W^{ab}_{F^{S_{\cit}}_{\o\omega }}\times
W^{ab}_{F^{S_{\cit}}_{\omega }})$ is the irreducible $2n_\ell $-dimensional shifted global Weil-Deligne representation of
$
(W^{ab}_{F^{S_{\cit}}_{\o\omega }}\times
W^{ab}_{F^{S_{\cit}}_{\omega }})$.
\Ee}
\vskip 11pt

\begin{proof}
\Bi
\item The shifted bilinear complete semigroup
$G^{(\2nfl)}((
F_{\o\omega  }\otimes \cit)\times
(F_{\omega  }\otimes \cit))
=\Repsp( \GL_{\nfl}((
F_{\o\omega  }\otimes \cit)\times
(F_{\omega }\otimes \cit)))$, isomorphic to the bilinear semigroup of matrices
$\GL_{\nfl}((
F_{\o\omega }\otimes \cit)\times
(F_{\omega  }\otimes \cit))$, 
is a
$\GL_{\nfl}((
F_{\o\omega  }\otimes \cit)\times
(F_{\omega  }\otimes \cit))
$-bisemimodule
$M^{2n_\ell }_R[2\fl]\otimes M^{2n_\ell }_L[2\fl]$.

\item The   representation of 
$\GL_{\nfl}((
F_{\o\omega  }\otimes \cit)\times
(F_{\omega  }\otimes \cit))
$ into 
$M^{2n_\ell }_R[2\fl]\otimes M^{2n_\ell }_L[2\fl]$ corresponds to an algebraic morphism from
$\GL_{\nfl}((
F_{\o\omega  }\otimes \cit)\times
(F_{\omega  }\otimes \cit))$
into
$\GL(M^{2n_\ell }_R[2\fl]\otimes M^{2n_\ell }_L[2\fl])$
which denotes the group of automorphisms of
$M^{2n_\ell }_R[2\fl]\otimes M^{2n_\ell }_L[2\fl]$.

So, $\GL(M^{2n_\ell }_R[2\fl]\otimes M^{2n_\ell }_L[2\fl])$
constitutes the $2n_\ell $-dimensional equivalent of the product
$(W^{ab}_{F^{S_{\cit}}_{\o\omega }}\times
W^{ab}_{F^{S_{\cit}}_{\omega }})$ of shifted global Weil groups.

\item As $\GL(M^{2n_\ell }_R[2\fl]\otimes M^{2n_\ell }_L[2\fl])$ is isomorphic to
$\GL_{\nfl}((
F_{\o\omega }\otimes \cit)\times
(F_{\omega  }\otimes \cit))$, the shifted bilinear   semigroup
$G^{(\2nfl)}((
F_{\o\omega _\oplus}\otimes \cit)\times
(F_{\omega _\oplus}\otimes \cit))$ becomes the natural irreducible $2n_\ell $-dimensional shifted global Weil-Deligne representation of
$(W^{ab}_{F^{S_{\cit}}_{\o\omega }}\times
W^{ab}_{F^{S_{\cit}}_{\omega }})$.\qedhere
\Ei
\end{proof}
\vskip 11pt

\subsection{Proposition}

{\em
On the shifted bilinear complete semigroup
$G^{(\2nfl)}((
F_{\o\omega  }\otimes \cit)\times
(F_{\omega  }\otimes \cit))
$
there exists the {\bbf \em geometric-shifted global bilinear correspondence of Langlands:}
{ \begin{small}
\[\begin{psmatrix}[colsep=.5cm,rowsep=1cm]
\IrrRep^{(2n_\ell [2\fl])}_{W_{F\RL}}(W^{ab}_{F^{S_{\cit}}_{\o\omega }}\times
W^{ab}_{F^{S_{\cit}}_{\omega }}) 
& \raisebox{3mm}{$\sim$}
&
\Irrcusp(\GL_{n_\ell [\fl]}(({F^T_{\o\omega }}\otimes\cit)\times
({F^T_{\omega }}\otimes\cit)))
 \\
G^{(2n_\ell [2\fl])}((F_{\o\omega _{\bigoplus}}\otimes\cit)\times
(F_{\omega _{\bigoplus}}\otimes\cit))
 & &
\EIS\RL(2n_\ell [2\fl],j,m_j)
\\
 && G^{(2n_\ell [2\fl])}
( ( {F^T_{\o\omega }}\otimes\cit )\times
( {F^T_{\omega }}\otimes\cit ) )
\\
\noalign{\vspace*{-8mm}}
&& \rotatebox{90}{\Huge{$\approx$}}\\
\noalign{\vspace*{-8mm}}
&& H^{2n_\ell -2\fl}(\o S^{P_n[\fl]}_{\GL_n[\fl]}, 
\widehat M^{2n_\ell }_{T_{R_\oplus}}[2\fl]\otimes
\widehat M^{2n_\ell }_{T_{L_\oplus}}[2\fl])\\
\noalign{\vspace*{-8mm}}
&& \rotatebox{90}{\Huge{$\approx$}}\\
\noalign{\vspace*{-8mm}}
&& {\rm CY}^{2n_\ell }_{\oplus}(Y_R,[2\fl])\times
{\rm CY}^{2n_\ell }_{\oplus}(Y_L,[2\fl])
\ncline[arrows=->,nodesep=5pt]{1,1}{1,3}
\ncline[doubleline=true,nodesep=5pt]{1,1}{2,1}
\ncline[doubleline=true,nodesep=5pt]{1,3}{2,3}
\ncline[arrows=->,nodesep=5pt]{2,1}{3,3}^{\rotatebox{-30}{$\sim$}}
\ncline[arrows=->,nodesep=5pt]{3,3}{2,3}<{\rotatebox{90}{$\sim$}}
\end{psmatrix}\]
\end{small}}
\Bi
\item from the (infinite) sum of products, right by left, of the equivalence classes of the irreducible $2n_\ell $-dimensional shifted global Weil-Deligne representation\linebreak 
$\IrrRep^{(2n_\ell [2\fl])}_{W_{F\RL}}(W^{ab}_{F^{S_{\cit}}_{\o\omega }}\times
W^{ab}_{F^{S_{\cit}}_{\omega }}) $ of the shifted bilinear global Weil group
$W^{ab}_{F^{S_{\cit}}_{\o\omega }}\times
W^{ab}_{F^{S_{\cit}}_{\omega }}$ given by the shifted bilinear complete semigroup
$G^{(2n_\ell [2\fl])}((F_{\o\omega _{\bigoplus}}\otimes\cit)\times
(F_{\omega _{\bigoplus}}\otimes\cit))$
\

\item to the shifted irreducible supercuspidal representation
$\Irrcusp(\GL_{n_\ell [\fl]}(({F^T_{\o\omega }}\otimes\cit)\times
({F^T_{\omega }}\otimes\cit)))$ of
$\GL_{n_\ell [\fl]}(({F^T_{\o\omega }}\otimes\cit)\times
({F^T_{\omega }}\otimes\cit))$
given by the $2n_\ell $-dimensional solvable truncated shifted Eisenstein biserie
$\EIS\RL(2n_\ell [2\fl],j,m_j)$

\item which are in one-to-one correspondence with the (infinite) sum of the products, right by left, of the equivalence classes of the   irreducible $2n_\ell $-dimensional shifted cycle representatives
${\rm CY}^{2n_\ell }_{\oplus}(Y_R,[2\fl])\times
{\rm CY}^{2n_\ell }_{\oplus}(Y_L,[2\fl])$.
\Ei
}
\vskip 11pt

\begin{proof}
\Bi
\item This proposition is an adaptation of proposition 3.4.14 of \cite{Pie3} to the shifted case.

\item In lemma 5.1, it was proved that
$\IrrRep^{(2n_\ell [2\fl])}_{W_{F\RL}}(W^{ab}_{F^{S_{\cit}}_{\o\omega }}\times
W^{ab}_{F^{S_{\cit}}_{\omega }}) 
= G^{(\2nfl)}((
F_{\o\omega _\oplus}\otimes \cit)\times
(F_{\omega _\oplus}\otimes \cit))$ in such a way that:
\Bean
\item
$G^{(\2nfl)}((
F^T_{\o\omega _\oplus}\otimes \cit)\times
(F^T_{\omega _\oplus}\otimes \cit))\approx
G^{(\2nfl)}((
F_{\o\omega _\oplus}\otimes \cit)\times
(F_{\omega _\oplus}\otimes \cit))$ as mentioned in the proof of proposition 4.9.

\item $G^{(\2nfl)}((
F^T_{\o\omega _\oplus}\otimes \cit)\times
(F^T_{\omega _\oplus}\otimes \cit))
=( M^{2n_\ell }_{T_{R_\oplus}}[2\fl]\otimes
  M^{2n_\ell }_{T_{L_\oplus}}[2\fl])$ has an analytical development given by
  $\EIS\RL(2n_\ell [2\fl],j,m_j)$ according to theorem 4.9 which leads to the bijection:
  \[
  G^{(\2nfl)}((
F_{\o\omega _\oplus}\otimes \cit)\times
(F_{\omega _\oplus}\otimes \cit))
  \overset{\sim}{\To} \EIS\RL(2n_\ell [2\fl],j,m_j)\;.\]
  \Ee
  

\item Finally, from proposition 4.5, it results that:
  \[ 
G^{(2n_\ell [2\fl])}(({F^T_{\o\omega }}\otimes\cit)\times
({F^T_{\omega }}\otimes\cit)))
\approx
{\rm CY}^{2n_\ell }_{\oplus}(Y_R,[2\fl])\times
{\rm CY}^{2n_\ell }_{\oplus}(Y_L,[2\fl])\]
where
\[{\rm CY}^{2n_\ell }_{\oplus}(Y_L,[2\fl])=
\bigoplus_j\bigoplus_{m_j}{\rm CY}^{2n_\ell }_{\oplus}(Y_L,[2\fl])\;.\qedhere\]
\Ei
\end{proof}
\vskip 11pt

\subsection{Definition: Partially reducible shifted representation}

Similarly as in \cite{Pie3}, shifted global reducible correspondences of Langlands can be introduced.

Let us, for example, consider a partition $n=n_1+\dots+n_\ell +\dots+n_s$ of $n$ leading to the shifted partition:
\[ n[\fnl]=n_1[f_1\cdot\ell]+\dots+n_\ell [\fl]+\dots+n_s[f_s\cdot\ell]\;, \]
where it can be assumed that
$f_n=f_1+\dots+f_\ell +\dots+f_s$, in such a way that:
\begin{multline*}
 \qquad \qquad \Rep(\GL_{n[\fnl]=n_1[f_1\cdot\ell]+\dots+n_s[f_s\cdot\ell]}
((F_{\o\omega }\otimes\cit)\times
(F_{\omega }\otimes\cit)))\\
= 
 \mathop{\boxplus}\limits^{n_s}_{n_\ell =n_1}
\Irr\Rep(\GL_{n_\ell [\fl]}
((F_{\o\omega }\otimes\cit)\times
(F_{\omega }\otimes\cit)))\qquad 
\end{multline*}
constitutes a partially reducible shifted representation of
$ \GL_{n[\fnl]}
((F_{\o\omega }\otimes\cit)\times
(F_{\omega }\otimes\cit))$.
\vskip 11pt

\subsection{Proposition}

{\em
If
\[
G^{(2n[2\fnl])}
((F_{\o\omega }\otimes\cit)\times
(F_{\omega }\otimes\cit))
= 
 \mathop{\boxplus}\limits^{n_s}_{n_\ell =n_1}
G^{(2n_\ell [2\fl])}
((F_{\o\omega }\otimes\cit)\times
(F_{\omega }\otimes\cit))\]
represents the decomposition of the shifted $2n$-dimensional bilinear complete semigroup into irreducible components of dimension $2n_\ell $ shifted in $2\fl$ dimensions, then we have that:
\begin{align*}
G^{(2n[2\fnl])}
((F_{\o\omega }\otimes\cit)\times
(F_{\omega }\otimes\cit))
&= 
\RedRep^{(2n[2\fnl])}_{W_{F\RL}}(W^{ab}_{F^{S_{\cit}}_{\o\omega }}\times W^{ab}_{F^{S_{\cit}}_{\omega }}) \\
&= \bigoplus^{n_s}_{n_\ell =n_1}
\Irr\Rep^{(2n_\ell [2\fl])}_{W_{F\RL}}(W^{ab}_{F^{S_{\cit}}_{\o\omega }}\times W^{ab}_{F^{S_{\cit}}_{\omega }}) \end{align*}
where
$\RedRep^{(2n [2\fl])}_{W_{F\RL}}(W^{ab}_{F^{S_{\cit}}_{\o\omega }}\times
W^{ab}_{F^{S_{\cit}}_{\omega }}) $
denotes the $2n$-dimensional reducible shifted global Weil-Deligne representation of the shifted bilinear global Weil group
$(W^{ab}_{F^{S_{\cit}}_{\o\omega }}\times
W^{ab}_{F^{S_{\cit}}_{\omega }})$.
}
\vskip 11pt

\begin{proof}
According to lemma 5.1, we have that:
\[
\Irr\Rep^{(2n_\ell [2\fl])}_{W_{F\RL}}(W^{ab}_{F^{S_{\cit}}_{\o\omega }}\times W^{ab}_{F^{S_{\cit}}_{\omega }})
=G^{(2n_\ell [2\fl])}
((F_{\o\omega }\otimes\cit)\times
(F_{\omega }\otimes\cit))\]
from which the thesis follows if definition 5.3 is taken into account.\end{proof}
\vskip 11pt

\subsection{Proposition}
{\em
The toroidal compactification of
 $G^{(2n[2\fnl])}
((F_{\o\omega _\oplus}\otimes\cit)\times
(F_{\omega _\oplus}\otimes\cit))$ generates
\[ 
G^{(2n[2\fnl])}
((F^T_{\o\omega_\oplus }\otimes\cit)\times
(F^T_{\omega_\oplus}\otimes\cit))
= \bigoplus^{n_s}_{n_\ell =n_1} G^{2n_\ell [2\fl]}
((F^T_{\o\omega _\oplus}\otimes\cit)\times
(F^T_{\omega_\oplus }\otimes\cit))\]
whose supercuspidal representation is given by:
\[
\Redcusp(\GL_{n [\fnl]}(({F^T_{\o\omega }}\otimes\cit)\times
({F^T_{\omega }}\otimes\cit)))
= \bigoplus^{n_s}_{n_\ell =n_1}
\EIS\RL(2n_\ell [2\fl],j,m_j)\]
where
$\EIS\RL(2n_\ell [2\fl],j,m_j)$ is the $2n_\ell $-dimensional truncated shifted cuspidal biserie.
}
\vskip 11pt

\begin{proof}
This directly results from definition 5.3 and proposition 4.13.\end{proof}
\vskip 11pt

\subsection{Proposition}

{\em
Let
$\o S^{P_{n [\fnl]}}_{\GL_{n[\fnl]}}=
\mathop{\boxplus}\limits^{n_s}_{n_\ell =n_1}
\o S^{P_{n_\ell [\fl]}}_{\GL_{n_\ell [\fl]}}$
be the decomposition of the reducible shifted $2n$-\linebreak dimensional bisemispace
$\o S^{P_{n [\fnl]}}_{\GL_{n[\fnl]}}$ into irreducible components
$\o S^{P_{n_\ell [\fl]}}_{\GL_{n_\ell [\fl]}}$ introduced in section 4.1.

Then, the cohomology of this reducible shifted bisemispace
$\o S^{P_{n[\fnl]}}_{\GL_{n[\fnl]}}$ decomposes according to:
\begin{multline*}
\qquad H^*(\o S^{P_{n[\fnl]}}_{\GL_{n[\fnl]}}
,
\widehat M^{2n}_{T_{R}}[2\fnl]\otimes
\widehat M^{2n }_{T_{L}}[2\fnl])
\\
\begin{aligned}
&= \bigoplus_{n_\ell=n_1 }^{n_s}
H^{2n_\ell -2\fl}
(\o S^{P_{n[\fnl]}}_{\GL_{n[\fnl]}}
, 
\widehat M^{2n_\ell }_{T_{R}}[2\fl]\otimes
\widehat M^{2n_\ell }_{T_{L}}[2\fl])
\\
&\approx
\bigoplus_{n_\ell }
{\rm CY}^{2n_\ell }(Y_R,[2\fl])\times
{\rm CY}^{2n  }(Y_L,[2\fl])
\end{aligned}\qquad
\end{multline*}
where
$(\widehat M^{2n}_{T_{R}}[2\fnl]\otimes
\widehat M^{2n_\ell }_{T_{L}}[2\fnl])$
is the bisemisheaf over the ``partially reducible''\linebreak
$\GL_{n[\fnl]}
((F^T_{\o\omega }\otimes\cit)\times
(F^T_{\omega }\otimes\cit))$-bisemimodule.
}
\vskip 11pt

\subsection{Proposition}

{\em
On the reducible shifted pseudoramified bilinear complete semigroup
\[
G^{2n[2\fnl]}
((F_{\o\omega }\otimes\cit)\times
(F_{\omega }\otimes\cit))=
\mathop{\boxplus}\limits^{n_s}_{n_\ell =n_1}
G^{2n_\ell [2\fl]}
((F_{\o\omega }\otimes\cit)\times
(F_{\omega }\otimes\cit))\;,\]
there exists the {\bbf \em geometric-shifted global bilinear ``reducible'' correspondence of Langlands:}
{ \begin{small}
\[\begin{psmatrix}[colsep=.5cm,rowsep=1cm]
\RedRep^{(2n [2\fnl])}_{W_{F\RL}}
(W^{ab}_{F^{S_{\cit}}_{\o\omega }}\times
W^{ab}_{F^{S_{\cit}}_{\omega }}) 
& \raisebox{3mm}{$\sim$}
&
\Redcusp(\GL_{n[\fnl]}(({F^T_{\o\omega }}\otimes\cit)\times
({F^T_{\omega }}\otimes\cit)))
 \\
G^{(2n [2\fnl])}((F_{\o\omega _{\bigoplus}}\otimes\cit)\times
(F_{\omega _{\bigoplus}}\otimes\cit))
 &  \qquad \raisebox{3mm}{$\sim$}&
\bigoplus_{n_\ell }\EIS\RL(2n_\ell [2\fl],j,m_j)
\\
 && G^{(2n [2\fnl])}
( ( {F^T_{\o\omega }}\otimes\cit )\times
( {F^T_{\omega }} \otimes\cit) )
\\
\noalign{\vspace*{-8mm}}
&& \rotatebox{90}{\Huge{$\simeq$}}\\
\noalign{\vspace*{-8mm}}
&&H^*(\o S^{P_{n[\fnl]}}_{\GL_{n[\fnl]}}
,
\widehat M^{2n}_{T_{R_\oplus}}[2\fnl]\otimes
\widehat M^{2n }_{T_{L_\oplus}}[2\fnl])
\\
\noalign{\vspace*{-8mm}}
&& \rotatebox{90}{\Huge{$\approx$}}\\
\noalign{\vspace*{-8mm}}
&&\bigoplus_{n_\ell }
({\rm CY}^{2n_\ell }_{\oplus}(Y_R,[2\fl])\times
{\rm CY}^{2n_\ell }_{\oplus}(Y_L,[2\fl]))
\ncline[arrows=->,nodesep=5pt]{1,1}{1,3}
\ncline[arrows=->,nodesep=5pt]{2,1}{2,3}
\ncline[doubleline=true,nodesep=5pt]{1,1}{2,1}
\ncline[doubleline=true,nodesep=5pt]{1,3}{2,3}
\ncline[arrows=->,nodesep=5pt]{2,1}{3,3}^{\rotatebox{-30}{$\sim$}}
\ncline[arrows=->,nodesep=5pt]{3,3}{2,3}<{\raisebox{-2mm}{\rotatebox{90}{$\sim$}}}
\end{psmatrix}\]
\end{small}}
\epr}
\vskip 11pt

The objective consists now in establishing shifted global bilinear correspondences on the boundary
$\partial \o S^{P_{n[\fnl]}}_{\GL_{n[\fnl]}}$
of the shifted bisemispace
$\o S^{P_{n[\fnl]}}_{\GL_{n[\fnl]}}$ and to introduce by this way shifted Eisenstein cohomology.
\vskip 11pt

\subsection{Shifted real completions and shifted global Weil groups}

\Bi
\item From sections 1.2 and 1.4, {\bbf direct sums of shifted real completions\/} will be given by:
\[
F^{+,(S_\rit)}_{v_\oplus}
= \bigoplus_{j_\delta }\bigoplus_{m_{j_\delta }}
(F^+_{v_{j_\delta ,m_{j_\delta }}}\otimes\rit)
\rresp{
F^{+,(S_\rit)}_{v_{\o\oplus}}
= \bigoplus_{j_\delta }\bigoplus_{m_{j_\delta }}(F^+_{\o v_{j_\delta ,m_{j_\delta }}}\otimes\rit)}.
\]

 \item According to \cite{Pie3}, there is a one-to-one correspondence between the real shifted completions and their complex equivalents.
 
 \item {\bf The global Weil groups\/}
 $W^{ab}_{F^{+,(S_\rit)}_v}$
 \resp{$W^{ab}_{F^{+,(S_\rit)}_{\o v}}$}, shifted over $\rit$ and referring to pseudoramified extensions characterized by degrees $d=0\mod N$, will be introduced by:
 \[
 W^{ab}_{F^{+,(S_\rit)}_v}
 = \bigoplus_{j_\delta ,m_{j_\delta }}
 \Gal(\dot{\wt F}^{+,(S_\rit)}_{v_{j_\delta ,m_{j_\delta }}}\big/ k)
 \rresp{%
  W^{ab}_{F^{+,(S_\rit)}_{\o v}}
 = \bigoplus_{j_\delta ,m_{j_\delta }}
 \Gal(\dot{\wt F}^{+,(S_\rit)}_{\o v_{j_\delta ,m_{j_\delta }}}\big/ k)}\]
 where
 $\dot{\wt F}^{+,(S_\rit)}_{\o v_{j_\delta ,m_{j_\delta }}}$ denote the shifted pseudoramified extensions with degrees\linebreak $d=0\mod (N)$.
 \Ei
 \vskip 11pt
 
 \subsection[The boundary of the Borel-Serre compactification shifted over $\rit$]{\bbf The boundary of the Borel-Serre compactification shifted over $\rit$}
 
 \Bi
 \item As developed in \cite{Pie3}, the boundary
 $\partial \o Y^{(2n_\ell )}_{S^T\RL}$ of the Borel-Serre compactification
 $\o Y^{(2n_\ell )}_{S^T\RL}= \GL_{n_\ell }(F^T_R\times F^T_L)\big/\GL_{n_\ell }((\zit/N\ \zit)^2)$ of $Y^{(2n_\ell)}_{S\RL}$ (see section 1.12) is given by:
 \begin{align*}
 \partial \o Y^{(2n_\ell )}_{S^T\RL}
 &= \GL_{n_\ell }(F^{+,T}_R\times F^{+,T}_L)\big/ \GL_{n_\ell }((\zit/N\ \zit)^2)\\
 &= \GL_{n_\ell }(F^{+,T}_{\o v}\times F^{+,T}_v)
 \end{align*}
 where:
 \Bi
 \item $F^{+,T}_R$ and $F^{+,T}_L$ are real toroidal compactifications of $\wt F^{+}_R$ and $\wt F^{+}_L$ respectively;
 
 \item $F^{+,T}_v=\{F^{+,T}_{v_1},\dots,
 F^{+,T}_{v_{j_\delta ,m_{j_\delta }}},\dots,
 F^{+,T}_{v_{r_\delta ,m_{r_\delta }}}\}$ is the set of real toroidal completions.
 \Ei
 
{\bbf The boundary $\partial \o Y^{(2n_\ell )}_{S^T\RL}$ shifted over $\rit$ in $(2\fl)$ real dimensions\/} is given by:
\begin{align*}
\partial \o Y^{\2nfl}_{S^T\RL}
&= \GL_{\nfl}((F^{+,T}_R\otimes\rit)\times
(F^{+,T}_L\otimes\rit)) \big/
\GL_{\nfl}((\zit/N\  \zit)^2\otimes\rit^2)\\
&\approx \GL_{\nfl}((F^{+,T}_{\o v}\otimes\rit)\times
(F^{+,T}_v\otimes\rit))
\;.\end{align*}

\item The corresponding {\bbf double coset decomposition of the shifted bilinear semigroup\/}
$\GL_{\nfl}((F^{+,T}_{\o v}\otimes\rit)\times
(F^{+,T}_{v}\otimes\rit))$ can be introduced by:
\begin{multline*}
\partial\o S^{P_{\nfl}}_{\GL_{\nfl}}=
P_{\nfl}((F^{+,T}_{\o v^1}\otimes\rit)\times
(F^{+,T}_{\o v^1}\otimes\rit))\\
\setminus \GL_{\nfl}((F^{+,T}_R\otimes\rit)\times
(F^{+,T}_L\otimes\rit))
\big/ \GL_{\nfl}((\zit/N\ \zit)^2\otimes\rit^2)
\end{multline*}
where $F^{+,T}_{\o v^1}$ and $F^{+,T}_{v^1}$ denote the set of irreducible subcompletions characterized by a degree $N$ (see section 1.1).
\Ei
\vskip 11pt

\subsection[Action of the differential bioperator
$(D^{2\fl}_R\otimes D^{2\fl}_L)$]{\bbf Action of the differential bioperator
$(D^{2\fl}_R\otimes D^{2\fl}_L)$}
 
 The differential bioperator
 $(D^{2\fl}_R\otimes D^{2\fl}_L)$ maps the bisemisheaf
 $(\widehat M^{2n_\ell }_{T_{v_R}}\otimes\widehat M^{2n_\ell }_{T_{v_L}})$ over the
 $\GL_{n_\ell }(F^{+,T}_{\o v}\times (F^{+,T}_{\o v})$-bisemimodule
 $( M^{2n_\ell }_{T_{v_R}}\otimes M^{2n_\ell }_{T_{v_L}})$
 into the corresponding perverse bisemisheaf
 $(\widehat M^{2n_\ell }_{T_{v_R}}[2\fl]\otimes\widehat M^{2n_\ell }_{T_{v_L}}[2\fl])$ over the
 $\GL_{\nfl }(F^{+,T}_{\o v}\times \rit)
 \times (F^{+,T}_{ v}\times \rit))$-bisemimodule
 $( M^{2n_\ell }_{T_{v_R}}[2\fl]\otimes M^{2n_\ell }_{T_{v_L}}[2\fl])$ according to
 \[
 D^{2\fl}_R\otimes D^{2\fl}_L: \qquad
 \widehat M^{2n_\ell }_{T_{v_R}}\otimes\widehat M^{2n_\ell }_{T_{v_L}}\To
 \widehat M^{2n_\ell }_{T_{v_R}}[2\fl]\otimes\widehat M^{2n_\ell }_{T_{v_L}}[2\fl]\;.\]
 
\vskip 11pt
 
 \subsection{Proposition}
 
 {\em
 The  bilinear cohomology of the shifted Shimura bisemivariety
 $\partial \o S^{P_{n[\fnl]}}_{\GL_{n[\fnl]}}$, $n>n_\ell $, is the bilinear shifted Eisenstein cohomology:
\begin{multline*}
\qquad \qquad H^{2n_\ell -2\fl}
(\partial \o S^{P_{n[\fnl]}}_{\GL_{n[\fnl]}})
,
\widehat M^{2n_\ell }_{T_{v_R}}[2\fl]\otimes\widehat M^{2n_\ell }_{T_{v_L}}[2\fl])\\
\simeq
\FRepsp (\GL_{\nfl}((F^{+,T}_{\o v}\otimes \rit)\times
(F^{+,T}_{v}\otimes \rit))\qquad \qquad
\end{multline*}
in such a way that the functional representation space
$\FRepsp (\GL_{\nfl}((F^{+,T}_{\o v}\otimes \rit)\times
(F^{+,T}_{v}\otimes \rit))$ of the shifted complete bilinear semigroup
$\GL_{\nfl}((F^{+,T}_{\o v}\otimes \rit)\times
(F^{+,T}_{v}\otimes \rit)$ is:
\begin{multline*}
\FRepsp (\GL_{\nfl}((F^{+,T}_{\o v}\otimes \rit)\times
(F^{+,T}_{v}\otimes \rit))\\
\begin{aligned}
&\equiv \widehat G^{\2nfl}((F^{+,T}_{\o v}\otimes \rit)\times
(F^{+,T}_{v}\otimes \rit))\\
&= \{
(\phi (g^{n_\ell }_{T_R}([\fl],(j_\delta ,m_{j_\delta }))\times
\phi (g^{n_\ell }_{T_L}([\fl],(j_\delta ,m_{j_\delta })))\}_{j_\delta ,m_{j_\delta }}\\
&\simeq \{
({\rm CY}^{n_\ell }(\partial\o Y_R,[\fl]),(j_\delta ,m_{j_\delta })\times
{\rm CY}^{n_\ell }(\partial\o Y_L,[\fl]),(j_\delta ,m_{j_\delta })\}_{j_\delta ,m_{j_\delta }}
\end{aligned}
\end{multline*}
where:
\Bi
\item $\widehat G^{\2nfl}((F^{+,T}_{\o v}\otimes \rit)\times
(F^{+,T}_{v}\otimes \rit))$ is the bisemisheaf over
$G^{\2nfl}((F^{+,T}_{\o v}\otimes \rit)\times
(F^{+,T}_{v}\otimes \rit))$
 which decomposes according to the set of products, right by left,
 $g^{n_\ell }_{T\RL}([\fl],(j_\delta ,m_{j_\delta }))$ of conjugacy class representatives of real dimension $n_\ell$.
 
 \item the functions 
 $\phi (g^{n_\ell }_{T_L}([\fl],(j_\delta ,m_{j_\delta }))$ on these conjugacy class representatives are in one-to-one correspondence with the equivalent class representative
 ${\rm CY}^{n_\ell }(\partial\o Y_L,[\fl]),(j_\delta ,m_{j_\delta }))$ of cycles of codimension $n_\ell $ shifted in $\fl$ real dimension.
 \Ei
 }
 \vskip 11pt
 
 \begin{proof}
 This is an adaptation to the real case of propositions 4.5 and 4.6.\end{proof}
 \vskip 11pt
 
 \subsection{Shifted global elliptic semimodules}
 
 \Bi
 \item Every \lr function on the conjugacy class representative
 $g^{n_\ell }_{T_L}([\fl],(j_\delta ,m_{j_\delta }))$
 \resp{$g^{n_\ell }_{T_R}([\fl],(j_\delta ,m_{j_\delta }))$} of
 $\GL_{\nfl}((F^{+,T}_{\o v}\otimes \rit)\times
(F^{+,T}_{v}\otimes \rit))$ is a $n_\ell $-dimensional real semitorus
$T^{n_\ell }_{L}([\fl],(j_\delta ,m_{j_\delta }))$
 \resp{$T^{n_\ell }_{R}([\fl],(j_\delta ,m_{j_\delta }))$} shifted in $\fl$ real dimensions, localized in the upper \resp{lower} half space and having the following analytic development:
 \begin{align*}
 T^{n_\ell }_{L}([\fl],(j_\delta ,m_{j_\delta }))
 &\simeq E_{\fl}(n_\ell ,j_\delta ,m_{j_\delta })\
 \lambda ^{\half}(n_\ell ,j_\delta ,m_{j_\delta })\ e^{2\pi ij_\delta x_{n_\ell }}\\[11pt]
 \rresp{%
  T^{n_\ell }_{R}([\fl],(j_\delta ,m_{j_\delta }))
 &\simeq E_{\fl}(n_\ell ,j_\delta ,m_{j_\delta })\
 \lambda ^{\half}(n_\ell ,j_\delta ,m_{j_\delta })\ e^{-2\pi ij_\delta x_{n_\ell }}}
 \end{align*}
 where:
 \Bi
 \item $\lambda (n_\ell ,j_\delta ,m_{j_\delta })$ is a generalized Hecke  global character obtained from the product of the eigenvalues of the $(j_\delta ,m_{j_\delta })$-th coset representative of the Hecke bioperator (see proposition 4.7);
 
 \item $E_{\fl}(n_\ell ,j_\delta ,m_{j_\delta })$ is the shift in $\fl$ real dimensions of
 $\lambda ^{\half}(n_\ell ,j_\delta ,m_{j_\delta })$;
 
 \item $x_{n_\ell }=\sum\limits^{n_\ell }_{\beta =1}x_\beta \ \overrightarrow e^\beta $ is a vector of $\rit^{n_\ell }$.
 \Ei
 
 \item The analytic representation of the shifted bilinear complete semigroup\linebreak
 $G^{\2nfl}((F^{+,T}_{\o v_\oplus}\otimes \rit)\times
(F^{+,T}_{v_\otimes}\otimes \rit))$, which is also a supercuspidal representation of
$\GL_{\nfl}((F^{+,T}_{\o v}\otimes \rit)\times
(F^{+,T}_{v}\otimes \rit))$, is obtained by summing over $j_\delta $ and $m_{j_\delta }$ the analytic representations of the conjugacy class representatives
$g^{n_\ell }_{T\RL}([\fl],(j_\delta ,m_{j_\delta }))$ giving rise to the product, right by left, of the shifted global elliptic semimodules \cite{Drin}:
\begin{multline*}
\ELLIP_R(\2nfl,j_\delta ,m_{j_\delta })\otimes
\ELLIP_L(\2nfl,j_\delta ,m_{j_\delta })\\
= \L[
\bigoplus^r_{j_\delta =1}
\bigoplus_{m_{j_\delta }}
E_{\fl}(n_\ell ,j_\delta ,m_{j_\delta })\
 \lambda ^{\half}(n_\ell ,j_\delta ,m_{j_\delta })\ e^{-2\pi ij_\delta x_{n_\ell }}\R]\\
 \otimes
 \L[
\bigoplus^r_{j_\delta =1}
\bigoplus_{m_{j_\delta }}
E_{\fl}(n_\ell ,j_\delta ,m_{j_\delta })\
 \lambda ^{\half}(n_\ell ,j_\delta ,m_{j_\delta })\ e^{2\pi ij_\delta x_{n_\ell }}\R]\;.\end{multline*}
 \Ei

\subsection{Proposition}

{\em
The shifted bilinear Eisenstein cohomology has the following analytic development:
\begin{multline*}
\qquad \qquad H^{2n_\ell -2\fl}
(\partial \o S^{P_{n[\fnl]}}_{\GL_{n[\fnl]}})
,
\widehat M^{2n_\ell }_{T_{v_{R_\oplus}}}[2\fl]\otimes\widehat M^{2n_\ell }_{T_{v_{L_\oplus}}}[2\fl])\\
\simeq \ELLIP_R(\2nfl,j_\delta ,m_{j_\delta })\otimes
\ELLIP_L(\2nfl,j_\delta ,m_{j_\delta })\qquad \qquad
\end{multline*}
and gives rise to the eigen(bi)value equation:
\begin{multline*}
\qquad (D^{2\fl}_R\otimes D^{2\fl}_L)
(\ELLIP_R(2n_\ell ,j_\delta^{\rm up}=j_\delta  ,m_{j_\delta })\otimes
\ELLIP_L(2n_\ell ,j_\delta^{\rm up}=j_\delta ,m_{j_\delta }))\\
= E^2_{\fl}(n_\ell ,j_\delta ,m_{j_\delta })_{\rm eig}
(\ELLIP_R(2n_\ell ,j_\delta^{\rm up}=j_\delta ,m_{j_\delta })\otimes
\ELLIP_L(2n_\ell ,j_\delta^{\rm up}=j_\delta ,m_{j_\delta }))\;.
\end{multline*}
}
\vskip 11pt

\begin{proof}
This is an adaptation to the real case of proposition 4.9.
\end{proof}
\vskip 11pt

\subsection{Proposition}

{\em Taking into account that the irreducible $2n_\ell $-dimensional shifted global Weil-Deligne representation 
$\Irr\Rep^{(\2nfl)}_{W_{F^+\RL}}(W^{ab}_{F^{+,(S_{\rit)}}_{\o v }}\times W^{ab}_{F^{+,(S_{\rit)}}_{v }})$ of the shifted bilinear global Weil group
$(W^{ab}_{F^{+,(S_{\rit)}}_{\o v }}\times W^{ab}_{F^{+,(S_{\rit)}}_{v }})$ is given by the shifted bilinear complete semigroup
$G^{\2nfl}((F^{+,T}_{\o v_\oplus}\otimes \rit)\times
(F^{+,T}_{v_\oplus}\otimes \rit))$, we have on the  shifted Shimura bisemivariety
$\partial \o S^{P_{\nfl}}_{\GL_{\nfl}}$ the following {\bbf \em geometric-shifted global bilinear correspondence of Langlands} {\em \cite{Lan}}:
{ \begin{small}
\[\begin{psmatrix}[colsep=.2cm,rowsep=1cm]
\Irr\Rep^{(\2nfl)}_{W_{F^+\RL}}
(W^{ab}_{F^{+,(S_{\rit})}_{\o v }}\times
W^{ab}_{F^{+,(S_{\rit})}_{v }}) 
& \raisebox{3mm}{$\sim$}
&
\Irr\ELLIP(\GL_{\nfl}(({F^{+,T}_{\o v }}\otimes\rit)\times
({F^{+,T}_{v }}\otimes\rit)))
 \\
G^{(\2nfl)}((F^+_{\o v _{\bigoplus}}\otimes\rit)\times
(F^+_{v _{\bigoplus}})\otimes\rit)
 &  &
\ELLIP\RL(\2nfl,j_\delta ,m_{j_\delta })
\\
 && G^{(\2nfl)}
( ( {F^{+,T}_{\o v }}\otimes\rit )\times
( {F^{+,T}_{v }}\otimes\rit ) )
\\
\noalign{\vspace*{-8mm}}
&& \rotatebox{90}{\Huge{$\simeq$}}\\
\noalign{\vspace*{-8mm}}
&&H^{2n_\ell -2\fl}(\partial\o S^{P_{n[\fnl]}}_{\GL_{n[\fnl]}}
, \widehat M^{2n_\ell }_{T_{v_R}}[2\fl]\otimes
\widehat M^{2n_\ell  }_{T_{v_L}}[2\fl])
\\
\noalign{\vspace*{-8mm}}
&& \rotatebox{90}{\Huge{$\approx$}}\\
\noalign{\vspace*{-8mm}}
&&
{\rm CY}^{2n_\ell }(\partial\o Y_R,[2\fl])\times
{\rm CY}^{2n_\ell }(\o Y_L,[2\fl])
\ncline[arrows=->,nodesep=5pt]{1,1}{1,3}
\ncline[doubleline=true,nodesep=5pt]{1,1}{2,1}
\ncline[doubleline=true,nodesep=5pt]{1,3}{2,3}
\ncline[arrows=->,nodesep=5pt]{2,1}{3,3}^{\rotatebox{-30}{$\sim$}}
\ncline[arrows=->,nodesep=5pt]{3,3}{2,3}<{\raisebox{-2mm}{\rotatebox{90}{$\sim$}}}
\end{psmatrix}\]
\end{small}}
where $\Irr\ELLIP(\GL_{\nfl}(({F^{+,T}_{\o v }}\otimes\rit)\times
({F^{+,T}_{v }}\otimes\rit)))$ is the shifted irreducible elliptic representation of
$\GL_{\nfl}(({F^{+,T}_{\o v }}\otimes\rit)\times
({F^{+,T}_{v }}\otimes\rit))$ given by the $2n_\ell $-dimensional solvable elliptic bisemimodule
$\ELLIP\RL(\2nfl,j_\delta ,m_{j_\delta })$ shifted in 
$(2\fl)$ dimensions.
}
\vskip 11pt

\subsection[The complex and real $2n_\ell $-dimensional irreducible geometric-shifted global bilinear correspondences of Langlands]{\bbf The complex and real  $2n_\ell $-dimensional irreducible\linebreak geometric-shifted global bilinear correspondences of\linebreak Langlands}

can be summarized in the following diagram:
{ \begin{small}
\[\begin{psmatrix}[colsep=.5cm,rowsep=1cm]
\Irr\Rep^{(\2nfl)}_{W_{F\RL}}
(W^{ab}_{F^{S_{\cit}}_{\o \omega  }}\times
W^{ab}_{F^{S_{\cit}}_{\omega  }}) 
& &
\Irr\cusp(\GL_{\nfl}(({F^{T}_{\o \omega  }}\otimes\cit)\times
({F^{T}_{\omega  }}\otimes\cit)))
 \\
\Irr\Rep^{(\2nfl)}_{W_{F^+\RL}}
(W^{ab}_{F^{+,(S_{\rit})}_{\o v  }}\times
W^{ab}_{F^{(+,(S_{\rit})}_{v  }}) 
& &
\Irr\ELLIP(\GL_{\nfl}(({F^{+,T}_{\o v  }}\otimes\rit)\times
({F^{+,T}_{v  }}\otimes\rit)))
\ncline[arrows=->,nodesep=5pt]{1,1}{1,3}
\ncline[arrows=->,nodesep=5pt]{2,1}{2,3}
\ncline[arrows=->,nodesep=5pt]{1,1}{2,1}
\ncline[arrows=->,nodesep=5pt]{1,3}{2,3}
\end{psmatrix}\]
\end{small}}
\vskip 11pt

\subsection{Definition}

The partially reducible shifted representation of
$\GL_{n[\fnl]}(F^+_{\o v}\otimes\rit)\times
(F^+_{v}\otimes\rit))$ can be introduced, as in definition 5.3, on the basis of the shifted partition (real dimensions):
$ n\cdot[\fnl]=
n_1[f_1\cdot\ell ]+\dots+
n_\ell [f_\ell \cdot\ell ]+\dots+
n_s [f_s \cdot\ell ]$ by:
\begin{multline*}
 \qquad \qquad \Rep(\GL_{n[\fnl]=n_1[f_1\cdot\ell]+\dots+n_s[f_s\cdot\ell]}
((F^+_{\o v }\otimes\rit)\times
(F^+_{v }\otimes\rit)))\\
= 
 \mathop{\boxplus}\limits^{n_s}_{n_\ell =n_1}
\Irr\Rep(\GL_{n_\ell [\fl]}
((F^+_{\o v }\otimes\rit)\times
(F^+_{v }\otimes\rit)))\;.\qquad \qquad
\end{multline*}
\vskip 11pt

\subsection{Proposition}

{\em
If
\[ 
G^{(2n[2\fnl])}
((F^+_{\o v }\otimes\rit)\times
(F^+_{v }\otimes\rit)))
= 
 \mathop{\boxplus}\limits^{n_s}_{n_\ell =n_1}
G^{(2n_\ell [2\fl])}
((F^+_{\o v }\otimes\rit)\times
(F^+_{v }\otimes\rit)))\]
represents the decomposition of the shifted $2n$-dimensional real bilinear complete semigroup into irreducible components of dimension $2n_\ell $  shifted in $2\fl$ dimensions, then the $2n$-dimensional reducible shifted global Weil 
representation of the shifted bilinear global Weil
group
$(W^{ab}_{F^{+,(S_{\rit})}_{\o v  }}\times
W^{ab}_{F^{(+,(S_{\rit})}_{v  }})$ is given by:
\begin{align*}
\RedRep^{(2n [2\fnl])}_{W_{F^+\RL}}
(W^{ab}_{F^{+,(S_{\rit})}_{\o v  }}\times
W^{ab}_{F^{(+,(S_{\rit})}_{v  }})
&= 
G^{(2n[2\fnl])}
((F^+_{\o v _\oplus}\otimes\rit)\times
(F^+_{v_\oplus }\otimes\rit)))\\
&= \bigoplus^{n_s}_{n_\ell =n_1}
\Irr\Rep^{(\2nfl)}_{W_{F^+\RL}}
(W^{ab}_{F^{+,(S_{\rit})}_{\o v  }}\times
W^{ab}_{F^{(+,(S_{\rit})}_{v  }})\;.\end{align*}
}
\vskip 11pt

\subsection{Proposition}

{\em
The toroidal compactification of
$G^{(2n[2\fnl])}
((F^+_{\o v _\oplus}\otimes\rit)\times
(F^+_{v_\oplus }\otimes\rit)))$ generates by decomposition:
\[
G^{(2n[2\fnl])}
((F^{+,T}_{\o v _\oplus}\otimes\rit)\times
(F^{+,T}_{v_\oplus }\otimes\rit)))
=\bigoplus^{n_s}_{n_\ell =n_1}
G^{(2n_\ell [2\fl])}
((F^{+,T}_{\o v _\oplus}\otimes\rit)\times
(F^{+,T}_{v_\oplus }\otimes\rit)))\]
whose elliptic representation is given by:
\[
\RedELLIP(\GL_{n[\fnl]}(({F^{+,T}_{\o v }}\otimes\rit)\times
({F^{+,T}_{v }}\otimes\rit)))
= \bigoplus^{n_s}_{n_\ell =n_1}
\ELLIP\RL(2n_\ell [2\fl],j_\delta ,m_{j_\delta })\]
where $\ELLIP\RL(2n_\ell [2\fl],j_\delta ,m_{j_\delta })$ is the product, right by left, of $2n_\ell $-dimensional shifted global elliptic semimodules as introduced in section 5.12.
}
\vskip 11pt

\subsection{Proposition}

{\em Let
\[
\partial \o S^{P_{n[\fnl]}}_{\GL_{n[\fnl]}}= \mathop{\boxplus}\limits^{n_s}_{n_\ell =n_1}
\partial\o S^{P_{n_\ell [\fl]}}_{\GL_{n_\ell [\fl]}}
\]
be the decomposition of the boundary
$ \partial \o S^{P_{n[\fnl]}}_{\GL_{n[\fnl]}}$
of the reducible shifted bisemispace\linebreak
$ \o S^{P_{n[\fnl]}}_{\GL_{n[\fnl]}}$ into irreducible components
$\partial \o S^{P_{n_\ell [\fl]}}_{\GL_{n_\ell [\fl]}}$.  Then, the Eisenstein cohomology of this reducible shifted bisemispace
$\partial \o S^{P_{n[\fnl]}}_{\GL_{n[\fnl]}}$ decomposes following:
\begin{multline*}
\qquad \qquad H^*(\partial\o S^{P_{n[\fnl]}}_{\GL_{n[\fnl]}}
,
\widehat M^{2n}_{T_{R_\oplus}}[2\fnl]\otimes
\widehat M^{2n }_{T_{L_\oplus}}[2\fnl])\\
\begin{aligned}
&=\bigoplus^{n_s}_{n_\ell =n_1}
H^{2n_\ell -2\fl}(\partial\o S^{P_{n[\fnl]}}_{\GL_{n[\fnl]}}
, \widehat M^{2n_\ell }_{T_{v_R}}[2\fl]\otimes
\widehat M^{2n_\ell  }_{T_{v_L}}[2\fl])\\
&\simeq \bigoplus_{n_\ell}
{\rm CY}^{2n_\ell }(\partial\o Y_R,[2\fl])\times
{\rm CY}^{2n_\ell }(\o Y_L,[2\fl])
\end{aligned}\end{multline*}
where $( \widehat M^{2n }_{T_{v_R}}[2\fnl]\otimes
\widehat M^{2n  }_{T_{v_L}}[2\fnl])$ is the bisemisheaf
over the partially reducible\linebreak
$\GL_{n[\fnl]}((F^{+,T}_{\o v }\otimes\rit)\times
(F^{+,T}_{v }\otimes\rit))$-bisemimodule.
}\vskip 11pt

\subsection{Proposition}

{\em
On the reducible shifted pseudoramified bilinear complete semigroup
\[
G^{(2n[2\fnl])}((F^{+}_{\o v }\otimes\rit)\times
(F^{+}_{v }\otimes\rit))
= \mathop{\boxplus}\limits^{n_s}_{n_\ell =n_1}
G^{(2n_\ell [2\fl])}((F^{+}_{\o v }\otimes\rit)\times
(F^{+}_{v }\otimes\rit))\]
there exists the {\bbf \em geometric-shifted global bilinear ``reducible'' correspondence of Langlands:}
{ \begin{small}
\[\begin{psmatrix}[colsep=.25cm,rowsep=1cm]
\RedRep^{(2n[2\fnl])}_{W_{F^+\RL}}
(W^{ab}_{F^{+,(S_{\rit})}_{\o v }}\times
W^{ab}_{F^{+,(S_{\rit})}_{v }}) 
& \raisebox{3mm}{$\sim$}
&
\RedELLIP(\GL_{n[\fnl]}(({F^{+,T}_{\o v }}\otimes\rit)\times
({F^{+,T}_{v }}\otimes\rit)))
 \\
G^{(2n[2\fnl])}((F^+_{\o v _{\oplus}}\otimes\rit)\times
(F^+_{v _{\oplus}})\otimes\rit)
 &  &
\bigoplus_{n_\ell }\ELLIP\RL(\2nfl,j_\delta ,m_{j_\delta })
\\
 && G^{(2n[2\fnl])}
( ( {F^{+,T}_{\o v }}\otimes\rit )\times
( {F^{+,T}_{v }}\otimes\rit ) )
\\
&&H^{*}(\partial\o S^{P_{n[\fnl]}}_{\GL_{n[\fnl]}}
, \widehat M^{2n }_{T_{v_{R_\oplus}}}[2\fnl]\otimes
\widehat M^{2n   }_{T_{v_{L_\oplus}}}[2\fnl])
\ncline[arrows=->,nodesep=5pt]{1,1}{1,3}
\ncline[doubleline=true,nodesep=5pt]{3,3}{4,3}
\ncline[doubleline=true,nodesep=5pt]{1,1}{2,1}
\ncline[doubleline=true,nodesep=5pt]{1,3}{2,3}
\ncline[arrows=->,nodesep=5pt]{2,1}{3,3}^{\rotatebox{-30}{$\sim$}}
\ncline[arrows=->,nodesep=5pt]{3,3}{2,3}<{\raisebox{-2mm}{\rotatebox{90}{$\sim$}}}
\end{psmatrix}\]
\end{small}}
}